%% file: self.tex
\documentclass[11pt,makeidx]{article}
\input amssym.def
\input amssym.tex
\usepackage{epsfig,psfig}
\newtheorem{theorem}{Theorem}
\newtheorem{lemma}{Lemma}

\newtheorem{coro}{Corollary}
\newtheorem{defn}{Definition}

\def\biblitem#1{\bibitem{#1}}
\def\binom#1#2{{#1}\choose{#2}}

\def\slfrac#1#2{\hbox{\kern.1em %
 \raise.5ex\hbox{\the\scriptfont0 #1}\kern-.11em %
 /\kern-.15em\lower.25ex\hbox{\the\scriptfont0 #2}}}

\newcommand{\eqn}[1]{(\ref{#1})}
\newcommand{\hsp}{\hspace*{\parindent}}

\newcommand{\eeq}{\end{equation}}
\newcommand{\beql}[1]{\begin{equation}\label{#1}}
\newcommand{\bsq}{{\vrule height .9ex width .8ex depth -.1ex }}

\newcommand{\Eone}{\ensuremath{2}}
\newcommand{\EoneI}{\ensuremath{2_{\rm I}}}
\newcommand{\EoneII}{\ensuremath{2_{\rm II}}}
\newcommand{\Etwo}{\ensuremath{3}}
\newcommand{\Ethree}{\ensuremath{4^{\rm H}}}
\newcommand{\Efour}{\ensuremath{4^{\rm E}}}
\newcommand{\Efive}{\ensuremath{4^{\rm H+}}}
\newcommand{\EfiveI}{\ensuremath{4^{\rm H+}_{\rm I}}}
\newcommand{\EfiveII}{\ensuremath{4^{\rm H+}_{\rm II}}}
\newcommand{\Eff}{\ensuremath{4^{\rm E+}}}

\newcommand{\Esix}{\ensuremath{q^{\rm H}}}
\newcommand{\Eseven}{\ensuremath{q^{\rm E}}}
\newcommand{\Eeight}{\ensuremath{4^{\ZZ}}}
\newcommand{\EeightI}{\ensuremath{4^{\ZZ}_{\rm I}}}
\newcommand{\EeightII}{\ensuremath{4^{\ZZ}_{\rm II}}}
\newcommand{\Enine}{\ensuremath{m^{\ZZ}}}
\newcommand{\EnineI}{\ensuremath{m^{\ZZ}_{\rm I}}}
\newcommand{\EnineII}{\ensuremath{m^{\ZZ}_{\rm II}}}

\newcommand{\hA}{\hat{A}}

\newcommand{\sH}{{\cal H}}

\newcommand{\sE}{{\cal E}}
\newcommand{\sA}{{\cal A}}

\newcommand{\ZZ}{{\Bbb Z}}
\newcommand{\RR}{{\Bbb R}}
\newcommand{\FF}{{\Bbb F}}
\newcommand{\QQ}{{\Bbb Q}}
\newcommand{\CC}{{\Bbb C}}

\newcommand{\bx}{{\bf x}}

\newcommand{\af}{\alpha}
\newcommand{\la}{\lambda}

\newcommand{\om}{\omega}
\newcommand{\oom}{\overline{\omega}}
\newcommand{\oga}{\overline{\gamma}}
\newcommand{\oaf}{\overline{\alpha}}
\newcommand{\obe}{\overline{\beta}}
\newcommand{\of}{\overline{f}}

\newcommand{\ox}{\overline{x}}
\newcommand{\oy}{\overline{y}}
\newcommand{\ou}{\overline{u}}

\newcommand{\La}{\Lambda}

\newcommand{\ep}{\epsilon}
\newcommand{\Om}{\Omega}

\newcommand{\dd}{\ldots}

\newcommand{\sG}{{\cal G}}
\newcommand{\sJ}{{\cal J}}

\newcommand{\sD}{{\cal D}}
\newcommand{\sB}{{\cal B}}

\newcommand{\sK}{{\cal K}}

\makeatletter
\def\@sect#1#2#3#4#5#6[#7]#8{\ifnum #2>\c@secnumdepth
     \def\@svsec{}\else
     \refstepcounter{#1}\edef\@svsec{\csname the#1\endcsname.\hskip .75em }\fi
     \@tempskipa #5\relax
      \ifdim \@tempskipa>\z@
        \begingroup #6\relax
          \@hangfrom{\hskip #3\relax\@svsec}{\interlinepenalty \@M #8\par}%
        \endgroup
       \csname #1mark\endcsname{#7}\addcontentsline
         {toc}{#1}{\ifnum #2>\c@secnumdepth \else
                      \protect\numberline{\csname the#1\endcsname}\fi
                    #7}\else
        \def\@svsechd{#6\hskip #3\@svsec #8\csname #1mark\endcsname
                      {#7}\addcontentsline
                           {toc}{#1}{\ifnum #2>\c@secnumdepth \else
                             \protect\numberline{\csname the#1\endcsname}\fi
                       #7}}\fi
     \@xsect{#5}}
\def\@begintheorem#1#2{\it \trivlist \item[\hskip \labelsep{\bf #1\ #2.}]}

\def\plain{plain}\ifx\fmtname\plain\csname fi\endcsname
     
     \input docstrip
     \preamble

     Do not distribute the stripped version of this file.
     The checksum in the header refers to the documented version.

     \endpreamble
     \generateFile{here.sty}{t}{\from{here.doc}{}}
     \endinput
\fi
\ifcat a\noexpand @\let\next\relax\else\def\next{%
    \documentstyle[here,doc]{article}\MakePercentIgnore}\fi\next
\ifx\@Hxfloat\@Hundef\else\expandafter\endinput\fi
\let\@Hxfloat\@xfloat
\def\@xfloat#1[{\@ifnextchar{H}{\@HHfloat{#1}[}{\@Hxfloat{#1}[}}
\def\@HHfloat#1[H]{%
\expandafter\let\csname end#1\endcsname\end@Hfloat
\vskip\intextsep\vbox\bgroup\def\@captype{#1}\parindent\z@
\ignorespaces}
\def\end@Hfloat{\egroup\vskip \intextsep}

\makeatother

\begin{document}
\begin{center}
{\Large {\bf Self-Dual Codes}} \\
\vspace{1.5\baselineskip}
{\em E. M. Rains} and {\em N. J. A. Sloane} \\
\vspace*{1\baselineskip}
Information Sciences Research, AT\&T Labs-Research \\
180 Park Avenue, Florham Park, NJ 07932-0971 \\
\vspace{1.5\baselineskip}
May 19 1998 \\
\vspace{1.5\baselineskip}
{\bf ABSTRACT}
\vspace{.5\baselineskip}
\end{center}

A survey of self-dual codes, written for the {\em Handbook of Coding Theory}.

\vspace*{+.1in}
Self-dual codes are important because many of the best codes known are of this type and they have a rich
mathematical theory.
Topics covered in this chapter include codes over $\FF_2$, $\FF_3$, $\FF_4$,
$\FF_q$, $\ZZ_4$, $\ZZ_m$, shadow codes, weight enumerators, Gleason-Pierce theorem,
invariant theory, Gleason theorems, bounds, mass formulae, enumeration,
extremal codes,
open problems.
There is a comprehensive bibliography.
\clearpage
\thispagestyle{empty}
\setlength{\baselineskip}{1.5\baselineskip}

\section{Self-dual codes over rings and fields}\label{SD}
\hsp
\subsection{Inner products}\label{SD1}
\hsp
There are several different kinds of self-dual codes.
Let $\FF$ be a finite set called the {\em alphabet} (e.g.
$\FF= \{0,1\}$ for binary codes).
A {\em code $C$ over $\FF$} of {\em length} $n$ is any subset of $\FF^{\,n}$.
If $\FF$ has the structure of an additive group then $C$ is {\em additive} if it is an additive subgroup of $\FF^{\,n}$.
\index{code!linear}\index{code!additive}
If $\FF$ has a ring structure then $C$ is {\em linear over $\FF$} if it is additive and also closed under multiplication by elements of $\FF$.
(We will always assume that multiplication in $\FF$ is commutative.)

In order to define dual codes we must equip $\FF$ with an {\em inner product}\index{inner product}
(cf. \cite{Lang65}, \cite{MiHu}).
We denote this by $(~,~)$ and require that it satisfy the following conditions:
$$
\begin{array}{l}
(x+y,z) = (x,z) + (y,z) ~, \\
(x,y+z) = (x,y) + (x,z) ~, \\
{\rm if} ~~ (x,y) =0 ~~\mbox{for all $x$ then $y=0$}~, \\
{\rm if} ~~ (x,y) =0 ~~\mbox{for all $y$ then $x=0$}~.
\end{array}
$$
To define the dual of a linear code we impose the further condition that
$\FF$ has a {\em conjugacy}\index{conjugation} operation, or ``involutory anti-automorphism''
(which may be the identity), denoted by a bar, which satisfies
$$\overline{\overline{x}} = x , ~
\overline{x+y} = \overline{x} + \overline{y} , ~
\overline{xy} = \overline{x}\, \overline{y} ~.
$$
The inner product must then satisfy
$$(x,y) = \overline{(y,x)} , ~ (ax,y) = (x, \overline{a} y) ~.$$

The inner product of vectors $x= (x_1, \dd , x_n )$, $y= (y_1, \dd , y_n )$ in $\FF^{\,n}$ is defined by
$$(x,y) = \sum_{i=1}^n (x_i , y_i) ~.$$
\subsection{Families of self-dual codes}\label{SD2}
\hsp
Families (\Eone) through (\Enine) include the most important families of codes we will
consider in this chapter.

\noindent
(\Eone) Binary linear codes:\index{binary linear code}\index{\Eone\ (family of binary self-dual codes)}
$\FF = \FF_2 = \{ 0,1 \}$, with inner product
$(x,y) = xy$, $C = \mbox{subspace of $\FF_2^{\,n}$}.$

\noindent
(\Etwo) Ternary linear codes:\index{ternary linear code}\index{\Etwo\ (family of ternary self-dual codes)}
$\FF = \FF_3 = \{ 0,1,2 \}$,
$(x,y) = xy$, $C = \mbox{subspace of $\FF_3^{\,n}$}.$

\noindent
(\Ethree)\index{quaternary linear code}\index{\Ethree\ (Hermitian self-dual over $\FF_4$)} Quaternary linear codes:
$\FF = \FF_4 = \{ 0,1, \om , \om^2 \}$,
where $\om^2 + \om +1 =0$, $\om^3 =1$, $\ox = x^2$ for $x \in \FF_4$,
with the Hermitian inner product\index{Hermitian inner product}\index{inner product!Hermitian}
$(x,y) = x \oy$, $C= \mbox{subspace of $\FF_4^{\,n}$}.$
Note that for $x,y \in \FF_4$, $(x+y)^2 = x^2 + y^2$, $x^4 = x$.

\noindent
(\Efour)\index{quaternary linear code}\index{\Efour\ (Euclidean self-dual over $\FF_4$)}
Quaternary linear codes: $\FF = \FF_4$,
but with
the Euclidean inner product\index{Euclidean inner product}\index{inner product!Euclidean} $(x,y) = xy$.

\noindent
(\Efive)\index{quaternary additive code}\index{inner product!trace}\index{\Efive\ (additive trace-self-dual over $\FF_4$)}
Quaternary additive codes: $\FF = \FF_4$,
with $(x,y) = xy^2 + x^2 y = {\rm trace} (x \oy )$ (the trace from $\FF_4$ to $\FF_2$);
$C = \mbox{additive subgroup of $\FF_4^{\,n}$}$.

For completeness we should also mention family \Eff, quaternary additive codes
with the Euclidean trace inner product:
\index{inner product!Euclidean trace}
$\FF= \FF_4$, with $(x,y) = xy + (xy)^2 = {\rm trace} (xy)$ (the trace
from $\FF_4$ to $\FF_2$);
$C =$ additive subgroup of $\FF_4^{\,n}$.
However, the map
$$x= \omega x_1 + \oom x_2 \in \FF_4^{\,n} \leftrightarrow
x_1 x_2 \in \FF_2^{\,2n}$$
shows that these codes are equivalent to binary codes from family
\Eone~with a particular pairing of the coordinates.  Since we don't know
any interesting examples of this family other than linear codes, we shall
say no more about them.

\noindent
(\Esix)\index{\Esix\ (Hermitian self-dual over $\FF_q$)}
Linear codes over $\FF_q$\index{linear codes over $\FF_q$} (or $q$-ary linear codes), where $q$ is an even
power of an arbitrary prime $p$, with $\ox = x^{\sqrt{q}}$ for $x \in \FF_q$,
$(x,y) = x \oy$, $C= \mbox{subspace of $\FF_q^{\,n}$}$.
Note that for $x,y \in \FF_q$, $(x+y)^{\sqrt{q}} = x^{\sqrt{q}} + y^{\sqrt{q}}$,
$x^q =x$.

\noindent
(\Eseven)\index{\Eseven\ (Euclidean self-dual over $\FF_q$)} Linear codes over $\FF_q$,
but with $(x,y) = xy$.
If $q$ is a square, family \Esix\ is generally preferred to \Eseven.

\noindent
(\Eeight)\index{\Eeight\ (self-dual over $\ZZ_4$)} $\ZZ_4$-linear codes:\index{$\ZZ_4$-linear code}
$\FF = \ZZ_4 = \{0,1,2,3\}$,
with $(x,y) = xy$ $(\bmod ~4)$, $C=$ linear subspace\footnote{Strictly speaking, a $\ZZ_4$-submodule.} of $\ZZ_4^n$.

\noindent
(\Enine)\index{\Enine\ (self-dual over $\ZZ_m$)}\index{$\ZZ_m$-linear code} $\FF= \ZZ_m = \ZZ/m \ZZ$, where $m$ is an integer $\ge 2$, with
$(x,y) = xy$ $(\bmod~m)$,
$C =$ linear subspace\footnote{Strictly speaking, a $\ZZ_m$-submodule.} of $\ZZ_m^n$.

Note
that for the families \Eone, \Etwo, \Eeight, \Enine, an additive code is automatically linear.

The following families are less important for our present purposes:

\noindent
(F1) Linear codes over $\FF_q [u]/ (u^2)$, where $u$ is an indeterminate, with $\ou = -u$, $(x,y) = x \oy$.
(References \cite{Bach97} and \cite{Gab96} consider
such codes, as well as a noncommutative variant.)

\noindent
(F2) Additive codes over $\FF_4$, with $(x,y) = x \oy$.

\noindent
If we relax the requirement that $\FF$ be commutative and finite, we can add:

\noindent
(F3) Linear codes over the $p$-adic integers.
\index{code!p-adic} \index{p-adic code}

\noindent
(F4) Codes over Frobenius rings.\index{codes!over Frobenius rings}
\index{rings!codes over}

\noindent
(F5) Lattices\index{lattice} in $\RR^n$ (see Section~\ref{Latt}).

\subsection{The dual code}\label{SD3}
\hsp
Once we have specified a family of codes by giving $\FF$ and an inner product we can define the {\em dual}\index{dual}\index{code!dual} of a code $C$ to be
$$C^\perp = \{ u \in \FF^{\,n}: (u,v) =0 \quad \mbox{for all} \quad v \in C \} ~.$$
The dual of a binary linear code (family~\Eone) is again a binary linear code.
Similarly, the dual of a code in any of families \Etwo\ through \Enine\
is again a code of the same family.
For family \Efive, the dual of an additive code is additive;
if $C$ is also linear so is $C^\perp$, and then $C^\perp$ coincides with the dual in family \Ethree.
The dual in family \Efour\ is the conjugate of the dual in family \Ethree.

For families \Eone\  through \Enine\ it is easily checked that we have
\beql{Eq1}
|C| \, |C^\perp | = | \FF|^n ~,
\eeq
which implies
\beql{Eq2}
(C^\perp )^\perp = C ~.
\eeq
In general, however, we can say only that
$$C \subseteq (C^\perp )^\perp ~.$$
In particular, \eqn{Eq2} does not necessarily hold for family F2
(consider, for example, the code $\{00,11\}$ which has dual $\{00,11,\om\om,\oom\oom\}$, containing only 4 words).

\subsection{Self-dual codes}\label{SD4}\index{code!self-dual} \index{code!self-orthogonal} \index{code!weakly self-orthogonal}
\hsp
If $C= C^\perp$ then $C$ is said to be {\em self-dual}.\index{self-dual}
If $C \subseteq C^\perp$, $C$ is {\em self-orthogonal}.\index{self-orthogonal}
(In the past, some authors have used ``self-orthogonal'' and ``weakly self-orthogonal''\index{weakly self-orthogonal} for these two concepts.)

In families \Eone\ through \Enine, if $C$ is self-dual then
\beql{Eq3}
|C| = | \FF |^{n/2} ~,
\eeq
and if $| \FF |$ is not a square then $n$ must be even.
In particular, if $C$ is linear over a field, then $n$ is even and $C$
is a subspace of dimension $n/2$.
The only families from \Eone\ through \Enine\ that contain self-dual codes of odd
length are \Efive, \Eeight\ and \Enine\ with $m$ a square.

\paragraph{Remarks about the final three families.}

(F3):\index{code!p-adic} Let $C$ be a code of length $n$ over the $p$-adic integers $\ZZ_{p^\infty}$ (such codes have been studied in \cite{CLP95}, \cite{Me191}).
In general it is not clear how one should define $C^\perp$.
However, if when we reduce $C \bmod~p$ it has the same dimension over
$\FF_p$ as $C$ had over $\FF_{p^\infty}$, then there is a natural way to define the dual so that it satisfies
$$(C^\perp )^\perp = C, \quad \dim~C + \dim~C^\perp =n ~.$$
Namely, let $D = \QQ_{p^\infty} \otimes C$ be the code over
the $p$-adic {\em rationals} $\QQ_{p^\infty}$ generated by $C$.
Since $D$ is a linear code over a field, $D^\perp$ exists and satisfies
$(D^\perp )^\perp =D$, $\dim ~D + \dim D^\perp =n$.
Now set $C^\perp = D^\perp \cap \ZZ_{p^\infty}^n$.

(F4): J. A. Wood\index{Wood, J. A.} (\cite{Wood96}, see also \cite{WaWo96}) has investigated
codes over noncommutative finite rings $\FF$, and has shown that the
two fundamental MacWilliams theorems\index{MacWilliams theorem}\index{theorem!MacWilliams} (Theorem~\ref{thM1} below and Theorems 10.4 and 10.6 of Chapter~1) hold precisely when $\FF$ is a Frobenius ring.\index{Frobenius ring}
At present however no interesting examples of self-dual codes over
noncommutative rings are known.

(F5):\index{lattice!unimodular}
Unimodular lattices are analogues of self-dual codes in $\RR^n$ --- see Section~\ref{Latt}.

\section{Equivalence of codes}\label{AG}
\hsp
\subsection{Equivalent codes}\label{AG1}\index{equivalent}\index{code!equivalent}
\hsp
Codes that differ only in minor ways, such as in the order in which the coordinates are arranged, are said to be {\em equivalent}.
The transformations that we allow in defining equivalence for the above families of codes are as follows (these are precisely the transformations that commute with the process of forming the dual).

\noindent
(\Eone) Permutations of the coordinates.

\noindent
(\Etwo) Monomial transformations of the coordinates (that is, a permutation of the coordinates followed by multiplication of the coordinates by nonzero field elements).

\noindent
(\Ethree) Monomials; global conjugation.

\noindent
(\Efour) Permutations; global conjugation.

\noindent
(\Efive) Monomials; conjugation of individual coordinates.

\noindent
(\Esix) Monomials over the subgroup
$$\{ x \in \FF_q : x \ox = 1\} \cong \FF_q^\ast / \FF_{\sqrt{q}}^\ast ~,$$
where the star denotes the set of nonzero field elements;
global multiplication by elements of $\FF_q^\ast$;
global action of Galois group $Gal ( \FF_q / \FF_p)$

\noindent
(\Eseven) Monomials over $\{ \pm 1 \}$;
global multiplication by units; global action of Galois group.

\noindent
(\Eeight)
Monomials over $\{ \pm 1 \}$.

\noindent
(\Enine)
Monomials over square roots of unity; global multiplication by units of $\ZZ_m$.

\subsection{Automorphism groups}\index{automorphism group}\index{group!automorphism}\index{group!symmetry}\label{AG2}
\hsp
In each case, the subset of such transformations that preserves the code forms the {\em automorphism group} $Aut (C)$ of the code.

Let $G$ denote the full group of all transformations listed.
The order of $G$ in the above cases is:

\noindent
(\Eone) $n!$

\noindent
(\Etwo) $2^n n!$

\noindent
(\Ethree) $2.3^n n!$

\noindent
(\Efour) $2.n!$

\noindent
(\Efive) $6^n n!$

\noindent
(\Esix) $\log_p (q) ( \sqrt{q} -1) ( \sqrt{q} +1)^n n!$

\noindent
(\Eseven) $\log_p (q) \frac{q-1}{2} 2^n n!$

\noindent
(\Eeight) $2^n n!$

\noindent
(\Enine) For $m=5,6,7,8,9$ the orders are
$$\frac{5-1}{2} 2^n n! ,~
2^n n! ,~ \frac{7-1}{2} 2^n n! ,~ 4^n n! ,~ 3.2^n n!
$$
respectively.

The number of codes\index{codes!number of} that are equivalent to a given code $C$ is then
$$\frac{|G|}{|Aut (C)|} ~.$$
In most cases it is possible to determine the total number $T_n$
(say) of distinct self-dual codes of length $n$ in one of our families.
Then
$$T_n = \sum_{{\rm inequivalent} \atop {C}}
\frac{|G|}{|Aut (C)|}$$
where the sum is over all {\em inequivalent} codes.
In other words
\beql{EqAGa}
\sum_{{\rm inequivalent} \atop {C}}
\frac{1}{|Aut (C)|} = \frac{T_n}{|G|} ~.
\eeq
Equation~(\ref{EqAGa}) is called a {\em mass formula}.\index{mass formula}
The appropriate values of $T_n$ are:

\noindent
(\Eone)
\beql{EqAGb}
\prod_{i=1}^{\frac{1}{2} n -1} (2^i +1) \quad (n\equiv 0~( \bmod~2))
\eeq

\noindent
(\EoneII) (weights divisible by 4):
\beql{EqAGc}
2 \prod_{i=1}^{\frac{1}{2} n-2} (2^i +1) \quad (n \equiv 0~( \bmod~8))
\eeq

\noindent
(\Etwo)
\beql{EqAGd}
2 \prod_{i=1}^{\frac{1}{2} n-1} (3^i +1) \quad (n \equiv 0~( \bmod~4))
\eeq

\noindent
(\Ethree)
\beql{EqAGe}
\prod_{i=0}^{\frac{1}{2} n-1}
(2^{2i+1} +1) \quad (n \equiv 0~( \bmod~2))
\eeq

\noindent
(\Efour)
\beql{EqAGf}
\prod_{i=1}^{\frac{1}{2} n-1} (4^i +1) \quad (n \equiv 0~( \bmod~2))
\eeq

\noindent
(\Efive)
\beql{EqAGg}
\prod_{i=1}^n (2^i +1)
\eeq

\noindent
(\EfiveII)
(all weights even):
\beql{EqAGh}
2 \prod_{i=1}^{n-1} (2^i +1) \quad (n \equiv 0~( \bmod~2))
\eeq

\noindent
(\Esix)
\beql{EqAGi}
\prod_{i=0}^{\frac{1}{2} n-1}
(q^{i + \frac{1}{2} } +1) \quad (n \equiv 0~( \bmod~2))
\eeq

\noindent
(\Eseven)
\beql{EqAGj}
b \prod_{i=1}^{\frac{1}{2} n -1} (q^i +1) \quad (n \equiv 0~( \bmod~2))
\eeq
where $b=1$ if $q$ is even, 2 if $q$ is odd

\noindent
(\Eeight)
\beql{EqAGk}
\sum_{k=0}^{n/2} \sigma (n,k) 2^{k(k+1)/2} ~,
\eeq
where $\sigma (n,k)$, the number of binary self-orthogonal $[n,k]$
codes with all weights divisible by 4, is equal to 1 if $k=0$, and otherwise is given by
$$
\prod_{i=0}^{k-1}
\frac{2^{n-2i-2} + 2^{\left[ \frac{n}{2} \right] -i-1} -1}{2^{i+1} -1} ~, \quad {\rm if} \quad
n \equiv \pm 1~( \bmod~8 ) ~,
$$
$$
\prod_{i=0}^{k-1}
\frac{2^{n-2i-2} -1}{2^{i+1} -1} ~,
\quad {\rm if} \quad
n \equiv \pm 2~( \bmod~8) ~,
$$
$$
\prod_{i=0}^{k-1}
\frac{2^{n-2i-2} - 2^{\left[ \frac{n}{2} \right] -i-1} -1}{2^{i+1} -1} ~, \quad {\rm if} \quad
n \equiv \pm 3~( \bmod~8) ~,
$$
$$
\left[ \prod_{i=0}^{k-2}
\frac{2^{n-2i-2} + 2^{\frac{n}{2} -i-1} -2}{2^{i+1} -1} \right]
\cdot \left[ \frac{1}{2^{k-1}} +
\frac{2^{n-2k}+ 2^{\frac{n}{2} -k} -2}{2^k -1} \right]~,
\quad {\rm if} \quad n \equiv 0~( \bmod~8) ~,
$$
$$
\left[ \prod_{i=0}^{k-2}
\frac{2^{n-2i-2} - 2^{\frac{n}{2} -i-1} -2}{2^{i+1} -1} \right] 
\cdot \left[ \frac{1}{2^{k-1}} +
\frac{2^{n-2k}- 2^{\frac{n}{2} -k} -2}{2^k -1} \right]~,
\quad {\rm if} \quad n \equiv 4~( \bmod~8) ~.
$$
There is a similar but even more complicated formula for $T_n$ for self-dual codes over $\ZZ_4$ with Euclidean norms divisible by 8,
see \cite{Gab96}.

Formulae \eqn{EqAGb}--\eqn{EqAGj} are based on various sources including
\cite{Hoh96}, \cite{Me24}, \cite{Ple7}, \cite[Chap.~19]{MS}.
Equation \eqn{EqAGk}
is due to Gaborit\index{Gaborit, P.} \cite{Gab96}.

Here are two proofs of \eqn{EqAGb}.
(i)~Let $\sigma_{n,k}$ denote the number of $[n,k]$ self-orthogonal
codes $C$ containing {\bf 1}.
Any such $C$ can be extended to an $[n,k+1]$ self-orthogonal
code $D$ by adjoining any vector of $C^\perp \setminus C$, and any $D$ will
arise $2^k -1$ times from different $C$'s.
So we have $\sigma_{n,1} =1$,
$$\frac{\sigma_{n,k+1}}{\sigma_{n,k}} =
\frac{2^{n-2k} -1}{2^k -1} ~,
$$
and $\sigma_{n,n/2}$ gives \eqn{EqAGb}.  (ii)~A more sophisticated proof
can be obtained by observing that the Euclidean inner product induces a
symplectic geometry\index{symplectic geometry} structure on the space of even weight vectors modulo
{\bf 1}.  A self-dual code is then a maximally isotropic subspace.\index{isotropic subspace}  The
number of maximally isotropic subspaces of a symplectic geometry of
dimension $2k$ is $\Pi_{i=1}^k (2^i +1)$ \cite[\S9.4]{BCN}, and we obtain
\eqn{EqAGb} by noting that our symplectic geometry has dimension
$n-2$.~~~$\bsq$

Similarly, a binary self-dual code with weights divisible by 4 is a maximally totally singular subspace\index{totally singular subspace} of the orthogonal geometry\index{orthogonal geometry}\index{geometry!orthogonal} of dimension $n-2$
induced by $\frac{1}{2} wt(v)$, which leads to \eqn{EqAGc}.
Equations \eqn{EqAGd}, \eqn{EqAGf}, \eqn{EqAGh}, \eqn{EqAGj} are also obtained via
orthogonal geometry, \eqn{EqAGg} via symplectic geometry,\index{symplectic geometry}\index{geometry!symplectic} and \eqn{EqAGe} and \eqn{EqAGi} via unitary geometry.\index{unitary geometry}\index{geometry!unitary}

These mass formulae are useful when one is attempting to find all inequivalent codes of a given length (compare Section~\ref{GU}).
For example, suppose we are trying to find all binary self-dual
codes of length 8.
We immediately find two codes, $i_2 \oplus i_2 \oplus i_2 \oplus i_2$, where $i_2 = [11]$, and the Hamming
code\index{Hamming code}\index{Hamming code!$e_8$}\index{$e_8$, Hamming code} $e_8$, and then it appears that there are no others.
To prove this, we compute the automorphism groups of these two codes:
they have orders $2^4 4! = 384$ and $8.7.6.4 = 1344$, respectively.
We also calculate $T_8 /|G| = 3.5.9/8! = 3/896$ from \eqn{EqAGb}, and
see that indeed
$$\frac{1}{384} + \frac{1}{1344} = \frac{3}{896} ~,$$
verifying that this enumeration is complete.
We will return to this in Section~\ref{GU}.

There are also formulae that give the total number of
self-dual codes containing a fixed self-orthogonal vector or code --- see
\cite[Chapter~19]{MS}.
\subsection{Codes over $\ZZ_4$}\label{AG3}\index{codes!over $\ZZ_4$}
\hsp
Codes over rings are probably less familiar to the reader than codes over
fields, and so we will add some remarks here about the first such case,
codes over $\ZZ_4$, family \Eeight.

Any code over $\ZZ_4$ is equivalent to one with generator matrix of the form
\beql{EqA3a}
\left[
\begin{array}{ccc}
I_{k_1} & X & Y_1 + 2Y_2 \\
0 & 2I_{k_2} & 2Z
\end{array}
\right]
\eeq
where $X$, $Y_1$, $Y_2$, $Z$ are binary matrices.
Then $C$ is an elementary abelian group of type $4^{k_1}  2^{k_2}$, containing
$2^{2k_1 + k_2}$ words.
We indicate this by writing $|C| = 4^{k_1} 2^{k_2}$.
The dual code $C^\perp$ has generator matrix
$$
\left[
\begin{array}{ccc}
(-Y_1 + 2Y_2)^{tr} - Z^{tr} X^{tr} & Z^{tr} & I_{n-k_1 - k_2} \\
2X^{tr} & 2I_{k_2} & 0
\end{array}
\right]
$$
and $|C^\perp | = 4^{n-k_1 -k_2} 2^{k_2}$.

There are two binary codes $C^{(1)}$ and $C^{(2)}$ associated with $C$,
having generator matrices
\beql{EqA3b}
[I_{k_1} ~ X ~ Y_1 ] \quad {\rm and} \quad
\left[ \begin{array}{ccc}
I_{k_1} & X & Y_1 \\
0 & I_{k_2} & Z
\end{array}
\right]
\eeq
and parameters $[n,k_1]$ and $[n, k_1 + k_2]$ respectively.
If $C$ is self-orthogonal then $C^{(1)}$ is doubly-even and
$C^{(1)} \subseteq C^{(2)} \subseteq C^{(1) \perp}$.
If $C$ is self-dual then $C^{(2)} = C^{(1) \perp}$.
The next two theorems give the converse assertions.
\begin{theorem}\label{thAG1}
If $A$, $B$ are binary codes with $A \subseteq B$ then there is a code $C$ over $\ZZ_4$ with $C^{(1)} = A$, $C^{(2)} =B$.
If in addition $A$ is doubly-even and $B \subseteq A^\perp$ then $C$ can be
made self-orthogonal.
If $B= A^\perp$ then $C$ is self-dual.
\end{theorem}
\paragraph{Proof.}
Suppose $A$, $B$ have generator matrices as shown in \eqn{EqA3b}.
Then
\beql{EqA3c}
\left[
\begin{array}{ccc}
I & X & Y \\
0 & 2I & 2Z
\end{array}
\right]
\eeq
is a generator matrix for a code $C$ with $C^{(1)} = A$, $C^{(2)} = B$.
To establish the second assertion we must modify \eqn{EqA3c} to make $C$ self-orthogonal.
This is accomplished by replacing the $(j,i)$th entry of \eqn{EqA3c} by the inner
product modulo 4 of rows $i$ and $j$,
for $1 \le i \le k_1$, $1 \le j \le k_1 + k_2$, $i < j$.~~~$\bsq$

In this way every self-{\em orthogonal} doubly-even binary code corresponds to one or more self-{\em dual} codes over $\ZZ_4$.
\begin{theorem}\label{thAG2}
{\rm \cite{Gab96}}
A code $C$ over $\ZZ_4$ with generator matrix \eqn{EqA3a} is self-dual
if and only if $C^{(1)}$ is doubly-even, $C^{(2)} = C^{(1) \perp}$,
and $Y_2$ is chosen so that if $M= Y_1 Y_2^{tr}$, then $M_{ij} + M_{ji} \equiv \frac{1}{2} wt (v_i \cap v_j)$, where $v_1, \ldots, v_{k_1}$ are the generators of $C^{(1)}$.
\end{theorem}

In contrast to self-dual codes over fields, self-dual codes over $\ZZ_4$ exist for all lengths, even or odd.
Furthermore, a self-dual code $C$ over $\ZZ_4$ of length $n$ can be shortened
to a self-dual code of length $n-1$ by deleting any one of its coordinates.
This is accomplished as follows.
If the projection of $C$ onto the $i$th coordinate contains all of $\ZZ_4$, the shortened code is obtained by taking those words of $C$ that are
0 or 2 in the $i$th coordinate and omitting that coordinate.
If the projection of $C$ onto the $i$th coordinate contains only 0 and 2,
we take the words of $C$ that are 0 on the $i$th coordinate and omit that
coordinate.

In this way all self-dual codes over $\ZZ_4$ belong to a common ``family tree'', with
$i_1 = \{0,2\}$ at the root.
The beginning of this tree, showing all self-dual codes of lengths $n \le 8$,
is given in Fig.~2 of \cite{Me168}.

\section{Weight enumerators and MacWilliams theorem}\label{WE}
\hsp
\subsection{Weight enumerators}\label{Wee1}
\hsp
The {\em Hamming weight}\index{Hamming weight}\index{weight!Hamming} of a vector $u= (u_1, \ldots, u_n) \in \FF^n$, denoted by $wt (u)$, is the number of nonzero components $u_i$.

Two other types of ``weight'' are useful for studying
nonbinary codes.
For the codes in families \Eeight, \Enine\ (and hence for \Eone, \Etwo, and, if $q$
is a prime, \Eseven) we define the {\em Lee weight}\index{Lee weight}\index{weight!Lee} and {\em Euclidean norm}\index{Euclidean norm}\index{norm!Euclidean} of $u \in \FF$ by
\begin{eqnarray*}
{\rm Lee} (u) & = & \min
\{| u|, |\FF |- |u| \} ~, \\
{\rm Norm} (u) & = & ( {\rm Lee} (u))^2 ~.
\end{eqnarray*}
For a vector $u= (u_1, \ldots, u_n) \in \FF^n$, we set
\begin{eqnarray*}
{\rm Lee} (u) & = & \sum_{i=1}^n {\rm Lee} (u_i) ~, \\
{\rm Norm} (u) & = & \sum_{i=1}^n {\rm Norm} (u_i) ~.
\end{eqnarray*}
Of course, if $u$ is a binary vector, $wt (u) = {\rm Lee} (u) = {\rm Norm} (u)$.

It is customary to use the symbol $A_i$ to denote the number of vectors 
in a code $C$ having Hamming weight (or Lee weight, or Euclidean norm, depending
on context) equal to $i$.
Then $\{A_0, A_1, A_2, \ldots \}$ is called the {\em weight distribution}\index{weight distribution}
of the code.
The {\em Hamming weight enumerator}\index{Hamming weight enumerator}\index{weight enumerator!Hamming}
(abbreviated hwe)\index{hwe} of $C$ is defined to be
\beql{Eq4}
W_C (x,y) = \sum_{u \in C} x^{n-wt (u)} y^{wt(u)} =
\sum_{i=0}^n A_i x^{n-i} y^i ~.
\eeq
(The adjective ``Hamming'' is often omitted.)
There are good reasons for taking the Hamming weight enumerator\index{weight enumerator} to be a homogeneous polynomial of degree $n$ (see below).
However, no information is lost if we set $x=1$,
and write it as a polynomial in the single variable $y$.

There is an analogous definition for nonlinear codes:
for $v \in \FF^{\,n}$, let $A_i (v)$ be the number of codewords at Hamming
distance $i$ from $v$.
The {\em average Hamming weight distribution} for a nonlinear or nonadditive
code is then
$$A_i = \frac{1}{|C|} \sum_{c \in C} A_i (c) ~,$$
with associated Hamming weight enumerator
$$W_C (x,y) = \sum_{i=0}^n A_i x^{n-i} y^i ~.$$

Much more information about a code $C$ is supplied by its
{\em complete weight enumerator}\index{complete weight enumerator}\index{weight enumerator!complete} (abbreviated cwe)\index{cwe} and defined as follows.
Let the elements of the alphabet $\FF$ be $\xi_0$, $\xi_1 , \ldots, \xi_a$, and introduce
corresponding indeterminates $x_0$, $x_1 , \ldots, x_a$.
Then
\beql{Eq7}
cwe_C (x_0, \ldots, x_a) = \sum_{u \in C} x_0^{n_0 (u)} x_1^{n_1 (u)} \cdots x_a^{n_a (u)} ~,
\eeq
where $n_\nu (u)$ is the number of components of $u$ that take the value
$\xi_\nu$.

If there is a natural way to pair up some of the symbols in $\FF$ then we can often reduce
the number of variables in the cwe without losing any essential
information, by identifying indeterminates corresponding to paired symbols.
The result is a {\em symmetrized weight enumerator}\index{symmetrized weight enumerator}\index{weight enumerator!symmetrized}\index{weight enumerator!Lee} (abbreviated swe).\index{swe}
Some examples will make this clear.
For linear codes over $\FF_4$ the symmetrized weight enumerator is
\beql{Eq9}
swe_C (x,y,z) = \sum_{u \in C} x^{n_0 (u)} y^{n_1 (u)} z^{N_w (u)}
= cwe_C (x,y,z,z) ~,
\eeq
where $n_0 (u)$, $n_1 (u)$ are as above and $N_w (u)$ is the number of
components in $u$ that are equal to either $\om$ or $\overline{\om}$.
For linear codes over $\ZZ_4$, the appropriate symmetrized weight
enumerator is
\beql{Eq10}
swe_C (x,y,z) = \sum_{u \in C} x^{n_0 (u)} y^{n_\pm (u)} z^{n_2 (z)} =
cwe_C (x,y,z,y) ~,
\eeq
where
$n_\pm (u)$ is the number of components of $u$ that are equal to either $+1$ or $-1$.
There is an obvious generalization of \eqn{Eq10} to linear codes
over $\ZZ_m$.

The swe contains only about half as many variables as the complete
weight enumerator, and yet still contains enough information to determine
the Lee weight or norm distribution of a code.

All the weight enumerators mentioned so far can be obtained from the ``full weight enumerator''\index{full weight enumerator}\index{weight enumerator!full} of the code.
This is a generating function, or formal sum (not a polynomial), listing all the codewords:
$$\sum_{u \in C} z_1^{u_1} z_2^{u_2} \cdots z_n^{u_n} ~,$$
where we use a different indeterminate $z_i$ for each coordinate position.
To obtain the symmetrized weight enumerator of a code over $\FF_4$, for example,
we replace each occurrence of $z_i^0$ by $x$,
each $z_i^1$ by $y$, and each $z_i^\om$ or $z_i^{\overline{\om}}$ by $z$.

Still further weight enumerators that have proved useful can also be obtained from the full weight enumerator.
For example, the {\em split Hamming weight enumerator}\index{split Hamming weight enumerator}\index{weight enumerator!split} of a code of length $n=2m$ is
$$split_C (x,y,X,Y) =
\sum_{y \in C} x^{m- l(u)} y^{l(u)} X^{m-r(u)} Y^{r(u)} ~,
$$
where $l(u)$ (resp. $r(u)$) is the Hamming weight of the left half (resp. right half) of $u$.
Split weight enumerators have been investigated in
\cite{Me35}, for example.
\index{weight enumerator!left} \index{weight enumerator!right}
Of course, the split need not be into equal parts.
Multiply-split weight
enumerators have been extensively used in \cite{Jaffe96}.

One may also define weight enumerators for {\em translates} of codes: if $C$
is a translate of a
linear or additive code, its weight enumerator\index{weight enumerator!of translate} is
$$
W_{C}(x,y) = \sum_{c\in C} x^{n-wt(c)} y^{wt(c)}.
$$
We will use such weight enumerators later in this chapter when studying
the ``shadow'' of a self-dual code.

The {\em biweight enumerator}\index{biweight enumerator}\index{weight enumerator!biweight} of
a code generalizes the weight enumerator to consider the overlaps of pairs of codewords, and 
the joint weight enumerator\index{joint weight enumerator}\index{weight enumerator!joint}
of two codes $C$ and $D$ considers the overlaps of pairs of codewords $u \in C$ and $v \in D$.
More generally, the $k$-fold
{\em multiple weight enumerator}\index{multiple weight enumerator}\index{weight enumerator!multiple} of
a code considers the composition of $k$ codewords chosen simultaneously from the code.
Again there are generalizations of the MacWilliams and Gleason theorems
(\cite{Me27}, \cite[Chap.~5]{MS}, \cite{Shiro96}).
The connections between multiple weight enumerators of self-dual codes and Siegel
modular forms\index{modular forms} have been investigated by Duke \cite{duke}, Ozeki \cite{Oze76}, \cite{Oze4}, \cite{Oze97} and Runge \cite{run1}--\cite{run4}.

Ozeki \cite{Oze97} has recently introduced another generalization of the weight enumerator of a code $C$, namely its {\em Jacobi polynomial}.\index{Jacobi polynomial}\index{polynomial!Jacobi}
For a fixed vector $v \in \FF^n$, this is defined by
$$Jac_{C,v} (x,z) = \sum_{u \in C} x^{wt (u)} z^{wt(u \cap v)} ~,$$
which is essentially a split weight enumerator.
These polynomials have been studied in \cite{BaMO96}, \cite{BaO96},
\cite{BoMS97}.
They have the same relationship
to Jacobi forms \cite{EiZ85} as weight enumerators do to modular forms (cf.
the remarks in Section~\ref{Latt}).

For future reference we note the following relations between inner products
and weights or norms
for four of our families:

\noindent
(\Eone):
\beql{Eq12a}
(u,v) = \frac{1}{2} \{ wt(u+v) - wt(u) - wt(v) \}
\eeq

\noindent
(\Efive):
\beql{Eq12b}
(u,v) = wt (u+v) - wt (u) - wt(v)
\eeq

\noindent
(\Eeight), (\Enine):
\beql{Eq12c}
(u,v) = \frac{1}{2} \{ {\rm Norm} (u+v) - {\rm Norm} (u) - {\rm Norm} (v) \} ~.
\eeq

\subsection{Examples of self-dual codes and their weight enumerators}\label{Wee2}
\hsp
The following are some key examples of self-dual codes of the different
families mentioned in Section~\ref{SD},
together with their weight enumerators.
Some of these weight enumerators will be labeled for later reference.
Unless indicated otherwise, all the codes mentioned are self-dual codes of the appropriate kind.

We write $[n,k,d]_q$\index{$[n,k,d]_q$} to indicate a linear code of length $n$,
dimension $k$ and minimal distance $d$ over the field $\FF_q$, omitting $q$ when it is equal to 2.
$[n,k,d]_{4+}$\index{$[n,k,d]_{4+}$} indicates an additive code over $\FF_4$ containing $4^k$ vectors (so $k \in \frac{1}{2} \ZZ$).
Usually the subscript on the symbol for a code (e.g. $e_8$) gives its length.
We adopt the convention that parentheses in a vector mean that all permutations indicated by the parentheses are to be applied to that vector.\index{parentheses notation} \index{cyclic code!parenthesis notation for}
For example, in the definition of $e_8$ below, $1(1101000)$ stands for the seven vectors $11101000$, $10110100$, $10011010$, etc.
The generators for the hexacode in \eqn{Eq16b} could have been
abbreviated as $(100) (1 \om \om )$.

The following codes are all self-dual.

\noindent
(\Eone)
The first example of a binary self-dual code is the $[2,1,2]$ repetition
code\index{repetition code} $i_2 = \{ 00, 11\}$, with weight enumerator
\beql{Eq11}
W_{i_2} (x,y) = x^2 + y^2 = \phi_2 ~~{\rm (say)} ~,
\eeq
and $|Aut (i_2)| =2$.

The $[8,4,4]$ {\em Hamming} code\index{Hamming code!$e_8$}\index{code!Hamming, q.v.} $e_8$ (see Section 12 of Chapter~1;
\cite{SPLAG}, p.~80)
generated by $1(1101000)$,
is self-dual with weight enumerator
\beql{Eq12}
W_{e_8} (x,y) = x^8 + 14x^4 y^4 + y^8 = \phi_8 ~,
\eeq
and group $GA_3 (2)$ of order $8.7.6.4 = 1344$.

The $[24,12,8]$ {\em binary Golay} code $g_{24}$\index{$g_{24}$, binary Golay code}\index{code!Golay, q.v.}\index{Golay code!binary ($g_{24}$)}
(Section 12 of Chapter 1; \cite{SPLAG}, Chaps.~3, 11), generated by
\beql{Eq12aa}
1(10101110001100000000000) ~,
\eeq
or equivalently by the idempotent generator
\beql{Eq12aab}
1(00000101001100110101111)~,
\eeq
has weight enumerator
\beql{Eq13}
W_{g_{24}} (x,y) = x^{24} + 759 x^{16} y^8 + 2576 x^{12} y^{12} + 759 x^8 y^{16} + y^{24} = \phi_{24} ~.
\eeq
$Aut (g_{24} )$ is the Mathieu group\index{Mathieu group}\index{group!Mathieu} $M_{24}$, of
order $24.23.22.21.20.48 = 244823040$.

All three codes
$i_2$, $e_8$, $g_{24}$ are unique in the sense that any linear or nonlinear
code with the same length, size and minimal distance and containing the zero vector is linear and equivalent to the code given above \cite{Ple8}
(see also \cite{snover}).

\noindent
(\Etwo) Self-dual codes over $\FF_3$ exist if and only if the length $n$ is
a multiple of 4 (this follows from Gleason's theorem,
see \eqn{EqG11a}, and is also a
consequence of the argument used to prove \eqn{EqAGd} \cite{Ple7}).
We use indeterminates $x,y$ for the Hamming weight enumerator $W(x,y)$ and
$x,y,z$ for the cwe.

The $[4,2,3]_3$ {\em tetracode} $t_4$,\index{tetracode $t_4$}\index{$t_4$, tetracode}\index{code!tetracode, q.v.} generated by $\{1110, 0121\}$ (Section 7 of Chapter 1; \cite{SPLAG} p.~81) has
\beql{Eq14}
W_{t_4} (x,y) = x^4 + 8xy^3
\eeq
and
cwe $x \{x^3 + (y+z)^3 \}$.
$Aut (t_4) = 2.S(4)$, where $S(n)$\index{group!symmetric $S(n)$} denotes a symmetric group of order $n!$.

The $[12,6,6]_3$ {\em ternary Golay} code\index{ternary Golay code}\index{$g_{12}$, ternary Golay code}\index{Golay code!ternary ($g_{12}$)}
$g_{12}$ (Section 12 of Chapter 1; \cite{SPLAG}, p.~85),
generated by
$1(11210200000)$, has
\beql{Eq15}
W_{g_{12}} (x,y) =
x^{12} + 264 x^6 y^6 + 440 x^3 y^9 + 24y^{12}
\eeq
and (assuming the all-ones codeword is present)
\beql{Eq16}
cwe (x,y,z) = x^{12} + y^{12} + z^{12} + 22 (x^6 y^6 + y^6 z^6 + z^6 x^6 ) + 220
(x^6 y^3 z^3 + x^3 y^6 z^3 + x^3 y^3 z^6) ~.
\eeq
$Aut (g_{12} ) = 2.M_{12}$ (where $M_{12}$ is a Mathieu group),\index{Mathieu group}\index{group!Mathieu} of order 190080.

These two codes are unique in the same sense as our binary examples \cite{Ple8}.

\noindent
(\Ethree) We use indeterminates $x,y$ for the Hamming weight enumerator,
$x,y,z$ for the swe and $x,y,z,t$ (corresponding to the symbols $0,1, \om, \overline{\om}$) for the cwe, so that
$swe (x,y,z) = cwe (x,y,z,z)$.

The $[2,1,2]_4$ repetition code $i_2 = \{ 00, 11, \om \om, \oom \oom \}$ has
\begin{eqnarray}
\label{Eq16a}
W_{i_2} (x,y) & = & x^2 + 3y^2~, \nonumber \\
swe & = & x^2 + y^2 + 2z^2 ~, \nonumber \\
cwe & = & x^2 + y^2 + z^2 + t^2 ~,
\end{eqnarray}
and a group of order 12.

The $[6,3,4]_4$ {\em hexacode} $h_6$\index{hexacode $h_6$}\index{$h_6$, hexacode}\index{code!hexacode, q.v.} (Section 12 of Chapter 1, \cite[p.~82]{SPLAG}) in the form with generator matrix
\beql{Eq16b}
\left[
\matrix{
1& 0 & 0 & 1 & \om & \om \cr
0 & 1 & 0 & \om & 1 & \om \cr
0 & 0 & 1 & \om & \om & 1 \cr
}
\right]
\eeq
has
\begin{eqnarray}\label{Eq17}
W_{h_6} (x,y) & = &
x^6 + 45 x^2 y^4 + 18 y^6 , \\
\label{Eq18}
swe & = & x^6 + y^6 + 2z^6 + 15 (2x^2 y^2 z^2 + x^2 z^4 + y^2 z^4 ) ~, \\
\label{Eq19}
cwe & = & x^6 + y^6 + z^6 + t^6 + 15
(x^2 y^2 z^2 + x^2 y^2 t^2 + x^2 z^2 t^2 + y^2 z^2 t^2 )
\end{eqnarray}
and $Aut (h_6) = 3.S(6)$, of order 2160.

Again these codes are unique.

Of course this $i_2$ is simply the $\FF_4$-span of the binary code $i_2$ defined above.
\index{tensor product}
In general, if $C$ is defined over an alphabet $\FF$, and $\FF' \supseteq \FF$ is a larger alphabet, we write $C \otimes \FF'$ to indicate this process.

If $C$ is a binary self-dual code then $C \otimes \FF_4$ is a self-dual code belonging to both families \Ethree\ and \Efour.
Conversely, it is not difficult to show that if $C$ is self-dual over $\FF_4$ with respect to both the Hermitian and Euclidean inner products, then $C= B \otimes \FF_4$ for some self-dual binary code $B$.

\noindent
(\Efour) The $[4,2,3]_4$ Reed-Solomon code\index{Reed-Solomon code}
\index{code!Reed-Solomon}
$$
\left[
\matrix{
1 & 1 & 1 & 1 \cr
0 & 1 & \om & \oom \cr
}
\right]
$$
has
\begin{eqnarray*}
W(x,y) & = & x^4 + 12xy^3 + 3y^4 ~, \\
swe & = & x^4 + y^4 +~ 2z^4 + 12 xyz^2 ~, \\
cwe & = & x^4 + y^4 + z^4 + t^4 + 12 xyzt ~.
\end{eqnarray*}
The automorphism group is $3.S(4)$, of order 72.

\noindent
(\Efive) The smallest example is the $[1, \frac{1}{2}, 1]_{4+}$
code $i_1 =\{0,1\}$, with automorphism group of order 2 (conjugation).
The $[12,6,6]_{4+}$ {\em dodecacode}\index{dodecacode $z_{12}$}\index{code!dodecacode, q.v.}
\index{$z_{12}$, dodecacode}
\index{$Z(n)$, cyclic group}
$z_{12}$
can be defined as the cyclic code with generator
$\om 10100100101$
(\cite{Me223}, see also \cite{Hoh96}).
\index{group!cyclic $Z(n)$}
\index{cyclic group $Z(n)$}
$Aut (z_{12} )$ is a semi-direct product of $Z(3)^3$ with $S(4)$ (where $Z(n)$ denotes a cyclic group of order $n$) and has order 648.

\noindent
(\Esix) Since the norm map from $\FF_q$ to $\FF_{\sqrt{q}}$ is surjective, there is an element $a \in \FF_q$ with $a \overline{a} = -1$.
Then $[1a]$ is self-dual.

\noindent
(\Eseven) As in family \Ethree,
there is a restriction on $n$:
if $q \equiv 3$ $(\bmod ~4)$ then self-dual codes exist if and only if $n$ is a multiple of 4;
for other values of $q$, $n$ need only be even \cite{Ple7}.
Provided $q \not\equiv (3)$ $(\bmod~4)$, $\FF_q$ contains an element $i$ such that $i^2 = -1$, and then $[1i]$ is self-dual.

\noindent
(\Eeight) The smallest example is the self-dual code $i_1 =\{ 0,2\}$
of length 1.
The {\em octacode} $o_8$\index{octacode $o_8$}\index{$o_8$, octacode}\index{code!octacode, q.v.} (\cite{SPLAG}, \cite{Me168}) is the length 8 code generated by the vectors $3(2001011)$, or equivalently with generator matrix
\beql{Eq19a}
\left[
\matrix{
1& 0 & 0 & 0 & 2 & 1 & 1 & 1 \cr
0 & 1 & 0 & 0 & 3 & 2 & 1 & 3 \cr
0 & 0 & 1 & 0 & 3 & 3 & 2 & 1 \cr
0 & 0 & 0 & 1 & 3 & 1 & 3 & 2 \cr
}
\right]~,
\eeq
having minimal Lee weight 6 and minimal norm 8,
$$swe = x^8 + 16y^8 + z^8 + 14 x^4 z^4 + 112 xy^4 z (x^2 + z^2) ~,$$
and $|Aut (o_8)| =2.1344$.

The most interesting property of the octacode is that when mapped to a binary code under the Gray map\index{Gray map}\index{map!Gray}
\beql{EqLift2}
0 \to 00 , \quad 1 \to 01, \quad 2 \to 11, \quad 3 \to 10~,
\eeq
$o_8$ becomes the Nordstrom-Robinson code,\index{Nordstrom-Robinson code}\index{code!Nordstrom-Robinson}
a nonlinear binary code of length 16,
minimal distance 6, containing 256 words
(Section~14 of Chapter 1, Chapter~xx (Helleseth-Kumar), \cite{Me183}, \cite{Me184}).
The latter is therefore a formally self-dual binary code, see Section~\ref{SDM}.

The octacode\index{octacode $o_8$} reduces $\bmod~2$ to the Hamming code $e_8$.
\index{Hamming code!$e_8$}
There is another lift of $e_8$ to $\ZZ_4$, namely the code $\sE_8$, with generator matrix
\beql{EqLift1}
\left[ \matrix{
1 & 0 & 0 & 0 & 0 & 1 & 1 & 1 \cr
0 & 1 & 0 & 0 & 3 & 0 & 1 & 3 \cr
0 & 0 & 1 & 0 & 3 & 3 & 0 & 1 \cr
0 & 0 & 0 & 1 & 3 & 1 & 3 & 0 \cr
}
\right]~,
\eeq
but the minimal Lee weight and norm are now both only 4.
However, not all binary self-dual codes lift to self-dual codes over
$\ZZ_4$, e.g. $\{00, 11\}$ does not.
\begin{theorem}\label{Lift1}
(a) Let $C$ be a binary self-dual code of length $n$.
A necessary and sufficient condition for $C$ to be lifted to a self-dual code
$\hat{C}$ over $\ZZ_4$ is that all weights in $C$ are divisible by 4.
(b)~If this condition is satisfied, $\hat{C}$ can be chosen so that all norms
are divisible by 8.
(c)~More generally,
a self-dual code over $\ZZ_m$, $m$ even,
that reduces to a self-dual code $\bmod~2$ lifts to $\ZZ_{2m}$ precisely when all norms are divisible by $2m$, and in that
case all norms in the lifted code can be arranged to be divisible by $4m$.
Thus if a code lifts from $\ZZ_m$ to $\ZZ_{2m}$ then it lifts to
$\ZZ_{2^k m}$ for all $k$.
In particular, if a binary code lifts to $\ZZ_4$ then it lifts to a self-dual code over the 2-adic integers.
\end{theorem}
\paragraph{Proof.}

(a) (Necessity) Suppose $v \in C$ has weight $wt (v) \not\equiv 0$ $(\bmod~4)$, and let $\hat{v} \in \hat{C}$ be any lift of $v$.
Then ${\rm Norm} ( \hat{v}) \equiv {\rm Norm} (v)$ $(\bmod~4)$ because
for integers $x,y$ if $x \equiv y$ $(\bmod~2)$ then $x^2 \equiv y^2$ $(\bmod~4)$.

\noindent
(Sufficiency) Without loss of generality $C$ has a generator matrix of the form $[IA]$ where $AA^{tr} \equiv -I$ $(\bmod~2)$.
Let $B$ be any lift of $A$ to $\ZZ_4$.
We wish to find $\hat{A} = B +2M$ such that $\hat{A} \hat{A}^{tr} \equiv -I$ $(\bmod~4)$, since then we can take $\hat{C} = [I \hat{A}]$.
We have
$$\hat{A} \hat{A}^{tr} \equiv B B^{tr}
+2 (MB^{tr} + BM^{tr}) ~ (\bmod~4) ~.
$$
The condition on $C$ implies that $BB^{tr} +I$ has even coefficients and is zero on the diagonal.
But then there exists a binary matrix $M'$ such that $2(M' + M'^{tr}) = BB^{tr} +I$,
and we take
$M=M' (B^{-1})^{tr}$.
This completes the proof of (a).

(b) We need to show that we can choose $\hat{A}$ so that the diagonal entries of $\hat{A} \hat{A}^{tr} +I$ are zero $\bmod~8$.
Set $\hat{A}' = \hat{A} - 2L \hat{A}$, where $L$ is symmetric,
so that
$$\hat{A}' (\hat{A}')^{tr} = \hA \hA^{tr} + 4L + 4L^2 ( \bmod~8) ~.
$$
Let $\Delta = \frac{1}{4} ( \hA \hA^{tr} +I)$.
Then we need $L^2 + L + \Delta$ $(\bmod~2)$ to be symmetric with zero diagonal.
It is easy to see that we can accomplish this provided ${\rm trace} (\Delta) \equiv 0$ $(\bmod~2)$ (consider, for instance, $L= \left( {1 \atop 1}~ {1 \atop 0} \right)$.)
In fact, we have
$$1 \equiv \det ( \hA \hA^{tr} ) \equiv 1 + 4 ~ {\rm trace} ~ \Delta ~ (\bmod~8)$$
so trace $\Delta$ is even.

The proof of (c) is analogous.~~~$\bsq$

It follows from Theorem~\ref{Lift1} that the Golay code\index{Golay code} $g_{24}$
\index{lifting to $\ZZ_4$}
\index{Golay code!over $\ZZ_4$}
can be lifted to $\ZZ_4$.
Since $g_{24}$ is an extended cyclic code, the lift can be easily
performed by Graeffe's method\index{Graeffe's method} \cite{Me168}, \cite{Usp}.
Suppose $g_2 (x)$ divides $x^n -1$ $(\bmod~2)$, and we wish to find a monic polynomial $g(x)$ over $\ZZ_4$ such that $g(x) \equiv g_2 (x)$ $(\bmod~2)$ and
$g(x)$ divides $x^n -1$ $(\bmod~4)$.
Let $g_2 (x) = e(x) - d(x)$, where $e(x)$ contains only even powers
and $d(x)$ only odd powers.
Then $g(x)$ is given by
$g(x^2) = \pm (e^2 (x) -d^2 (x))$.
Applying this technique to the generator polynomial for $g_{24}$,
that is, to $g_2 (x) = 1+x+x^5 +x^6 + x^7 + x^9 + x^{11}$ (see
\eqn{Eq12aa}), we obtain
$g(x) = -1 + x + 2x^4 -x^5 - x^6 - x^7 - x^9 + 2x^{10} + x^{11}$, and so
\beql{Eq12ab}
3(31002333032100000000000)
\eeq
generates a self-dual code $G_{24}$ of length 24 which is the Golay code lifted to $\ZZ_4$.
Iterating this process enables us to lift cyclic or extended cyclic codes to
$\ZZ_{2^m}$ for arbitrarily large $m$.

\noindent
(F1)
Let $q=5$.
Then
\beql{Eq34a}
\left[
\matrix{ 1& -v & v & -1 & 0 & 0 \cr
0 & 1 & -v & v & -1 & 0 \cr
0 & 0 & 1 & -v & v & -1 \cr
}
\right]
\eeq
where $v= (1+u)/2$, generates a self-dual code of length 6 over $\FF_5 [u] / (u^2)$.

The matrix \eqn{Eq34a} also generates self-dual codes from family \Esix.
Suppose $q$ is a prime power such that $v^2 - v-1$ has no solution in $\FF_q$, and let $v$ be a solution in $\FF_{q^2}$.
Then \eqn{Eq34a} defines a Hermitian self-dual code over $\FF_{q^2}$ with minimal distance 4.
In the case $q=2$ we get the hexacode.\index{hexacode $h_6$}

\noindent
(F3) The 2-adic Hamming code\index{2-adic Hamming code}\index{Hamming code!2-adic}
\cite{Me191} is the self-dual code of length 8 with generator matrix
$$
\left[
\matrix{
1 & \la & \la -1 & -1 & 0 & 0 & 0 & 1 \cr
0 & 1 & \la & \la-1 & -1 & 0 & 0 & 1 \cr
0 & 0 & 1 & \la & \la-1 & -1 & 0 & 1 \cr
0 & 0 & 0 & 1 & \la & \la -1 & -1 & 1 \cr
} \right] ~,
$$
where $\la$ is the 2-adic integer $(1+ \sqrt{-7})/2$.
The 2-adic expansion of $\la$ is
$$\la = 2+4+32+128+256+512+1024+2048+4096+32768+ \cdots$$
This is the cyclic code with generator
$$1, \la , \la -1, -1, 0,0,0$$
with a 1 appended to each of the generators.

Similarly, the 2-adic self-dual Golay code of length 24 is the cyclic
\index{Golay code!2-adic}
\index{2-adic Golay code}
code with generator
$$1,1- \la , -2-\la, -4, \la-4, 2\la-3, 2\la+1, \la+3, 4, 3- \la, - \la, -1, 0,0,0,0,0,0,0,0,0,0,0~,
$$
where now $\la = (1+ \sqrt{-23})/2$,
with a 1 appended to each of the 12 generators.

The 3-adic self-dual Golay code of length 12 is the cyclic code with generator
\index{Golay code!3-adic}
$$1, \la , -1 , 1, \la -1 , -1, 0,0,0,0,0~,$$
where $\la = (1+ \sqrt{-11})/2$, again with a 1 appended to each generator.

\noindent
(F4)
We shall not discuss these codes here, but refer the reader to Wood \cite{Wood96}.

\subsection{MacWilliams Theorems}\label{SDM}
\hsp
MacWilliams\index{MacWilliams, F. J.} (\cite{MacW64}; see also \cite{MS}) discovered that the Hamming
\index{MacWilliams theorem}\index{theorem!MacWilliams}
weight distribution of the dual of a linear code is determined just by the Hamming weight distribution of the code.
There are versions of this theorem for most of our families of codes.
Although there are several ways to state these identities, the simplest
formulation is always in terms of the weight enumerator polynomials
(it is for this reason that we insist that the weight enumerator should be
a homogeneous polynomial).
\begin{theorem}\label{thM1}
{\rm (MacWilliams and others.)}

\noindent
{(\Eone)}
Three equivalent formulations of the result for binary self-dual codes are:
\begin{eqnarray}
\label{EqM1}
W_{C^\perp} (x,y) & = & \frac{1}{|C|} W_C (x+y, x-y) ~, \\
\label{EqM2}
\sum_{u \in C^\perp} x^{n-wt(u)} y^{wt(u)} & = &
\frac{1}{|C|} \sum_{u \in C} (x+y)^{n-wt(u)} (x-y)^{wt(u)} ~,
\end{eqnarray}
and, if $\{A_0^\perp , A_1^\perp , \ldots \}$ is the weight distribution of $C^\perp$,
\beql{EqM3}
A_k^\perp = \frac{1}{|C|} \sum_{i=0}^n A_i P_k (i)
\eeq
where
$$P_k (x) = \sum_{j=0}^k (-1)^j {\binom{x}{j}} {\binom{n-x}{k-j}}, ~~
k=0, \ldots, n ~,
$$
is a Krawtchouk polynomial\index{Krawtchouk polynomial}\index{polynomial!Krawtchouk}
(\cite{MS}, Chap.~5; etc.).
There are analogous Krawtchouk polynomials for any alphabet, see \cite{MS}, p.~151.
For the remaining cases we give just the formulation in terms of weight enumerators.

\noindent
{\rm (\Etwo)}
\begin{eqnarray*}
W_{C^\perp} (x,y) & = & \frac{1}{|C|} W_C (x+2y, x-y)~, \\
cwe_{C^\perp} (x,y,z) & = & \frac{1}{|C|} cwe_C (x+y+z, x+\om y + \oom z , x+ \oom y + \om z) ~.
\end{eqnarray*}

\noindent
{\rm (\Ethree)} and {\rm (\Efive)}
\begin{eqnarray*}
W_{C^\perp} (x,y) & = & \frac{1}{|C|} W_C (x+3y , x-y) ~, \\
swe_{C^\perp} (x,y,z) & = & \frac{1}{|C|} swe_C (x+y+2z, x+y-2z, x-y) ~, \\
cwe_{C^\perp} (x,y,z,t) & = & \frac{1}{|C|} cwe_C (x+y+z+t,
x+y-z-t, x-y+z-t, x-y-z+t)~.
\end{eqnarray*}

\noindent
{\rm (\Efour)}
\begin{eqnarray*}
W_{C^\perp} (x,y) & = & \frac{1}{|C|} W_C (x+3y , x-y) ~, \\
swe_{C^\perp} (x,y,z) & = & \frac{1}{|C|} swe_C (x+y+2z, x+y-2z, x-y) ~, \\
cwe_{C^\perp} (x,y,z,t) & = &
\frac{1}{|C|} cwe_{C^\perp} (x+y+z+t, x+y-z-t, x-y-z+t, x-y+z-t) ~.
\end{eqnarray*}

\noindent
{\rm (\Esix)}
\beql{EqM4}
W_{C^\perp} (x,y) = \frac{1}{|C|} W_C (x+(q-1) y, x-y) ~.
\eeq
Let $\la$ be a nontrivial linear functional from $\FF_q$ to $\FF_p$, and set
\beql{EqM5}
\chi_\beta (x) = e^{2 \pi i \la ( \beta \ox) / p} ~.
\eeq
The cwe for $C^\perp$ is obtained from the cwe for $C$ by replacing each $x_j$ by
$$\sum_{k=0}^{q-1} \chi_{\xi_j} (\xi_k) x_k ~.$$
(We omit discussion of the swe, since there are several different ways in which it might be defined.)

\noindent
{\rm (\Eseven)} Same as for \Esix, but omitting the bar in \eqn{EqM5}.

\noindent
{\rm (\Eeight)}
\begin{eqnarray*}
W_{C^\perp} (x,y) &=& \frac{1}{|C|} W_C (x+3y , x-y) \\
swe_{C^\perp} (x,y,z) & = & \frac{1}{|C|} swe_C (x+2y+z, x-y, x-2y+z) \\
cwe_{C^\perp} (x,y,z,t) & = &
\frac{1}{|C|} cwe_{C} (x+y+z+t, x+iy-z-it, x-y+z-t, x-iy-z+it) ~.
\end{eqnarray*}

\noindent
{\rm (\Enine)}
$$
W_{C^\perp} (x,y) = \frac{1}{|C|} W_C (x+(m-1)y, x-y) ~.$$
The cwe for $C^\perp$ is obtained from the cwe for $C$ by replacing
each $x_j$ by
\beql{EqM9.1}
\sum_{k=0}^{m-1} e^{2 \pi i jk/m} x_k ~.
\eeq
\end{theorem}

\paragraph{Proof.}
We prove the result for family \Eone.
There are analogous proofs for the other cases,
cf. Section~10 of Chapter~1, Section~8 of Chapter~xx (Helleseth-Kumar),
\cite{Lint}, \cite[Chap. 5]{MS}.

\index{Fourier transform}
\index{transform!Fourier}
\index{Hadamard transform}
\index{transform!Hadamard}
Let $f$ be a polynomial-valued function on $\FF_2^n$.
Define the {\em Fourier} (or {\em Hadamard}) {\em transform} of $f$ by
$$\hat{f} (u) = \sum_{v \in \FF_2^n} (-1)^{u.v} f(v), \quad
u \in \FF_2^n ~.$$
If $C$ is a linear code it is straightforward to verify that
\beql{EqM6}
\sum_{u \in C^\perp} f(u) = \frac{1}{|C|} \sum_{u \in C}
\hat{f} (u) ~.
\eeq
(This is a version of the Poisson summation formula\index{Poisson summation formula} --- cf. \cite{Dym}.)
Now we set $f(u) = x^{n-wt(u)} y^{wt(u)}$, and after some algebra
(the details can be found on p.~126 of \cite{MS}) discover that
\beql{EqM7}
\hat{f} (u) = (x+y)^{n-wt(u)} (x-y)^{wt(u)} ~.
\eeq
Equations \eqn{EqM6} , \eqn{EqM7} together imply \eqn{EqM2}.~~~$\bsq$
\subsection*{Examples}
\hsp
(a) The {\em repetition} code $C$ over a field $\FF_q$ has Hamming weight
enumerator
$$W_C (x,y) = x^n + (q-1) y^n ~,$$
so from \eqn{EqM4} we deduce that the dual code $C^\perp$, the {\em zero-sum} code,\index{zero-sum code}\index{code!zero-sum}
has weight enumerator
$$W_{C^\perp} (x,y) = \frac{1}{q} \left\{
(x+(q-1)y)^n + (q-1) (x-y)^n \right\} ~.
$$
Note that when $n=2$, $W_{C^\perp} = W_C$ (compare case (e) of Theorem~\ref{thGP1}).

(b) The binary codes $i_2$ and $e_8$ are self-dual, and indeed one easily
verifies that their weight enumerators $x^2 + y^2$ \eqn{Eq11} and
$x^8 + 14 x^4 y^4 + y^8$ \eqn{Eq12} are left unchanged if $x$ and $y$
are replaced by $(x+y) / \sqrt{2}$ and $(x-y) / \sqrt{2}$.
\subsection*{Remarks}
\hsp
1. The map that sends $W_C (x,y)$ to $\frac{1}{|C|} W_C (x+y, x-y)$, or that sends
$\{ A_0 , A_1 , \ldots \}$ to $\{A_0^\perp , A_1^\perp, \ldots \}$ as in \eqn{EqM3}, is often called the {\em MacWilliams}\index{MacWilliams transform} or
\index{transform!Krawtchouk}
\index{transform!MacWilliams}
{\em Krawtchouk} transform.\index{Krawtchouk transform}
A remarkable theorem of Delsarte \cite{MS350} --- see
Chapters xx (Brouwer), yy (Camion), zz (Levenshtein) --- shows that this transform is useful even for nonlinear codes.

2. For the families \Eone, \Ethree, \Efour\ and \Efive\ all the MacWilliams transforms
have order 2, as they do for the Hamming weight enumerators for families \Etwo\ and
\Eeight\ and the swe for \Eeight.
For the cwe in families \Etwo\ and \Eeight\ the square of the MacWilliams
transform takes $x_j$ to $x_{-j}$.
However, this does not change the cwe of the code, and so, in all cases,
if the MacWilliams transform is applied twice, the original weight enumerator is recovered.

3. The identity for the swe in family \Enine\ is left to the reader.
For (F1) we refer to Bachoc \cite{Bach97} and for
(F4) to Wood \cite{Wood96}.
Duality fails for (F2) and weights are undefined in case (F3).

4. Shor and Laflamme \cite{ShLa96} show that there is an analogue of the MacWilliams identity for quantum codes.\index{quantum code}\index{code!quantum}
There is also an analogue of the shadow \cite{RainsQSE}.

\subsection{Isodual and formally self-dual codes}\label{Iso}
\hsp
Following \cite{Me177}, we say that a linear code which
is equivalent to its dual is {\em isodual}.\index{isodual}
A (possibly nonlinear) code with the property that its weight
enumerator coincides with its MacWilliams transform is called {\em formally self-dual}.
\index{code!isodual}\index{code!formally self-dual}
An isodual code is automatically formally self-dual.\index{formally self-dual}

It is easy to prove that any self-dual code from family \Eeight\ produces a
formally self-dual binary code using the Gray map\index{Gray map} \eqn{EqLift2}
(\cite{Me183}, \cite{Me184}).
As already mentioned in Section~\ref{Wee2}, the octacode $o_8$ produces the (formally self-dual) Nordstrom-Robinson\index{Nordstrom-Robinson code} code in this way.
Similarly,\footnote{We are indebted to Dave Forney for these remarks.}
a self-dual code from family \Efive\ produces an isodual binary code using the map
\beql{EqLift3}
0 \to 00 , ~~ 1 \to 11 , ~~ \om \to 01, ~~
\oom \to 10 ~.
\eeq
We give several examples of this construction.

(i)~The code $d_3^+$ (see Section~\ref{GU7}) produces the isodual $[6,3,3]$ binary code with generator matrix
\beql{EqIso1}
\left[
\begin{array}{lll}
11 & 11 & 00 \\
11 & 00 & 11 \\
10 & 10 & 10
\end{array}
\right] ~.
\eeq
(The dual, which is a different code,
is obtained by interchanging the last two
columns.)

(ii)~The shortened hexacode, $h_5$, (see Section~\ref{EX4}) produces an isodual $[10,5,4]$ code.

(iii)~The hexacode\index{hexacode $h_6$} $h_6$ produces an isodual $[12,6,4]$ code.
There is an additive but not linear version of the hexacode, $h'_6$,
found by Ran and Snyders\index{code!Ran and Snyders}\index{Ran and Snyders code} \cite{RaS96}, generated by $(0011 \om \oom )$, which
under the map \eqn{EqLift3} produces a second, inequivalent, isodual $[12,6,4]$ code.
As members of the family \Efive, however, $h_6$ and $h'_6$ are equivalent.

Further examples of formally self-dual codes will be mentioned in Remark 4
following Theorem~\ref{thGP1}.  Isodual and formally self-dual codes have
also been studied in
\cite{DHo97}, \cite{gh0.1}, \cite{GuHa97a}, \cite{Hara96a},
\cite{KeP94}, \cite{MS}, \cite{Me64} (see also
\cite{Me177}).

\section{Restrictions on weights}\label{GP}
\hsp
\subsection{Gleason-Pierce Theorem}\label{GP1}\index{Gleason-Pierce theorem}\index{theorem!Gleason-Pierce}
\hsp
It is elementary that in a binary self-orthogonal code the weight of every vector is even, in a ternary self-dual code the weight of every vector is a multiple of 3, and in a Hermitian self-dual code over $\FF_4$ the weight of every vector is even.
\index{weights!divisibility of}
Furthermore, there are many well-known binary self-dual codes whose
weights are divisible by 4 --- see above.
The following theorem, due to Gleason and Pierce, shows that these four are essentially the only possible
nontrivial
divisibility restrictions that can be imposed on the weights of self-dual codes.
\begin{theorem}\label{thGP1}
{\rm (Gleason and Pierce \cite{AMT67}.)}
If $C$ is a self-dual code belonging to any of the families \Eone\ through
\Enine\ which has all its Hamming weights divisible by an integer $c >1$ then one of the following holds:
$$
\begin{array}{lll}
(a) & |\FF | =2 , \quad c=2 & (\mbox{so family \Eone}) \\ [+.1in]
(b) & |\FF| = 2, \quad c=4 & (\mbox{so family \Eone}) \\ [+.1in]
(c) & | \FF| =3 , \quad c=3 & (\mbox{so family \Etwo}) \\ [+.1in]
(d) & | \FF | =4, \quad c=2 & (\mbox{so families \Ethree, \Efour, \Efive, \Eeight}) \\ [+.1in]
(e) & \multicolumn{2}{l}{| \FF| = q, \quad \mbox{$q$ arbitrary}, \quad c=2, \quad {\rm and}}
\end{array}
$$
the Hamming weight enumerator of $C$ is
$$(x^2 + (q-1) y^2)^{n/2} ~.$$
\end{theorem}

\paragraph{Remarks.}

1. The theorem may be proved by considering how the Hamming weight enumerator
behaves under the MacWilliams transform --- see \cite{Me58} for details.
An alternative proof of a somewhat more general result is given in \cite{Wa81} ---
see Theorem~13.5 of Chapter~xx (Ward).

2. The same conclusion holds if ``$C$ is self-dual'' is replaced by ``$C$ is formally self-dual''.

3. Note that there are no nontrivial examples from families \Esix, \Eseven\ or \Enine.

4. There are several points to be mentioned concerning case (e).
{\em Linear} self-dual codes with weight enumerator $(x^2 + (q-1) y^2)^{n/2}$ always exist in families \Eone, \Ethree, \Efour, \Efive, \Esix;
exist in families \Eseven\ and \Enine\
precisely when there is a square root of $-1$ in $\FF_q$ or $\ZZ_m$ respectively; in particular, 
they never exist in families \Etwo\ or \Eeight.

Furthermore, it is easy to see that any linear code over $\FF_q$ for $q> 2$ with weight enumerator $(x^2 + (q-1) y^2)^{n/2}$ is a direct sum of codes of length 2.
However, in the binary case there are many examples of linear codes with weight enumerator $(x^2 + y^2)^{n/2}$ that are not self-dual:
these have been classified for $n \le 16$, see \cite{Me64}.
These are examples of formally self-dual codes: see Section~\ref{Iso}.
There are also examples from family \Efive, e.g. the additive code
$[1100, 0110, 0011, \om\om\om\om]$ with weight enumerator $(x^2+ 3y^2)^2$.

5. In some cases, analogous restrictions can be imposed on Euclidean norms of codewords.
In particular, suppose $C$ is a self-dual code over $\ZZ_m$ (that is, a code from
families \Eeight\ or \Enine) where $m$ is even.
Then the Euclidean norms of the codewords {\em must} be divisible by $m$, and {\em may} be divisible by $2m$ (\cite{BDHO}, \cite{BSBM97}, \cite{DGH97}, see also Theorem~\ref{Lift1}).

6. Codes from family (F1) with $q = 2$ can also satisfy case (d)
of the theorem, since they can be embedded in family \Efive\ via
the map $a+bu\to a+b\omega$.

\subsection*{Examples}
\hsp
Many of the examples given in Section~\ref{Wee2} satisfy one of these
divisibility conditions:

(\Eone): all self-dual codes satisfy (a), and $e_8$ and $g_{24}$ satisfy (b).
Note that any code satisfying (b) is self-orthogonal (from \eqn{Eq12a}).

(\Etwo): a code satisfies (c) precisely when it is self-orthogonal

(\Ethree): a code satisfies (d) precisely when it is self-orthogonal

(\Efour): a self-dual code satisfying (d) is a linearized binary code

(\Efive): The dodecacode $z_{12}$ satisfies (d).
Any code satisfying (d) is self-orthogonal.

\subsection{Type I and Type II codes}\label{GP2}
\index{code!Type I}
\index{code!Type II}
\index{code!doubly-even}
\index{code!singly-even}
\hsp
A binary self-dual code $C$ with all weights divisible by 4 is called
{\em doubly-even}\footnote{The unqualified
term ``even'' has been used to denote both Type I and Type II codes,
and is therefore to be avoided when speaking of self-dual codes.
\index{code!even}
Use
``singly-even'' or ``doubly-even'' instead.}
or of {\em Type II};\index{Type II}
if we do not impose this restriction then
$C$ is {\em singly-even} or of {\em Type I}.\index{Type I}
We denote these two families by \EoneI\ and \EoneII.
A Type I code may or may not also be of Type II:
the classes are not mutually exclusive.
We say a code is {\em strictly Type I}\index{strictly Type I} if it is not of Type II.

Similarly, we will say that a self-dual code over $\ZZ_m$, $m$ even, from the
\index{\EeightI\ (Type I self-dual over $\ZZ_4$)}
\index{\EeightII\ (even self-dual over $\ZZ_4$)}
\index{\EnineI\ (Type I self-dual over $\ZZ_m$)}
\index{\EnineII\ (Type I self-dual over $\ZZ_m$)}
families \Eeight\ or \Enine\ is of {\em Type~II} if the Euclidean norms
are divisible by $2m$,
or of {\em Type~I} if they are divisible by $m$.
(This terminology was introduced in \cite{BDHO}, \cite{BSBM97}, \cite{DGH97}.)
We denote these families by \EeightI\ (or \EnineI)
and \EeightII\ (or \EnineII).

There is one other situation where a similar distinction can be made.
\index{\EfiveI\ (Type I additive self-dual over $\FF_4$)}
\index{\EfiveII\ (Type II additive self-dual over $\FF_4$)}
An additive trace-Hermitian self-dual code over $\FF_4$ from the family
\Efive\ is of {\em Type II} if the Hamming weights are even,
or of {\em Type I} if odd weights may occur
(if odd weights do occur then the code cannot be linear).
We denote these two families by \EfiveI\ and \EfiveII.

More generally, we will say that a binary code is {\em doubly-even}
if all its weights are divisible by 4, or {\em singly-even} if its weights are even.
It follows from \eqn{Eq12a} that a doubly-even code is necessarily self-orthogonal (and from \eqn{Eq12b} and \eqn{Eq12c} that Type~II codes over $\ZZ_m$ and $\FF_4$ are necessarily self-orthogonal).

In view of Theorem~\ref{thGP1}, in the past self-dual codes over $\FF_3$ have been
\index{code!Type III}\index{code!Type IV}
\index{\EoneI\ (singly-even self-dual)}
\index{\EoneII\ (doubly-even self-dual)}
called {\em Type III} codes, and Hermitian self-dual codes over $\FF_4$ have been called {\em Type IV} codes.
However, we shall not use that terminology in this chapter.

\section{Shadows}\label{Shad}
\hsp
In the three cases where we can define a Type~II code (see the previous section) we can also define a certain canonical translate of a code called its shadow \cite{Me158}.
The weight enumerator of the shadow can be obtained from the weight enumerator of the code via a
transformation analogous to the MacWilliams transform of Theorem~\ref{thM1}.

We first discuss binary codes.
\begin{lemma}\label{GP3}
Let $C$ be a self-orthogonal singly-even binary code, and let $C_0$ be the subset of doubly-even codewords.
Then $C_0$ is a linear subcode of index 2 in $C$.
\end{lemma}
\paragraph{Proof.}
From \eqn{Eq12a}, $\frac{1}{2} wt (u)$ is a linear functional on $C$, and $C_0$ is its kernel.
\begin{defn}\label{dGP4}
{\rm \cite{Me158}.
\index{code!shadow}
The {\em shadow}\footnote{A somewhat more general definition of shadow has been proposed in \cite{BrPl91}, but since it fails to possess the crucial properties (i) and (iii) of Theorem~\ref{GP4} we shall not discuss it here.}\index{shadow} $S$ of a self-orthogonal binary code $C$ is
$$S = \left\{ \begin{array}{ll}
C_0^\perp \setminus C^\perp & \mbox{if~$C$ is singly-even} \\ [+.1in]
C^\perp & \mbox{if $C$ is doubly-even}
\end{array}
\right\}~.
$$
The weight enumerator of the shadow of $C$ will usually be denoted by $S_C ( x,y)$.
\index{$S_C (x,y)$, weight enumerator of shadow}
}
\end{defn}
\paragraph{Examples.}
(i) If $C$ is the repetition code $\{0^n , 1^n \}$ of even length $n$, then if $n \equiv 0$ $(\bmod~4)$, $S= C^\perp =$ all even weight vectors, but if $n \equiv 2$ $(\bmod~4)$, $S=$ all odd weight vectors.
(ii)~If $C = i_2 \oplus i_2 \oplus \cdots \oplus i_2$ then $S$ is the translate of $C$ by
$1010 \ldots 10$.
(iii)~Let $C$ be the $[22,11,6]$ shorter Golay code $g_{22}$,
obtained by ``subtracting''\index{subtraction} (see Section~\ref{GU3}) $i_2$ from $g_{24}$, so that $g_{22}$ consists of all words of $g_{24}$ that begin 00 or 11, with these two coordinates deleted.
Then $S$ consists of the remaining words of $g_{24}$ with the same two coordinates deleted.
\begin{theorem}\label{GP4}
{\rm \cite{Me158}}
The shadow $S$ has the following properties:

(i) $S$ is the set of ``parity vectors'' for $C$; that is,
\beql{EqR12a}
S = \{ u \in \FF_2^n : (u,v) \equiv \frac{1}{2} wt (v) \bmod~2 \quad \mbox{for all} ~~v \in C \}
\eeq

(ii)~$S$ is a coset of $C^\perp$

(iii)
\beql{EqR12b}
S_C (x,y) = \frac{1}{|C|} W_C (x+y, i(x-y)) ~. 
\eeq
\end{theorem}
\paragraph{Proof.}
If $C$ is doubly-even then (i) and (ii) are immediate, and (iii) follows from the
MacWilliams transform and the fact that the weights are divisible by 4.
Suppose $C$ is singly-even, let $C_0$ be the doubly-even subcode, and let $C_1 = C \setminus C_0$.
Then
\beql{EqR13a}
C_0 \subseteq C \subseteq C^\perp \subseteq C_0^\perp ~.
\eeq
The first and last inclusions have index 2, so $C_0^\perp =C^\perp \cup (a+ C^\perp)$, say, where
$(a,u) =0$ for $u \in C_0$, $(a,v) =1$ for $v \in C_1$.
Thus $S = C_0^\perp \setminus C_0 = a + C_0$ has the properties stated in (i) and (ii).
Also,
\begin{eqnarray*}
W_{C_0} (x,y) & = & \frac{1}{2}
\left\{
W_C (x,y) + W_C (x,iy) \right\} ~, \\
W_{C_0^\perp} (x,y) & = &
\frac{1}{|C|} \left\{
W_C (x+y, x-y) + W_C (x+y , i(x-y) \right\} ~,
\end{eqnarray*}
so
$$S_C(x,y) = W_{C_0^\perp} - W_{C^\perp} = \frac{1}{|C|} W_C (x+y, i(x-y)) ~.~~~\bsq$$

If $C$ is a singly-even self-dual code with doubly-even subcode $C_0$, then $C_0^\perp$ is the union of four translates of $C_0$, say $C_0$, $C_1$, $C_2$, $C_3$, with
\beql{EqR14a}
C= C_0 \cup C_2 , \quad S= C_1 \cup C_3 ~.
\eeq

When $n$ is a multiple of 8 then $C' = C_0 \cup C_1$ and $C'' = C_0 \cup C_3$ are both Type II codes (in the notation of Chapter xx (Pless), $C'$
and $C''$ are {\em neighbors}\index{neighbors} of $C$).
If $C$ has a weight 2 word then $C'$ and $C''$ are equivalent.

Similar definitions for the shadow can be given in the other two cases
mentioned.
If $C$ is an additive trace-Hermitian self-orthogonal code over $\FF_4$, let $C_0$ be the subcode with even Hamming weights, and secondly, if $C$ is a self-orthogonal code over $\ZZ_m$ ($m$ even) let $C_0$ be the subcode with Euclidean norms divisible
by $2m$.
In both cases the shadow is defined by:
$$S= \left\{
\begin{array}{ll}
C_0^\perp \setminus C_0 & \mbox{if}~~C \neq C_0 \\ [+.1in]
C^\perp & \mbox{if} ~~C = C_0 ~.
\end{array}
\right.
$$
If $C$ is self-dual from family \Efive\ then the quotient group $C_0^\perp / C_0$ is isomorphic to $Z(2) \times Z(2)$.
If $C$ is self-dual from family \Enine\ then
$C_0^\perp / C_0$ is isomorphic to $Z(2) \times Z(2)$ if $n$ is even and to $Z(4)$ if $n$ is odd.

There are analogues of Theorem~\ref{GP4}.
\begin{theorem}\label{GP5}
Let $C$ be a self-orthogonal additive code over $\FF_4$, with shadow $S$.

(i)~ $S= \{ u \in \FF_4^n : (u,v) \equiv wt(v) ~ ( \bmod~2) \quad \mbox{for all $v \in C$} \}$

(ii) $S$ is a coset of $C^\perp$

(iii) $$S_C (x,y) = \frac{1}{|C|} W_C (x+3y, y-x)~,$$
$$
swe_S (x,y,z) = \frac{1}{|C|} ~ swe_C (x+y+2z, -x-y+2z , y-x )
$$
$$
cwe_S (x,y,z,t) = \frac{1}{|C|} cwe_C ( x+y+z+t, -x-y+z+t, -x+y-z+t, -x+y+z -t) ~.
$$
\end{theorem}
\paragraph{Remark.}
It follows from Theorem~\ref{GP5} that there is a code equivalent to $C$ that has $1^n \in S$.
For the number of vectors of weight $n$ in $S$ is
$$S_C(0,1) = \frac{1}{|C|} W_C (3,1) > 0 ~.$$
All vectors of full weight are equivalent.
\begin{theorem}\label{GP6}
Let $C$ be a self-orthogonal linear code over $\ZZ_4$, with shadow $S$.

(i) $S= \{ u \in \ZZ_4^n : (u,v) \equiv \frac{1}{2} ~{\rm Norm} (v) ~(\bmod~4) \quad \mbox{for all $v \in C$} \}$

(ii) $S$ is a coset of $C^\perp$

(iii) $swe_S (x,y,z) = \frac{1}{|C|} swe_C (x+2y+z, \eta (x-y), -x+2y-z)$,
where $\eta = e^{\pi i/4}$,
$$cwe_S (x,y,z,t) = \frac{1}{|C|} cwe_C (x+y+z+t, \eta (x+iy -z-it), - (x-y+z-t), \eta (x-iy -z+it)) ~.$$
\end{theorem}
\paragraph{Remark.}
It follows that the shadow contains a vector of the form $\pm 1^n$.
(For $cwe_S(0,1,0,1) = \frac{1}{|C|} cwe_C (2,0,2,0) = cwe_C (1,0,1,0) > 0$, since $0^n \in C$.)
This observation, and a formula for the swe equivalent to ours, can be found in \cite{DHS97}.
In particular, a self-dual code from family \EeightII\ contains a vector of the form $\pm 1^n$.
\begin{theorem}\label{GP7}
Let $C$ be a self-orthogonal linear code over $\ZZ_m$, $m$ even, with shadow $S$.

(i) $S = \{ u \in \ZZ_m^n : (u,v) \equiv \frac{1}{2} ~{\rm Norm} (v) ~( \bmod~m) \quad \mbox{for all $v \in C$}\}$

(ii) $S$ is a coset of $C^\perp$

(iii) The cwe of $S$ is obtained from the cwe of $C$ by replacing each $x_j$ by
$$\sum_{k=0}^{m-1} e^{2 \pi i (j^2 + 2jk)/2m} x_k ~,$$
and then dividing by $|C|$.
\end{theorem}

The proofs are analogous to that of Theorem~\ref{GP4}.
\section{Invariant theory}\label{ITheory}
\subsection{An introduction to invariant theory}\label{IIT}\index{invariant theory}
\hsp
If $C$ is self-dual then its weight enumerator must be unchanged by the appropriate
transformation from Theorem~\ref{thM1}.
As we will see, this imposes strong restrictions on the weight enumerator.

We begin by discussing the particular case of the weight enumerator $W(x,y)$ of a binary doubly-even self-dual code $C$.
Since $C$ is self-dual, Theorem~\ref{thM1} implies
\begin{eqnarray}
\label{EqMS7}
W(x,y) & = & \frac{1}{2^{n/2}} W(x+y,x-y) \nonumber \\
& = & W \left( \frac{x+y}{\sqrt{2}} , \frac{x-y}{\sqrt{2}} \right)
\end{eqnarray}
(for $W(x,y)$ is homogeneous of degree $n$).
Since all weights are divisible by 4,
$W(x,y)$ only contains powers of $y^4$.
Therefore
\beql{EqMS8}
W(x,y) = W(x, iy) ~.
\eeq
The problem we wish to solve is to find all polynomials $W(x,y)$ satisfying \eqn{EqMS7} and \eqn{EqMS8}.

\paragraph{Invariants.}
Equation~\eqn{EqMS7} says that $W(x,y)$ is unchanged, or {\em invariant}, under the linear transformation
$$
\begin{array}{ll}
& \mbox{replace $x$ by} ~\displaystyle\frac{x+y}{\sqrt{2}} ~, \\
T_1: \\ 
& \mbox{replace $y$ by} ~ \displaystyle\frac{x-y}{\sqrt{2}} ~,
\end{array}
$$
or, in matrix notation,
$$T_1: \quad \mbox{replace}~
{\binom{x}{y}} ~\mbox{by}~
\frac{1}{\sqrt{2}} \left(
{1 \atop 1} ~ {1 \atop -1} \right)
{\binom{x}{y}} ~.
$$
Similarly, \eqn{EqMS8} says that $W(x,y)$ is also invariant under the transformation
$$
\begin{array}{ll}
& \mbox{replace $x$ by $x$} \\
T_2: \\ 
& \mbox{replace $y$ by $iy$}
\end{array}
$$
or
$$
T_2 : \quad \mbox{replace} ~
{\binom{x}{y}} ~\mbox{by}~ \left( {1 \atop 0} ~ {0 \atop i} \right)
{\binom{x}{y}} ~.
$$

Of course $W(x,y)$ must therefore be invariant under any combination $T_1^2$, $T_1 T_2$, $T_1 T_2 T_1 , \ldots$ of these transformations.
It is not difficult to show (as we shall see in the next section) that the matrices
$$\frac{1}{\sqrt{2}} \left( {1 \atop 1} ~ {1 \atop -1} \right) \quad {\rm and} \quad \left( {1 \atop 0} ~ {0 \atop i} \right)$$
when multiplied together in all possible ways produce a group $\sG_1$ containing 192 matrices.

So our problem now says: find the polynomials $W(x,y)$ which are invariant under all 192 matrices in the group $\sG_1$.

\paragraph{How many invariants?}
The first thing we want to know is how many invariants there are.
This isn't too precise, because of course if $f$ and $g$ are invariants, so is any
constant multiple $cf$ and also $f+g$, $f-g$ and the product $fg$.
Also it is enough to study the {\em homogeneous} invariants (in which all terms have the same degree).

So the right question to ask is:
how many linearly independent, homogeneous invariants are there of each degree $d$?
Let's call this number $a_d$.

A convenient way to handle the numbers $a_0 , a_1, a_2 , \ldots$ is by combining them into a power series or generating function
$$\Phi ( \la ) = a_0 + a_1 \la + a_2 \la^2 + \cdots ~.$$
Conversely, if we know $\Phi ( \la )$, the numbers $a_d$ can be recovered from the power series expansion of $\Phi ( \la )$.

At this point we invoke a beautiful theorem of T. Molien,\index{Molien, T.} published in 1897
\index{Molien theorem}\index{theorem!Molien}
(\cite{MS971};
see also \cite{Bens93}, p.~21;
\cite{MS190}, p.~110;
\cite{MS211}, p.~301;
\cite{MS955}, p.~259;
\cite{Smith95}, p.~86;
\cite{Stur93}, p.~29).
\begin{theorem}\label{ThMS2}
{\rm (Molien)}
For any finite group $\sG$ of complex $m \times m$ matrices,
$\Phi ( \lambda )$ is given by
\beql{EqMS9}
\Phi ( \la ) = \frac{1}{| \sG|} \sum_{A \in \sG}
\frac{1}{\det (I- \la A)} ~.
\eeq
\end{theorem}

We call $\Phi ( \la )$ the {\em Molien series}\index{Molien series} of $\sG$.
\index{series!Molien}
The proof of this theorem is given in the next section.

For our group $\sG_1$, from the matrices corresponding to $I$, $T_1$, $T_2, \ldots$ we get
\beql{EqMS10}
\Phi ( \la ) = \frac{1}{192} \left\{
\frac{1}{(1- \la)^2} +
\frac{1}{1- \la^2} +
\frac{1}{(1- \la) (1- i \la)} + \cdots \right\} ~.
\eeq
There are shortcuts, but it is quite feasible to work out the 192 terms directly (many are the same) and add them.
The result is a surprise:
everything collapses to give
\beql{EqMS11}
\Phi ( \la) = \frac{1}{(1- \la^8) (1- \la^{24})} ~.
\eeq

\paragraph{Interpreting $\Phi ( \la)$.}
The very simple form of \eqn{EqMS11} is trying to tell us something.
Expanding in powers of $\la$, we have
\begin{eqnarray}
\label{EqMS12}
\Phi ( \la ) & = & a_0 + a_1 \la + a_2 \la^2 + \cdots \nonumber \\
& = & (1+ \la^8 + \la^{16} + \la^{24} + \cdots ) ( 1+ \la^{24} + \la^{48} + \cdots ) ~.
\end{eqnarray}

We can deduce one fact immediately: $a_d$ is zero unless $d$ is a multiple of 8,
i.e. the degree of a homogeneous invariant must be a multiple of 8.
(This already proves that the length of
a doubly-even binary self-dual code must be a multiple of 8.)
But we can say more.
The right-hand side of \eqn{EqMS12} is exactly what we would find if
there were two ``basic'' invariants, of degrees 8 and 24, such that all invariants are formed from sums and products of them.

This is because two invariants, $\theta$, of degree 8, and $\phi$, of degree 24, would give rise to the following invariants.
\beql{EqMS13}
\begin{array}{ccc} \hline
\mbox{Degree d} & \mbox{Invariants} & \mbox{Number $a_d$} \\ \hline
0 & 1 & 1 \\
8 & \theta & 1 \\
16 & \theta^2 & 1 \\
24 & \theta^3, \phi & 2 \\
32 & \theta^4, \theta \phi & 2 \\
40 & \theta^5, \theta^2 \phi & 2 \\
48 & \theta^6, \theta^3 \phi, \phi^2 & 3 \\
\cdots & \cdots & \cdots
\end{array}
\eeq
Provided all the products $\theta^i \phi^j$ are linearly independent --- which is the same thing as saying that $\theta$ and $\phi$ are algebraically independent ---
the numbers $a_d$ in \eqn{EqMS13} are exactly the coefficients in
\begin{eqnarray}
\label{EqMS14}
\lefteqn{1+ \la^8 + \la^{16} + 2 \la^{24} + 2 \la^{32} + 2 \la^{40} + 3 \la^{48} + \cdots} \nonumber \\
&&= (1+ \la^8 + \la^{16} + \la^{24} + \cdots ) (1+ \la^{24} + \la^{48} + \cdots ) \nonumber \\
&& = \frac{1}{(1- \la^8) (1- \la^{24})} ~,
\end{eqnarray}
which agrees with \eqn{EqMS11}.
So if we can find two algebraically independent invariants of degrees 8 and 24, we will have solved our problem.
The answer will be that any invariant of this group is a polynomial in $\theta$ and $\phi$.
Now $\phi_8$ (Eq.~\eqn{Eq12}) and $\phi_{24}$ (Eq.~\eqn{Eq13}), the weight enumerators of the Hamming and Golay codes, have degrees 8 and 24 and are invariant under the group.
So we can take $\theta = \phi_8$ and $\phi = \phi_{24}$.
(It's not difficult to verify that they are algebraically independent.)
Actually, it is easier to work with
\beql{EqMS15}
\phi'_{24} = 
\frac{\phi_8^3 - \phi_{24}}{42} = x^4 y^4 (x^4 - y^4)^4
\eeq
rather than $\phi_{24}$ itself.
So we have proved the following theorem, discovered by Gleason in 1970.
\begin{theorem}
\label{ThMS3a}
Any invariant of the group $\sG_1$ is a polynomial in $\phi_8$ and $\phi'_{24}$.
\end{theorem}

This also gives us the solution to our original problem:
\begin{theorem}\label{ThMS3b}
Any polynomial which satisfies Equations $($\ref{EqMS7}$)$ and $($\ref{EqMS8}$)$ is a polynomial in $\phi_8$ and $\phi'_{24}$.
\end{theorem}

Finally, we have characterized the weight enumerator of a doubly-even binary self-dual code.
\begin{theorem}\label{ThMS3c}
{\rm (Gleason \cite{Gle70}.)}\index{Gleason, A. M.}
\index{Gleason theorem}\index{theorem!Gleason}
The weight enumerator of any Type II binary self-dual code is a polynomial in $\phi_8$ and $\phi'_{24}$.
\end{theorem}

Alternative proofs of this astonishing theorem are given by Berlekamp et~al. \cite{Me22},
and
Brou\'{e} and Enguehard \cite{MS201} (see also Assmus and Mattson \cite{MS47}).
But the proof given here seems to be the most informative, and the easiest
to understand and to generalize.

Notice how the exponents 8 and 24 in the denominator of \eqn{EqMS11} led us to guess the degrees of the basic invariants.

This behavior is typical, and is what makes the technique exciting to use.
One starts with a group of matrices $\sG$, computes the complicated-looking sum shown in \eqn{EqMS9}, and simplifies the result.
Everything miraculously collapses, leaving a final expression resembling \eqn{EqMS11} (although not always quite so simple ---
the precise form of the final expression is
given in \eqn{EqMS38a}, \eqn{EqMS38a}).
This expression then tells the degrees of the basic invariants to look for.

\paragraph{Finding the basic invariants.}
In general, finding the basic invariants is a simpler problem than finding $\Phi ( \la )$.
In our applications we can often use
the weight enumerators of codes having the appropriate
properties, as in the above example, or basic invariants can be found by {\em averaging},\index{averaging} using the following simple result (proved in Section~\ref{BTI}).
\begin{theorem}\label{ThMS4}
If $f( \bx ) = f( x_1 , \ldots , x_m )$
is any polynomial in $m$ variables, and $\sG$ is a finite group of $m \times m$ matrices, then
\beql{EqMS15a}
\of ( \bx ) = \frac{1}{|\sG|} \sum_{A \in \sG} A \circ f( \bx )
\eeq
is an invariant, where $A \circ f( \bx ) $ denotes the
polynomial obtained by applying the transformation $A$ to the variables in $f$.
\end{theorem}

Of course $\of ( \bx )$ may be zero.
An example of the use of this theorem is given below.

To illustrate the use of Theorem~\ref{ThMS3c}, we use it to find the weight enumerator of the $[48,24,12]$ extended quadratic residue code
$XQ_{47}$, using only the
fact that it is a doubly-even self-dual code with minimal distance 12.
This implies
that the weight enumerator of the code, which is a homogeneous polynomial of degree 48, has the form
\beql{EqMS17}
W(x,y) = x^{48} + A_{12} x^{36} y^{12} + \cdots ~.
\eeq
The coefficients of $x^{47} y$, $x^{46} y^2, \ldots, x^{37} y^{11}$ are zero.
Here $A_{12}$ is the unknown number of codewords of weight 12.
It is remarkable that, once we know Equation~\eqn{EqMS17}, the weight
enumerator is completely determined by Theorem~\ref{ThMS3c}.
For Theorem~\ref{ThMS3c} says that $W(x,y)$ must be a polynomial in $\phi_8$ and $\phi'_{24}$.
Since $W(x,y)$ is homogeneous of degree 48, $\phi_8$ is homogeneous of degree 8, and $\phi'_{24}$ is homogeneous of degree 24, this polynomial must be a linear combination of $\phi_8^6$, $\phi_8^3 \phi'_{24}$ and $\phi_{24}^{\prime 2}$.

Thus Theorem~\ref{ThMS3c} says that
\beql{EqMS18}
W(x,y) = a_0 \phi_8^6 + a_1 \phi_8^3 \phi'_{24} + a_2 \phi_{24}^{\prime 2} ~,
\eeq
for some real numbers $a_0$, $a_1$, $a_2$.
Expanding \eqn{EqMS18}, we have
\begin{eqnarray}\label{EqMS19}
W(x,y) & = & a_0 (x^{48} + 84x^{44} y^4 + 2946x^{40} y^8 + \cdots )
+ a_1 (x^{44} y^4 + 38x^{40} y^8 + \cdots ) \nonumber \\
&&~~~+ a_2 (x^{40} y^8 - \cdots ) ~,
\end{eqnarray}
and equating coefficients in \eqn{EqMS17}, \eqn{EqMS19} we get
$$a_0 = 1, \quad a_1 = - 84, \quad a_2 = 246 ~.$$
Therefore $W(x,y)$ is uniquely determined.
When these values of $a_0$, $a_1$, $a_2$ are substituted in \eqn{EqMS18} we find that
\begin{eqnarray}\label{EqMs20}
W(x,y) & = & x^{48} + 17296x^{36} y^{12} + 535095x^{32} y^{16} \nonumber \\
&&~~~+ 3995376x^{28} y^{20} + 7681680 x^{24} y^{24} + 3995376x^{20} y^{28} \nonumber \\
&& ~~~+ 535095x^{16} y^{32} + 17296x^{12} y^{36} + y^{48} ~.
\end{eqnarray}
This is certainly faster than computing $W$ by examining each of the $2^{24}$ codewords.

There is a fair amount of algebra involved in computing \eqn{EqMS11}.
Here is a second example, simple enough for the calculations to be shown in full.

For a self-dual code from family \Esix, from \eqn{EqM4} the Hamming weight enumerator
satisfies
\beql{EqMS21}
W \left( \frac{x+ (q-1) y}{\sqrt{q}} , \frac{x-y}{\sqrt{q}} \right) = W(x,y) ~.
\eeq
Let us consider the problem of finding all polynomials which satisfy \eqn{EqMS21}.

The solution proceeds as before.
Equation~\eqn{EqMS21} says that $W(x,y)$ must be invariant under the transformation
$$T_3 : \quad \mbox{replace} ~
{\binom{x}{y}} ~\mbox{by $A$}~
{\binom{x}{y}} ~,
$$
where
\beql{EqMs22}
A= \frac{1}{\sqrt{q}}
\left(
{1 \atop 1} ~ {q-1 \atop -1} \right) ~.
\eeq
Now $A^2 = I$, so $W(x,y)$ must be invariant under the group $\sG_2$ consisting of the two matrices $I$ and $A$.

To find how many invariants there are, we compute the Molien series $\Phi ( \la )$ from \eqn{EqMS9}.
We find
\begin{eqnarray}\label{EqMS23}
\det (I- \la I) & = & (1- \la)^2 ~, \nonumber \\
\det (I- \la A) & = & \det \left[
\begin{array}{cc}
1- \frac{\la}{\sqrt{q}}&- \frac{q-1}{\sqrt{q}} \la \\ [+.1in]
- \frac{\la}{\sqrt{q}}& 1 + \frac{\la}{\sqrt{q}}
\end{array}
\right] = 1- \la^2 ~, \nonumber \\
\Phi ( \la ) & = & \frac{1}{2} \left( \frac{1}{(1- \la)^2} + \frac{1}{1- \la^2} \right) \nonumber \\
& = & \frac{1}{(1- \la) (1- \la^2)} ~.
\end{eqnarray}
which is even simpler than \eqn{EqMS11}.
Equation \eqn{EqMS23} suggests that there might be two basic invariants, of degrees 1 and 2 (the exponents in the denominator).
If algebraically independent invariants of degrees 1 and 2 can be found, say $g$ and $h$, then \eqn{EqMS23} implies that any invariant of $\sG_2$ is a polynomial in $g$ and $h$.

This time we shall use the method of averaging to find the basic invariants.
Let us average $x$ over the group --- i.e., apply
Theorem~\ref{ThMS4} with $f(x,y) =x$.
The matrix $I$ leaves $x$ unchanged, of course, and the matrix $A$ transforms $x$ into $(1/ \sqrt{q}) (x+(q-1)y)$.
Therefore the average,
$$\of(x,y) = \frac{1}{2}
\left[ x + \frac{1}{\sqrt{q}} \{ x+ (q-1) y \} \right] =
\frac{( \sqrt{q} +1 ) \{ x+ ( \sqrt{q} -1) y \}}{2\sqrt{q}} ~,
$$
is an invariant.
Of course any scalar multiple of $\of (x,y)$ is also an invariant, so we may divide by
$( \sqrt{q} +1)/2 \sqrt{q}$ and take
\beql{EqMS24}
g = x+ ( \sqrt{q} -1) y
\eeq
to be the basic invariant of degree 1.
To get an invariant of degree 2 we average $x^2$ over the group, obtaining
$$\frac{1}{2} \left[ x^2 + \frac{1}{q} \{ x+ (q-1) y \}^2 \right] ~.
$$
This can be cleaned up by subtracting $((q+1)/2q)g^2$ (which of course is an invariant), and dividing by a suitable constant.
The result is
$$h= y(x-y) ~,$$
the desired basic invariant of degree 2.

Finally $g$ and $h$ must be shown to be algebraically independent:\index{algebraically independent}
it must be shown that no sum of the form
\beql{EqMS25}
\sum_{i,j} c_{ij} g^i h^j , \quad c_{ij} ~\mbox{complex and not all zero} ~,
\eeq
is identically zero when expanded in powers of $x$ and $y$.
This can be seen by looking at the leading terms.
The leading term of $g$ is $x$, the leading term of $h$ is $xy$, and the leading term of $g^i h^j$ is $x^{i+j} y^j$.
Since distinct summands in \eqn{EqMS25} have distinct leading terms,
\eqn{EqMS25} can only add to zero if all the $c_{ij}$ are zero.
Therefore $g$ and $h$ are algebraically independent.
So we have proved:
\begin{theorem}\label{ThMS5}
Any invariant of the group $\sG_2$, or equivalently any polynomial satisfying
$($\ref{EqMS21}$)$, or equivalently the Hamming weight enumerator of any self-dual code from family \Esix, is a polynomial in $g= x + ( \sqrt{q} -1) y$ and $h=y(x-y)$.
\end{theorem}

At this point the reader should
cry Stop!, and point out that self-dual codes from
family \Esix\ must have even
length, and so every term in the weight enumerator must have even degree.
But in Theorem~\ref{ThMS5}, $g$ has degree 1.

Thus we haven't made use of everything we know about the code.
$W(x,y)$ must also be invariant under the transformation
$$\mbox{replace} ~ {\binom{x}{y}} ~\mbox{by $B$}~
{\binom{x}{y}} ~,
$$where
$$B= \left( {-1 \atop 0} ~ {0 \atop -1} \right) = -I ~.$$
This rules out terms of odd degree.
So $W(x,y)$ is now invariant under the group $\sG_3$ generated by $A$ and $B$, which consists of
$I, ~~A, ~~-I , ~~-A$.
The reader can easily work out that the new Molien series is
\begin{eqnarray}\label{EqMS26}
\Phi_{\sG_3} ( \la) & = &
\frac{1}{2} \{ \Phi_{\sG_2} ( \la) + \Phi_{\sG_2} (- \la) \} \nonumber \\
& = & \frac{1}{2} \left\{
\frac{1}{(1- \la) (1- \la^2)} +
\frac{1}{(1+ \la)(1- \la^2)} \right\} \nonumber \\
& = & \frac{1}{(1- \la^2)^2} ~.
\end{eqnarray}
There are now two basic invariants, both of degree 2 (matching the exponents in the denominator of \eqn{EqMS26}), say $g^2$ and $h$, or the equivalent and slightly simpler pair
$g^\ast = x^2 + (q-1)y^2$ and $h=y (x-y)$.
Hence:
\begin{theorem}\label{ThMS6}
The Hamming weight enumerator of any Hermitian self-dual code over $\FF_q$ 
is a polynomial in $g^\ast$ and $h$.
\end{theorem}

\paragraph{The general plan of attack.}
As these examples have illustrated, there are two stages in using invariant theory to solve a problem.

\paragraph{Stage I.}
Convert the assumptions about the problem (e.g. the code) into algebraic constraints on polynomials (e.g. weight enumerators).
\paragraph{Stage II.}
Use the invariant theory to find all possible polynomials satisfying these constraints.

\subsection{The basic theorems of invariant theory}\label{BTI}
{\bf Groups of matrices.}
Given a collection $A_1, \ldots, A_r$ of $m \times m$ invertible matrices, we can form a group $\sG$ from them by multiplying them together in all possible ways.
Thus $\sG$ contains the matrices $I$, $A_1, A_2, \ldots, A_1 A_2, \ldots, A_2 A_1^{-1} A_2^{-1} A_3, \ldots$.
We say that $\sG$ is {\em generated} by $A_1, \ldots, A_r$.
We will suppose that $\sG$ is finite, which covers all the cases encountered in this chapter.
(For infinite groups, see for example Dieudonn\'{e} and Carroll \cite{MS375},
Rallis \cite{MS1087}, Springer \cite{Sprin77}, Sturmfels \cite{Stur93}, Weyl \cite{MS1410}.)
\paragraph{Example.}
Let us show that the group $\sG_1$ generated by the matrices
$$M= \frac{1}{\sqrt{2}} \left( {1 \atop 1} ~ {1 \atop -1} \right) \quad {\rm and} \quad J = \left( {1 \atop 0} ~ {0 \atop i} \right)
$$that was encountered in Section~\ref{IIT} does indeed have order 192.
The key is to discover (by randomly multiplying matrices together) that
$\sG_1$ contains
$$
\begin{array}{rllrll}
J^2 & = & \left( \matrix{
1 & 0 \cr
0 & -1 \cr
}
\right) ~, & ~~~ E & = &
(MJ)^3 = \frac{1+i}{\sqrt{2}}
\left( \matrix{
1& 0 \cr
0 & 1 \cr
}
\right) ~, \\ [+.2in]
E^2 & = & i \left(
\matrix{ 1& 0 \cr
0 & 1 \cr
}
\right)~, & R & = & MJ^2 M = \left( \matrix{
0& 1 \cr
1 & 0 \cr
} \right) ~.
\end{array}
$$
So $\sG_1$ contains the matrices
$$\af \left( \matrix{
1 & 0 \cr
0 & \pm 1 \cr
}
\right) ~, \quad
\af \left( \matrix{
0 & 1 \cr
\pm 1 & 0 \cr
}
\right)~, \quad \af \in \{1, i, -1, -i \} ~,
$$
which form a subgroup $\sH_1$ of order 16.
From this it is easy to see that $\sG_1$ consists of the union of 12 cosets of $\sH_1$:
\beql{EqMS27}
\sG_1 = \bigcup_{k=1}^{12} a_k \sH_1 ~,
\eeq
where $a_1, \ldots, a_6$ are respectively
$$
\left( \matrix{
1 & 0 \cr
0 & 1 \cr
} \right) , ~~\left( \matrix{
1& 0 \cr
0 & i \cr
} \right) , 
\frac{1}{\sqrt{2}} \left( \matrix{1 & 1 \cr
1 & -1 \cr
} \right) , ~
\frac{1}{\sqrt{2}}
\left( \matrix{
1 & 1 \cr
i & -i \cr
} \right) , ~
\frac{1}{\sqrt{2}} \left( \matrix{
1& i \cr
1 & -i \cr
} \right) , ~ \frac{1}{\sqrt{2}}
\left( \matrix{
1 & i \cr
i & 1 \cr
} \right),
$$
$a_7 = \eta a_1, \ldots, a_{12} = \eta a_6$, and $\eta = (1+i)/ \sqrt{2}$, an 8th root of unity.
Thus $\sG_1$ consists of the 192 matrices
\beql{EqMS28}
\eta^\nu \left( \matrix{
1 & 0 \cr
0 & \af \cr
} \right) , ~ \eta^\nu \left( \matrix{
0 & 1 \cr
\af & 0 \cr
} \right) , ~ \eta^\nu \frac{1}{\sqrt{2}} \left( \matrix{
1 & \beta \cr
\af &  - \af \beta \cr
} \right) ,
\eeq
for $0 \le \nu \le 7$ and $\af, \beta \in \{1, i, -1, -i \}$.

As a check, one verifies that every matrix in \eqn{EqMS28} can be written as a product of $M$'s and $J$'s;
that the product of two matrices in \eqn{EqMS28} is again in \eqn{EqMS28};
and that the inverse of every matrix in \eqn{EqMS28} is in \eqn{EqMS28}.
Therefore \eqn{EqMS28} {\em is} a group, and is the group generated by $M$ and $J$.
Thus $\sG_1$ is indeed equal to \eqn{EqMS28}.

We have gone into this example in some detail to emphasize that it is important to begin by understanding the group thoroughly.
(For an alternative way of studying $\sG_1$, see \cite[pp.~160--161]{MS201}.

\paragraph{Invariants.}
To quote Hermann Weyl \cite{MS1409},\index{Weyl, H.}
``the theory of invariants came into existence about the middle of the nineteenth century somewhat like
Minerva: a grown-up virgin, mailed in the shining armor of algebra, she
sprang forth from Cayley's\index{Cayley, A.} Jovian head.''
Invariant theory became one of the main branches of nineteenth century mathematics, but dropped out of fashion after
Hilbert's\index{Hilbert} work:
see Fisher \cite{MS431} and Reid \cite{MS1108}.
In the past thirty years, however, there has been a resurgence of interest,
with applications in algebraic geometry (Dieudonn\'{e} and Carroll \cite{MS375},
Mumford and Fogarty \cite{MS976}),
physics (see for example Agrawala and Belinfante \cite{MS7} and the references given there), combinatorics (Doubilet et~al. \cite{MS381}, Rota \cite{MS1125}, Stanley \cite{Stan79}) and
coding theory (\cite{Me27}, \cite{Me29}, \cite{Me35}, \cite{Me80}).
Recently a number of monographs
(Benson \cite{Bens93}, Bruns and Herzog \cite{BH93},
Smith \cite{Smith95},
Springer \cite{Sprin77}, Sturmfels \cite{Stur93}) and conference proceedings
(\cite{Foss86}, \cite{Gher82}, \cite{Koh87}, \cite{Stant90}) on invariant theory have appeared.

There are several different kinds of invariants,
but here an invariant is defined as follows.

Let $\sG$ be a group of $g$ $m \times m$ complex matrices $A_1, \ldots, A_g$, where the $(i,k)^{\rm th}$ entry of $A_\af$ is $a_{ik}^{( \af )}$.
In other words $\sG$ is a group of linear transformations on the variables
$x_1, \ldots, x_m$, consisting of the transformations
\beql{EqMS29}
T^{( \af)}:
\mbox{replace $x_i$ by $x_i^{( \af)} = \displaystyle\sum_{k=1}^m a_{ik}^{(\af)} x_k$},
\quad i=1, \ldots, m
\eeq
for $\af = 1,2, \ldots, g$.
It is worthwhile giving a careful description of how a polynomial $f( \bx ) = f(x_1, \ldots, x_m )$ is transformed by a matrix $A_\af$ in $\sG$.
The transformed polynomial is
$$A_\af \circ f( \bx ) = f(x_1^{(\af)} , \ldots, x_m^{(\af)})
$$
where each $x_i^{(\af)}$ is replaced by $\sum_{k=1}^m a_{ik}^{(\af)} x_k$.
Another way of describing this is to think of
$\bx = (x_1, \ldots, x_m)^T$ as a column vector.
Then $f( \bx )$ is transformed into
\beql{EqMS30}
A_\af \circ f( \bx ) = f( A_\af \bx ) ~,
\eeq
where $A_\af \bx$ is the usual product of a matrix and a vector.
One can check that
\beql{EqMS31}
B \circ (A \circ f( \bx )) = (AB) \circ f( \bx ) = f( AB \bx ) ~.
\eeq
For example,
$$A = \left( \matrix{
1 & 2 \cr
0 & -1 \cr
} \right)
$$
transforms $x_1^2 + x_2$ into $(x_1 +2x_2)^2 - x_2$.
\paragraph{Definition.}
An {\em invariant}\index{invariant} of $\sG$ is a polynomial $f( \bx )$ which is unchanged by
every linear transformation in $\sG$.
In other words, $f( \bx )$ is an invariant of $\sG$ if
\beql{EqMS31a}
A_\af \circ f( \bx ) = f( A_\af \bx ) = f( \bx )
\eeq
for all $\af = 1, \ldots , g$.
\paragraph{Example.}
Let
$$\sG_4 = \left\{ \left( \matrix{
1 & 0 \cr
0 & 1 \cr
} \right) , ~\left( \matrix{
-1 & 0 \cr
0 & -1 \cr
} \right) \right\} ~,
$$
a group of order $g=2$.
Then $x^2$, $xy$ and $y^2$ are homogeneous invariants of degree 2.

Even if $f(x)$ isn't an invariant, its average over the group,
\beql{EqMS32}
\of ( \bx ) = \frac{1}{g} \sum_{\af=1}^g A_\af \circ f( \bx )
\eeq
is, as was already stated in Theorem~\ref{ThMS4}.
To prove this, observe that any
$A_\beta \in \sG$ transforms the right-hand side of \eqn{EqMS32} into
\beql{EqMS33}
\frac{1}{g} \sum_{\af=1}^g (A_\af A_\beta ) \circ f( \bx )~,
\eeq
by \eqn{EqMS31}.
As $A_\af$ runs through $\sG$, so does $A_\af A_\beta$, if $A_\beta$ is fixed.
Therefore \eqn{EqMS33} is equal to
$$\frac{1}{g} \sum_{\gamma =1}^g A_\gamma \circ f( \bx ) ~, $$
which is $\of ( \bx )$.
Therefore $\of ( \bx )$ is an invariant.~~~$\bsq$

More generally, any symmetric function of the $g$ polynomials $A_1 \circ f( \bx ) , \ldots, A_g \circ f( \bx )$ is an invariant of $\sG$.

Clearly, if $f( \bx )$ and $h( \bx )$ are invariants of $\sG$, so are
$f( \bx ) + h( \bx )$, $f( \bx ) h ( \bx )$, and
$cf ( \bx )$ ($c$ complex);
or in other words the set of invariants of $\sG$,
which we denote by $\sJ ( \sG )$, forms a ring.
\index{ring!of invariants}

One of the main problems of invariant theory is to describe $\sJ ( \sG )$.
Since the transformations in $\sG$ do not change the degree of a polynomial, it is enough to describe the
homogeneous invariants (for any invariant is a sum of homogeneous invariants).

\paragraph{Basic invariants.}
Our goal is to find a ``basis'' for the invariants of $\sG$,
that is, a set of basic invariants such that any invariant can be expressed in terms of this set.
There are two different types of bases one might look for
\paragraph{Definition.}
Polynomials $f_1 ( \bx ) , \ldots, f_r ( \bx )$ are called
{\em algebraically dependent}\index{algebraically dependent} if there is a polynomial $p$ in $r$ variables
with complex coefficients, not all zero, such that
$p(f_1 ( \bx ) , \ldots, f_r ( \bx ))$ is identically zero.
Otherwise $f_1 ( \bx ) , \ldots, f_r ( \bx )$ are 
{\em algebraically independent}.\index{algebraically independent}
A fundamental result from algebra is (Jacobson \cite{MS687}, vol.~3, p.~154):
\begin{theorem}\label{ThMS7}
Any $m+1$ polynomials in $m$ variables are algebraically dependent.
\end{theorem}

The first type of basis we might look for is a set of $m$
algebraically independent invariants $f_1 ( \bx ) , \ldots f_m ( \bx )$.
Such a set is indeed a ``basis,'' for by Theorem~\ref{ThMS7} any invariant is
algebraically dependent on $f_1, \ldots, f_m$ and so is a root of a polynomial
equation in $f_1, \ldots, f_m$.
The following theorem guarantees the existence of such a basis.
\begin{theorem}\label{ThMS8}
{\rm \cite[p. 357]{MS211}}
There always exist $m$ algebraically independent invariants of $\sG$.
\end{theorem}
\paragraph{Proof.}
Consider the polynomial
$$\prod_{\af=1}^g (t- A_\af \circ x_i )$$
in the variables $t$, $x_1, \ldots, x_m$.
Since one of the $A_\af$ is the identity matrix, $t=x_1$ is a zero of this
polynomial.
When the polynomial is expanded in powers of $t$,
the coefficients are invariants, by the remark immediately following the proof
of Theorem~\ref{ThMS4}.
Therefore $x_1$ is an algebraic function of invariants.
Similarly each of $x_2 , \ldots, x_m$ is an algebraic function of invariants.
Now if the number of algebraically independent invariants were $m'$ $(< m)$, the $m$ independent variables
$x_1, \ldots, x_m$ would be algebraic functions of the $m'$ invariants, a
contradiction.
Therefore the number of algebraically independent invariants is at least $m$.
But by Theorem~\ref{ThMS7} this number cannot be greater than $m$.~~~$\bsq$
\paragraph{Example.}
For the preceding group $\sG_4$, we may take $f_1 = (x+y)^2$ and $f_2 = (x-y)^2$ as the algebraically independent invariants.
Then any invariant is a root of a polynomial equation in $f_1$ and $f_2$.
For example,
$$
\begin{array}{rll}
x^2 & = & \frac{1}{4} \left( \sqrt{f_1} + \sqrt{f_2} \right)^2 ~, \\ [+.1in]
xy & = & \frac{1}{4} (f_1 - f_2 ) ~,
\end{array}
$$
and so on.

However, by far the most convenient description of the invariants is a set $f_1, \ldots, f_l$ of invariants with the property that any invariant is a
{\em polynomial} in $f_1, \ldots, f_l$.
Then $f_1, \ldots, f_l$ is called a {\em polynomial basis}\index{polynomial basis} (or an
{\em integrity basis})\index{integrity basis} for the invariants of $\sG$.
Of course if $l > m$ then by Theorem~\ref{ThMS7} there will be polynomial equations, called {\em syzygies}, relating $f_1, \ldots, f_l$.
\index{syzygy}

For example, $f_1 = x^2$, $f_2 = xy$, $f_3 = y^2$ form a polynomial basis for the
invariants of $\sG_4$.
The syzygy relating them is
$$f_1 f_2 - f_2^2 =0 ~.$$
The existence of a polynomial basis, and a method of finding it, is given by the next theorem.
\begin{theorem}\label{ThMS9}
{\rm (Noether\index{Noether, E.} \cite{MS997}; \cite[p.~275]{MS1410}.)}
\index{Noether theorem}\index{theorem!Noether}
The ring of invariants of a finite group $\sG$ of complex $m \times m$ matrices has a polynomial basis consisting of not more than ${\binom{m+g}{m}}$ invariants, of degree not exceeding $g$, where $g$ is the order of $\sG$.
Furthermore this basis may be obtained by taking the average over $\sG$ of all monomials
$$x_1^{b_1} x_2^{b_2} \cdots x_m^{b_m}$$
of total degree $\sum b_i$ not exceeding $g$.
\end{theorem}
\paragraph{Proof.}
Let the group $\sG$ consist of the transformations \eqn{EqMS29}.
Suppose
$$f(x_1, \ldots, x_m) = \sum_{e} c_e x_1^{e_1} \cdots x_m^{e_m} ~,$$
$c_e$ complex, is any invariant of $\sG$.
(The sum extends over all $e= e_1 \cdots e_m$ for which there is nonzero term $x_1^{e_1} \cdots x_m^{e_m}$ in $f(x_1, \ldots, x_m )$.)
Since $f(x_1, \ldots, x_m )$ is an invariant, it is unchanged when we
average it over the group, so
\begin{eqnarray*}
f(x_1, \ldots, x_m) & = &
\frac{1}{g} \{ f( x_1^{(1)} , \ldots, x_m^{(1)} ) + \cdots +
f(x_1^{(g)} , \ldots, x_m^{(g)}) \} \\
& = & \frac{1}{g} \sum_e c_e \{ (x_1^{(1)})^{e_1} \cdots (x_m^{(1)} )^{e_m} + \cdots +
(x_1^{(g)} )^{e_1} \cdots (x_m^{(g)})^{e_m} \} \\
& = & \frac{1}{g} \sum_e c_e J_e \quad {\rm (say)} ~.
\end{eqnarray*}
Every invariant is therefore a linear combination of the (infinitely many)
special invariants
$$J_e = \sum_{\alpha =1}^g (x_1^{(\af )} )^{e_1} \cdots (x_m^{(\af)})^{e_m} ~.$$
Now $J_e$ is (apart from a constant factor) the coefficient of
$u_1^{e_1} \cdots u_m^{e_m}$ in
\beql{EqMS34}
P_e = \sum_{\af =1}^g (u_1 x_1^{(\af)} + \cdots + u_m x_m^{(\af)} )^{e}~,
\eeq
where $e= e_1 + \cdots + e_m$.
In other words, the $P_e$ are the power sums of the $g$ quantities
$$u_1 x_1^{(1)} + \cdots + u_m x_m^{(1)} , \ldots , u_1x_1^{(g)} + \cdots + u_m x_m^{(g)} ~.
$$
Any power sum $P_e$, $e=1,2, \ldots$, can be written as a polynomial with rational coefficients in the first $g$ power sums
$P_1$, $P_2 , \ldots, P_g$.
Therefore any $J_e$ for
$$e = \sum_{i=1}^m e_i > g$$
(which is a coefficient of $P_e$) can be written as a polynomial in the
special invariants
$$J_e \quad {\rm with} \quad e_1 + \cdots + e_m \le g$$
(which are the coefficients of $P_1, \ldots, P_g$).
Thus any invariant can be written as a polynomial in those $J_e$ with
$\sum_{i=1}^m e_i \le g$.
The number of such $J_e$ is the number of $e_1$, $e_2, \ldots, e_m$ with
$e_i \ge 0$ and $e_1 + \cdots + e_m \le g$, which is ${\binom{m+g}{m}}$.
Finally, $\deg J_e \le g$, and $J_e$ is obtained by averaging $x_1^{e_1} \cdots x_m^{e_m}$ over the group.~~~$\bsq$

\paragraph{Molien's theorem.}\index{Molien theorem}\index{theorem!Molien}
Since we know from Theorem~\ref{ThMS9} that a polynomial basis
always exists, we can go ahead with confidence and try to find it, using the methods described in Section~\ref{IIT}.
To discover when a basis has been found, we use Molien's theorem (Theorem~\ref{ThMS2}).
This states that if $a_d$ is the number of linearly independent homogeneous invariants of $\sG$ with degree $d$, and
$$\Phi_\sG ( \la ) = \sum_{d=0}^\infty a_d \la^d ~,
$$
then
\beql{EqMS35}
\Phi_\sG ( \la) = \frac{1}{g}
\sum_{\af =1}^g \frac{1}{\det (I - \la A_\af )} ~.
\eeq
The proof depends on the following theorem.
\begin{theorem}\label{ThMS10}
{\rm \cite[p.~258]{MS955}, \cite[p.~17]{MS1185}}
The number of linearly independent invariants of $\sG$ of degree 1 is
$$a_ 1 = \frac{1}{g} \sum_{\af =1}^g \mbox{{\em trace}} (A_\af ) ~.$$
\end{theorem}
\paragraph{Proof.}
Let
$$S = \frac{1}{g} \sum_{\af=1}^g A_\af ~. $$
Changing the variables on which $\sG$ acts from $x_1, \ldots, x_m$ to $y_1, \ldots, y_m$, where
$(y_1, \ldots, y_m) = (x_1, \ldots, x_m) T^{tr}$, changes $S$ to $S' = TST^{-1}$.
We may choose $T$ so that $S'$ is diagonal (see \cite[p.~252]{MS211}).
Now $S^2 = S$, $(S')^2 = S'$, hence the diagonal entries of $S'$ are 0 or 1.
So with a change of variables we may assume
$$S = \left[
\begin{array}{cccccc}
1 & ~ & ~ & ~ & ~ & 0 \\
~ & \ddots & \\
~ & ~ & 1 \\
~ & ~ & ~ & 0 \\
~ & ~ & ~ & ~ & \ddots \\
0 & ~ & ~ & ~ & ~ & 0
\end{array}
\right] 
$$
with say $r$ 1's on the diagonal.
Thus $S \circ y_i = y_i$ if $1 \le i \le r$, $S \circ y_i =0$ if $r+1 \le i \le m$.

Any linear invariant of $\sG$ is certainly fixed by $S$, so $a_1 \le r$.
On the other hand, by Theorem~\ref{ThMS4},
$$S \circ y_i = \frac{1}{g} \sum_{\af =1}^g A_\af \circ y_i$$
is an invariant of $\sG$ for any $i$, and so $a_1 \ge r$.~~~$\bsq$

Before proving Theorem~\ref{ThMS2} let us introduce some more notation.
Equation~\eqn{EqMS29} describes how $A_\af$ transforms the variables
$x_1, \ldots, x_m$.
The $d^{\rm th}$ {\em induced matrix},\index{induced matrix} denoted by
$A_\af^{[d]}$, describes how $A_\af$ transforms the products of the $x_i$ taken $d$ at a time,
namely $x_1^d , x_2^d , \ldots, x_1^{d-1} x_2 , \ldots$ (Littlewood \cite[p.~122]{MS857}).
E.g.
$$A_\af = {\binom{a~b}{c~d}}$$
transforms $x_1^2$, $x_1 x_2$ and $x_2^2$ into
\begin{eqnarray*}
&& a^2 x_1^2 + 2 abx_1 x_2 + b^2 x_2^2 ~, \\
&& acx_1^2 + (ad + bc) x_1 x_2 + bdx_2^2 ~, \\
&& c^2 x_1^2 + 2cdx_1 x_2 + d^2 x_2^2
\end{eqnarray*}
respectively.
Thus the $2^{\rm nd}$ induced matrix is
$$A_\af^{[2]} = \left[
\matrix{a^2 & 2ab & b^2 \cr
ac & ad+bc & bd \cr
c^2 & 2cd & d^2 \cr
} \right] ~,
$$
\paragraph{Proof of Theorem~\ref{ThMS2}.}
To prove \eqn{EqMS35}, note that $a_d$ is equal to the number of
linearly independent invariants of degree 1 of
$\sG^{[d]} = \{ A_\af^{[d]}: \af = 1, \ldots, g \}$.
By Theorem~\ref{ThMS10},
$$a_d = \frac{1}{g} \sum_{\af =1}^g {\rm trace} ~ A_\af^{[d]} ~.$$
Therefore, to prove Theorem~\ref{ThMS2}, it is enough to show that the
trace of $A_\af^{[d]}$ is equal to the coefficient of $\la^d$ in
\beql{EqMS36}
\frac{1}{\det (I- \la A_\af)} =
\frac{1}{(1- \la \om_1) \cdots (1- \la \om_m)} ~,
\eeq
where $\om_1 , \ldots, \om_m$ are the eigenvalues of $A_\af$.
By a suitable change of variables we can make
$$A_\af = \left[
\matrix{
\om_1 & ~ & 0 \cr
~ & \ddots & ~ \cr
0 & ~ & \om_m \cr
} \right] , \quad
A_\af^{[d]} =
\left[ \matrix{
\om_1^d & ~ & ~ & 0 & ~ \cr
~ & \om_2^d & ~ & ~ & ~ \cr
~ & ~ & \ddots \cr
~ & ~ & ~ & \om_1^{d-1} \om_2 \cr
~ & 0 & ~ & ~ & \ddots \cr
} \right] ~,
$$
and trace $A_\af^{[d]} =$ sum of the products of $\om_1, \ldots, \om_m$ taken $d$ at a time.
But this is exactly the coefficient of $\la^d$ in the expansion of \eqn{EqMS36}.~~~$\bsq$

It is worth remarking that the Molien series does not determine the group.
For example there are two groups of $2 \times 2$ matrices of order 8 having
$$\Phi ( \la ) = \frac{1}{(1- \la^2) (1- \la^4)}$$
(namely the dihedral group $D_8$ and the abelian group $Z(2) \times Z(4)$).
In fact there exist abstract groups $\sA$ and $\sB$ whose matrix representations can be paired in such a way that every
representation of $\sA$ has the same Molien series as the corresponding
representation of $\sB$ (Dade \cite{MS324}).

\paragraph{A standard form for the basic invariants.}
The following notation is very  useful in describing the ring
$\sJ ( \sG )$ of invariants of a group $\sG$.
The complex numbers are denoted by $\CC$, and if $p(\bx )$, $q(\bx) , \ldots$ are polynomials, $\CC [p( \bx ) , q( \bx ) , \ldots ]$ denotes the set of all polynomials in $p( \bx ) , q( \bx )$, $\ldots$ with complex coefficients.
For example Theorem~\ref{ThMS3a} just says that $\sJ( \sG_1 ) = \CC [\phi_8 , \phi'_{24} ]$.

Also, $\oplus$ will denote the usual direct sum operation.
For example a statement like $\sJ ( \sG) = R \oplus S$ means that every invariant of $\sG$ can be written uniquely in the form $r+s$ where $r \in R$, $s \in S$.

Using this notation we can now specify the most convenient form of polynomial basis for $\sJ ( \sG)$.

\paragraph{Definition.}
{\rm A {\em good polynomial basis}\index{good polynomial basis} for $\sJ ( \sG )$ consists of homogeneous invariants $f_1, \ldots, f_l$ $(l \ge m)$
where $f_1 , \ldots, f_m$ are algebraically independent and
\beql{EqMS38a}
\sJ ( \sG ) = \CC [f_1 , \ldots, f_m ] \quad {\rm if} \quad l=m ~,
\eeq
or, if $l> m$,
\beql{EqMS38b}
\sJ ( \sG) = \CC [f_1, \ldots, f_m ] \oplus f_{m+1} \CC [f_1, \ldots, f_m ] \oplus \cdots \oplus
f_l \CC [f_1, \ldots, f_m ] ~.
\eeq
}

In words, this says that any invariant of $\sG$ can be written as a polynomial in $f_1 , \ldots, f_m$ (if $l=m$), or as such a polynomial plus $f_{m+1}$ times another such polynomial plus $\cdots$ (if $l>m$).
$f_1 , \ldots, f_m$ are called {\em primary} invariants\index{primary invariants} and
$f_{m+1}, \ldots, f_l$ (if present) are {\em secondary} invariants.\index{secondary invariants}
Speaking loosely, \eqn{EqMS38a} and \eqn{EqMS38b} say
that when describing an arbitrary invariant,
$f_1, \ldots, f_m$ are ``free'' and can be used as often as needed,
while $f_{m+1} , \ldots, f_l$ are ``transients'' and can each be used at most once.
Equations \eqn{EqMS38a} and \eqn{EqMS38b} are sometimes called a
{\em Hironaka decomposition}\index{Hironaka decomposition} of $\sJ ( \sG )$ (\cite{Stur93}, p.~39).

For a good polynomial basis $f_1 , \ldots, f_l$ we can say exactly what
form the
syzygies must take.
If $l=m$ there are no syzygies.
If $l>m$ there are $\binom{l-m+1}{2}$ syzygies expressing the products
$f_i f_j$ $(m+1 \le i \le j \le l)$ in terms of $f_1, \ldots, f_l$.

It is important to note that the Molien series can be written down by
inspection from the degrees of a good polynomial basis.
Let $d_1 = \deg f_1, \ldots, d_l = \deg f_l$.
Then
\beql{EqMS39a}
\Phi_\sG ( \la ) = \frac{1}{\prod_{i=1}^m (1- \la^{d_i})} , \quad {\rm if} \quad l=m ~,
\eeq
or
\beql{EqMS39b}
\Phi_\sG ( \la ) =
\frac{1+ \sum_{j=l+1}^m \la^{d_j}}{\prod_{i=1}^m (1- \la^{d_i})} , \quad {\rm if} \quad l > m ~.
\eeq
(This is easily verified by expanding
\eqn{EqMS39a} and \eqn{EqMS39b}
in powers of $\la$ and comparing
with \eqn{EqMS38a} and \eqn{EqMS38b}.)

Some examples will make this clear.

(1) For the group $\sG_1$ of Section~\ref{IIT}, $f_1 = \phi_8$ and $f_2 = \phi'_{24}$ form a good polynomial basis, with degrees $d_1 =8$, $d_2 =24$.
Indeed, from Theorem~\ref{ThMS3a} and \eqn{EqMS11},
$$\sJ ( \sG_1) = \CC [\phi_8 , \phi'_{24} ]$$
and
$$\Phi_{\sG_1} ( \la) = \frac{1}{(1- \la^8) (1- \la^{24})} ~.$$

(2) For the group $\sG_4$ defined above, $f_1 = x^2$, $f_2 = y^2$, $f_3 = xy$ is a good polynomial basis, with $d_1 = d_2 = d_3 =2$.
The invariants can be described as
\beql{EqMS39}
\sJ ( \sG_4 ) = \CC [x^2, y^2] \oplus xy \CC [x^2, y^2] ~.
\eeq
In words, any invariant can be written uniquely as a polynomial in $x^2$ and $y^2$ plus $xy$ times another such polynomial.
E.g.
$$(x+y)^4 = (x^2)^2 + 6x^2 y^2 + (y^2)^2 + xy(4x^2 + 4y^2) ~.$$
The Molien series is
$$\Phi_{\sG_4} ( \la) =
\frac{1}{2} \left\{
\frac{1}{(1- \la)^2} + \frac{1}{(1+ \la)^2} \right\} =
\frac{1+ \la^2}{(1- \la^2)^2}
$$
in agreement with \eqn{EqMS39b} and \eqn{EqMS39}.
The single syzygy is $x^2 \cdot y^2 = (xy)^2$.
Note that $f_1 = x^2$, $f_2 = xy$, $f_3 =y^2$ is not a good polynomial basis, for the invariant $y^4$ is not in the ring $\CC [x^2, xy] \oplus y^2 \CC [x^2, xy]$.

Fortunately the following result holds.
\begin{theorem}\label{ThMS11}
{\rm (Hochster and Eagon \cite[Proposition 13]{MS657})}
A good polynomial basis exists for the invariants of any finite group
of complex $m \times m$ matrices.
\end{theorem}

For the proof see \cite{Bens93}, \cite{BH93}, \cite{MS657} or
\cite{Smith95}.

So we know that for any group the Molien series can be put into the standard form of \eqn{EqMS39a}, \eqn{EqMS39b} (with denominator consisting
of a product of $m$ factors
$(1- \la^{d_i} )$ and numerator consisting of sum of powers of $\la$ with positive coefficients);
and that a good polynomial basis \eqn{EqMS38a}, \eqn{EqMS38b} can be found whose degrees match the powers of $\la$
occurring in the standard form of the Molien series.

On the other hand the converse is not true.
It is not always true that when the Molien series has been put into the form
\eqn{EqMS39a}, \eqn{EqMS39b}
(by cancelling common factors and multiplying top and bottom
by new factors), then a good polynomial basis for $\sJ ( \sG )$ can be found whose
degrees match the powers of $\la$ in $\Phi ( \la )$.
This is shown by the following example, due to Stanley \cite{MS1262}.

Let $\sG_6$ be the group of order 8 generated by the matrices
${\rm diag} \{-1, -1, -1\}$ and
${\rm diag} \{1,1,i\}$.
The Molien series is
\begin{eqnarray}
\label{EqMS40}
\Phi_{\sG_6} ( \la ) & = & \frac{1}{(1- \la^2)^3} \\
\label{EqMS41}
& = & \frac{1+ \la^2}{(1- \la^2)^2(1- \la^4)} ~.
\end{eqnarray}
A good polynomial basis exists corresponding to \eqn{EqMS41},
namely
$$\sJ( \sG_6 ) = \CC [x^2 , y^2, z^4 ] \oplus xy \CC[x^2, y^2, z^4 ] ~,$$
but there is no good polynomial basis corresponding to \eqn{EqMS40}.

\paragraph{Remarks.}
(1)~Shephard and Todd \cite{MS1196}\index{Shephard and Todd classification} have characterized those groups for which
\eqn{EqMS38a} holds, i.e. for which a good polynomial basis exists consisting only of
algebraically independent invariants.
These are the groups known as ``unitary groups generated by reflections.''\index{unitary groups generated by reflections}
\index{reflection group}
\index{group!reflection}
A complete list of the 37 irreducible groups (or families of groups) of
this type is given in \cite{MS1196} and \cite{Smith95}, p.~199.

(2)~Sturmfels\index{Sturmfels, B.} \cite{Stur93} gives an algorithm for computing a good
polynomial basis for the ring of invariants of a finite group.
The computer language MAGMA\index{MAGMA} (\cite{Mag1}, \cite{Mag2}, \cite{Mag3}) has
commands for computing Molien series and finding a good polynomial basis (and many other things).

(3)~{\bf Relative invariants.}
If $\chi$ is a homomorphism from $\sG$ into the multiplicative group of the complex numbers (i.e. a {\em linear character}\index{linear character} of $\sG$), then a polynomial $f( \bx )$ is called a {\em relative invariant}\index{relative invariant} of $\sG$ with respect to $\chi$ if
\index{invariant!relative}
$$A \circ f( \bx ) = \chi(A) f( \bx ) \quad \mbox{for all} \quad A \in \sG~.$$
Molien's theorem for relative invariants states that the number of linearly
independent homogeneous relative invariants with respect to $\chi$ of degree $\nu$ is the coefficient of $\la^\nu$ in the expansion of
\index{Molien series}
$$\frac{1}{|\sG|}
\sum_{A \in \sG}
\frac{\overline{\chi} (A)}{\det | I - \la A|} ~.
$$

\section{Gleason's theorem and generalizations}\label{GL}\index{Gleason theorem}\index{theorem!Gleason}
\hsp
We now make use of the machinery developed in the previous section to give a series of results that characterize the rings to which the various weight enumerators of self-dual codes belong.
The first theorems of this type, for binary and ternary codes, were discovered by Gleason \cite{Gle70}.
The results can be proved by the generalizations of the arguments used to establish
Theorem~\ref{ThMS3c}.
We remind the reader that $hwe$, $swe$ and $cwe$ stand for Hamming,
symmetrized and complete weight enumerators, respectively.
The code under consideration is denoted by $C$ and its shadow by $S$.

In each case the conclusion is that the weight enumerator being considered must be an element of a certain ring $R$.
We describe $R$ by giving its {\em Molien series}\index{Molien series} (also called a {\em Hilbert series}\index{Hilbert series}\index{series!Hilbert} or {\em Poincar\'{e} series})\index{Poincar\'{e} series}\index{series!Poincar\'{e}}
$$\Phi ( \la ) = \sum_{n=0}^\infty ( \dim_\CC R_n) \la^n ~,$$
where $R_n$ is the subspace of homogeneous polynomials in $R$ of degree $n$.
We then give a good polynomial basis for $R$ 
(in the sense of \eqn{EqMS38a}, \eqn{EqMS38a}).

In many cases $R$ is obtained (as described in the previous sections) as the ring of invariants of a certain matrix group $G$.
If so then we start by giving generators for $G$, its order, and, if it is a well-known group, a brief description.
We have preferred to give natural generators for $G$, rather than attempting to find a minimal but less-intuitive set --- in most cases two generators would suffice.
If $G$ is a reflection group we give its number in Shephard and Todd's list (\cite{MS1196}; \cite{Smith95}, page~199).

In other cases (the symmetrized weight enumerator of a Hermitian self-dual code over $\FF_4$,
\eqn{EqG19a}, for example) the ring $R$ cannot be found directly as the ring of invariants of any group, but must be obtained by collapsing the ring of complete weight enumerators.

At the end of each subsection is a table that gives, for most of the rings mentioned, a list of
codes whose weight enumerators provide a polynomial basis for the ring.
The weight enumerators of the codes before the semicolon are primary invariants, those after the semicolon (if present) are secondary invariants.

For example, the first line of Table~\eqn{Eq5c} is equivalent to Theorem~\ref{ThMS3c}.
\subsection{Family $\EoneI$: Binary self-dual codes} \label{GL1a}
\paragraph{hwe of code $C$.}
{\rm (\cite{Gle70}, \cite{Me22}, \cite{MS201}, \cite{Me27})}
$G= G_{16} = \left\langle \frac{1}{\sqrt{2}} \left( {1 \atop 1} ~{1 \atop -1} \right) , ~ \left( {1 \atop 0} ~ {0 \atop -1} \right) \right\rangle \cong$
dihedral group $D_{16}$ (Shephard and Todd \#2b), order 16
$$\Phi = \frac{1}{(1- \la^2) (1- \la^8)}$$
\beql{EqG6a}
R: \frac{1}{\phi_2,~ \theta_8} ~,
\eeq
where $\phi_2 = x^2 + y^2$, $\theta_8 = x^2 y^2 (x^2 - y^2)^2$.
For example, since a Type II code is also a Type I code, the weight enumerator
of $g_{24}$, \eqn{Eq13}, must be in this ring.
It is:
$$
\phi_{24} = \phi_2^{12} - 12 \phi_2^8 \theta_8 + 6 \phi_2^4 \theta_8^2 - 64 \theta_8^2 ~.
$$
\paragraph{hwe of shadow $S$.}
It follows from Theorem \ref{GP4} that if $C$ has weight enumerator $W(x,y)$ then its shadow has weight enumerator $S(x,y) = W((x+y) / \sqrt{2} , i(x-y)/ \sqrt{2} )$.
This map from $W$ to $S$ preserves multiplication and addition, so to evaluate it it suffices to consider the images of the generators of the above ring.
We find that $x^2 + y^2$ becomes $2xy$ and $x^2 y^2 (x^2 - y^2)^2$ becomes $-4(x^4-y^4)^2$.
So $S(x,y)$ belongs to the ring
\beql{EqG6aa}
R: \frac{1}{xy,~ (x^4 - y^4)^2} ~.
\eeq

In particular, every element of the shadow has weight congruent to $n/2$ $\bmod~4$ (since this is true of the generators).

The shadow must satisfy an additional constraint.
If $C$ is Type I, let $W^{(j)} (x,y)$ be the weight enumerator of coset $C_j$, $j=0, \ldots, 3$ (see \eqn{EqR14a}).
Then $W^{(1)} (x,y) - W^{(3)} (x,y)$ is (up to sign) a multiplicative function
on codes, i.e., if $C$ is the direct sum of two Type I codes $C'$ and $C''$, then $W^{(1)} - W^{(3)}$ for $C$ is $\pm 1$ times the product of the polynomials $W^{(1)} - W^{(3)}$ for $C'$ and $C''$.
In order for this property to still hold when one (or both) of $C'$ and $C''$ is of Type~II, we
adopt the convention that for a Type~II code, $W^{(1)} - W^{(3)}$
is simply the weight enumerator of the code.

Then the additional condition satisfied by the shadow is that (if $C$ is Type I {\em or} Type II)
$W^{(1)} (x,y) - W^{(3)} (x,y)$ is a relative invariant for the group $G_{192}$ (see
\eqn{EqG5a}) with respect to the character
$$
\chi \left( \frac{1}{\sqrt{2}} \left( {1 \atop 1} ~ {1 \atop -1} \right) \right) = i^n ,~~
\chi \left( \left( {1 \atop 0} ~ {0 \atop i} \right) \right) = \eta^n
$$
where $\eta = (1+i) / \sqrt{2}$ \cite{Me158}.
An equivalent assertion is that $W^{(1)} - W^{(3)}$
is an absolute invariant for the subgroup of $G_{192}$ with determinant 1.

It follows (see \cite{Me158} for the proof) that for a Type I code $W^{(1)} (x,y) - W^{(3)} (x,y)$ lies in the following ring:
$$\Phi = \frac{1+ \la^{18}}{(1- \la^8) (1- \la^{12})}$$
\beql{EqG8a}
R: \frac{1,~ xy(x^8 - y^8) (x^8 -34 x^4 y^4 + y^8)}{x^8 + 14 x^4 y^4 + y^8, ~x^2 y^2 (x^4 - y^4)^2} ~.
\eeq

One of the differences between binary codes of Types I and II is that whereas the weight
enumerator of the former is invariant under a group of order only 16, the weight
enumerator of the latter is invariant under a group of order 192 (see Eq.~\eqn{EqG5a}).
The above result restores the balance to a certain extent, by requiring $W^{(1)} - W^{(3)}$ to be a
relative invariant for the larger group.

\subsection{Family \EoneII: Doubly-even binary self-dual codes}  \label{GL1b}
\paragraph{hwe of $C$} (\cite{Gle70}, \cite{Me22}, \cite{MS201}, \cite{Me27})
\beql{EqG5a}
G = G_{192} = \left\langle \frac{1}{\sqrt{2}} \left( {1 \atop 1} ~ {1 \atop -1} \right) , ~\left( {1 \atop 0} ~ {0 \atop i} \right) \right\rangle ~, \quad {\rm order~192}
\eeq
(Shephard and Todd \#9)
$$\Phi = \frac{1}{(1- \la^8) (1- \la^{24})}$$
\beql{EqG5b}
R: \frac{1}{x^8 + 14 x^4 y^4 + y^8 ,~ x^4 y^4 (x^4 - y^4)^4}~.
\eeq
Codes whose weight enumerators give generators for the above rings:
\beql{Eq5c}
\begin{array}{cl} \hline
\mbox{Ring} & \mbox{Codes} \\ \hline
\mbox{\eqn{EqG6a}} & \mbox{$i_2$ \eqn{Eq11}, $e_8$ \eqn{Eq12}} \\
\mbox{\eqn{EqG6aa}} & \mbox{$i_2$ \eqn{Eq11}, $e_8$ \eqn{Eq12}} \\
\mbox{\eqn{EqG8a}} & \mbox{$e_8$ \eqn{Eq12}, $g_{24}$ \eqn{Eq13}; $d_{12}^+$, $(d_{10} e_7 f_1)^+$ (\S\ref{GU3})} \\
\mbox{\eqn{EqG5b}} & \mbox{$e_8$ \eqn{Eq12}, $g_{24}$ \eqn{Eq13}}
\end{array}
\eeq
\paragraph{Remark.}
The above groups $G_{16}$ and $G_{192}$ are also the two-dimensional real and complex Clifford groups\index{Clifford group} occurring in quantum coding theory
\cite{Me222}, \cite{Me223}.
\index{group!Clifford}\index{code!quantum}\index{quantum code}
At present this appears to be nothing more than a coincidence.
However, in view of the other mysterious coincidences involving the Clifford groups,
there may be a deeper explanation that is presently hidden (compare the remarks in Section~\ref{GL7}).
\subsection{Family \Etwo: Ternary codes} \label{GL2}
\paragraph{hwe of $C$} {\rm (\cite{Gle70}; \cite{Me22}, \cite{Me27}, \cite[p.~620]{MS})}
\hsp
$$G= \left\langle \frac{1}{\sqrt{3}} \left[ {1 \atop 1} ~ {2 \atop -1} \right] ,~
\left[ {1 \atop 0} ~ {0 \atop \om} \right] , ~
\om = e^{2 \pi i/3} \right\rangle, \quad {\rm order~48}
$$
(Shephard \& Todd \#6)
\beql{EqG11a}
\Phi = \frac{1}{(1- \la^4 )(1- \la^{12})}
\eeq
\beql{EqG11b}
R: \frac{1}{x^4 + 8xy^3 , ~ y^3 (x^3 - y^3)^3}
\eeq
\paragraph{cwe of $C$, $1^n \in C$}
(\cite{Me80}; \cite[p.~617]{MS})
(This forces the length to be a multiple of 12.)
$$G =\left\langle \frac{1}{\sqrt{3}}
\left[ \begin{array}{ccc}
1 & 1 & 1 \\
1 & \om & \oom \\
1 & \oom & \om
\end{array}
\right] , ~~
\left[
\begin{array}{ccc}
1 & 0 & 0 \\
0 & 0 & 1 \\
0 & 1 & 0
\end{array}
\right] , ~~
\left[
\begin{array}{ccc}
0 & 1 & 0 \\
0 & 0 & 1 \\
1 & 0 & 0
\end{array}
\right], ~~
\left[
\begin{array}{ccc}
1 \\
~ & \om \\
~ & ~ & \om
\end{array}
\right] \right\rangle ~,
$$
order 2592,
$$\Phi = \frac{1+ \la^{24}}{(1- \la^{12})^2 (1- \la^{36})}$$
\beql{EqG13a}
R: \frac{1,~ \beta_6 \pi_9^2}{\af_{12},~ \beta_6^2 ,~ \pi_9^4}
\eeq
where
\begin{eqnarray*}
&&a = x^3 + y^3 + z^3 ,~
b= x^3 y^3 + y^3 z^3 + z^3 x^3 , \\
&& p= 3xyz, ~ \af_{12} = a(a^3 + 8p^3) , \\
&& \beta_6 = a^{12} - 12b , ~ \pi_9 = (x^3 - y^3) (y^3 - z^3) (z^3 - x^3) ~.
\end{eqnarray*}
\paragraph{cwe of $C$, not requiring that $1^n \in C$}
\cite{McEliece} (Now the length is just a multiple of 4.)
$$G =\left\langle \frac{1}{\sqrt{3}}
\left[ \begin{array}{ccc}
1 & 1 & 1 \\
1 & \om & \oom \\
1 & \oom & \om
\end{array}
\right] , ~~
\left[
\begin{array}{ccc}
1 & 0 & 0 \\
0 & 0 & 1 \\
0 & 1 & 0
\end{array}
\right] , ~~
\left[
\begin{array}{ccc}
1 \\
~ & \om \\
~ & ~ & \om
\end{array}
\right] \right\rangle ~,
$$
order 96
$$\Phi = \frac{1+ 4 \la^{12} + \la^{24}}{(1- \la^4) (1- \la^{12})^2} ~.$$
$$
R: \sum_{i=0}^6 f^{(i)} S
$$
where $S= \CC [\theta_4, \theta_6^2 , t^{12} ]$,
$s=y+z$, $t=y-z$,
$\theta_4 = x(x^3 + s^3)$, $\theta_6 = 8x^6 - 20x^3 s^3 -s^6$, $f^{(0)} =1$,
$f^{(1)} = t^2 \phi_4 \theta_6$, $\phi_4 = s(8x^3 -s^3)$, $f^{(2)} = t^4 \phi_4^2$,
$f^{(3)} = t^6 \theta_6$, $f^{(4)} = t^8 \phi_4$, $f^{(5)} = t^{10} \phi_4^2 \theta_6$.

\noindent{\bf Codes:}
\beql{EqG14aa}
\begin{array}{ll} \hline
\mbox{Ring} & \mbox{Codes} \\ \hline
\mbox{\eqn{EqG11b}} & \mbox{$t_4$ \eqn{Eq14}, $g_{12}$ \eqn{Eq15}} \\
\mbox{\eqn{EqG13a}} & \mbox{$e_3^{4+}$ (\S\ref{GU4}), $g_{12}$ \eqn{Eq15}, $S(36)$ (\S\ref{EX2}); $XQ_{23}$ (\S\ref{EX2})}
\end{array}
\eeq
\subsection{Family \Ethree: Self-dual codes over $\FF_4$ with Hermitian inner product} \label{GL3}
\paragraph{hwe of $C$}
(\cite{Me27}; \cite[p.~621]{MS}):
$$G = \left\langle \frac{1}{2} \left( \begin{array}{cc}
1 & 3 \\ 1 & -1
\end{array}
\right) , ~
\left( \begin{array}{cc}
1 & 0 \\ 0 & -1
\end{array}
\right) \right\rangle
= \mbox{Weyl group of type $G_2 \cong$ dihedral group $D_{12}$}
$$
(Shephard \& Todd \#2b)
$$\Phi = \frac{1}{(1- \la^2)(1- \la^6)}$$
\beql{EqG15a}
R: \frac{1}{x^2 + 3y^2,~ y^2 (x^2 - y^2)^2}
\eeq
\paragraph{cwe of $C$, $1^n \in C$}
(There must be some word of full weight, so this is not a severe restriction)
$$G = \left\langle
\frac{1}{2} \left[
\begin{array}{rrrr}
1 & 1 & 1 & 1 \\
1 & 1 & -1 & -1 \\
1 & -1 & 1 & -1 \\
1 & -1 & -1 & 1
\end{array} \right] ,~
\left[
\begin{array}{rrrr}
1 \\
~ & -1 \\
~ & ~ & -1 \\
~ & ~ & ~ & -1
\end{array} \right] ,~
\left[
\begin{array}{cccc}
0 & 1 & 0 & 0 \\
1 & 0 & 0 & 0 \\
0 & 0 & 0 & 1 \\
0 & 0 & 1 & 0
\end{array} \right] ,~
\left[
\begin{array}{cccc}
1 & 0 & 0 & 0 \\
0 & 0 & 1 & 0 \\
0 & 0 & 0 & 1 \\
0 & 1 & 0 & 0
\end{array} \right] \right\rangle
$$
order 576
$$
\Phi = \frac{1+ \la^{12}}{(1- \la^2 ) (1- \la^6) (1- \la^8) (1- \la^{12} )}
$$
\beql{EqG16a}
R: \frac{1,~(x^2 - y^2) (x^2 -z^2 ) (x^2 - t^2 ) (y^2 - z^2 ) (y^2 - t^2) (z^2 - t^2)}{x^2 + y^2 + z^2 + t^2 , ~\mbox{\eqn{Eq19}} , ~ f_8 , ~ f_{12}}
\eeq
where
\begin{eqnarray*}
f_8 &= &
x^8 + \cdots (\mbox{4 terms}) + 14x^4 y^4 + \cdots (\mbox{6 terms})
+ 168 x^2 y^2 z^2 t^2 = \mbox{cwe of $e_8 \otimes \FF_4$} \\
f_{12} & = & (s_4 - 3x^2 y^2 -3z^2 t^2 )(s_4 -3x^2 z^2 - 3y^2 t^2) (s_4 - 3x^2 t^2 - 3y^2 z^2) ~, \\
s_4 & = & x^2 y^2 + x^2 z^2 + \cdots ~ ( \mbox{6 terms}) ~.
\end{eqnarray*}

\paragraph{cwe of $C$, assuming $1^n \in C$ and $C$ and $\overline{C}$ have same cwe:} 
\begin{eqnarray*}
G & = & \mbox{previous $G$ together with} ~
\left[ \begin{array}{cccc}
1 & 0 & 0 & 0 \\
0 & 1 & 0 & 0 \\
0 & 0 & 0 & 1 \\
0 & 0 & 1 & 0
\end{array} \right]  \\
& \cong & \mbox{Weyl group of type $F_4$ (Shephard \& Todd \#28), order 1152}
\end{eqnarray*}
$$\Phi = \frac{1}{(1- \la^2) (1- \la^6) (1- \la^8) (1- \la^{12})}$$
\beql{EqG18a}
R: \frac{1}{x^2 + y^2 +z^2 + t^2, ~ \mbox{\eqn{Eq19}} , ~f_8 ,~ f_{12}}
\eeq
\paragraph{swe of $C$, $1^n \in C$:}
(Set $t=z$ in cwe)
$$\Phi = \frac{1+ \la^{12}}{(1- \la^2) (1- \la^6) (1- \la^8)}$$
\beql{EqG19a}
R: \frac{1, \{ (x^2 - z^2) (y^2 - z^2) \}^3}{x^2 + y^2 + 2z^2 , ~ \mbox{\eqn{Eq18}},~ \{ (x^2 - z^2) (y^2 -z^2) \}^2}
\eeq
\paragraph{Remark.}
If we try to apply invariant theory directly to the swe, we are led to the group
$$G = \left\langle
\frac{1}{2}
\left( \begin{array}{rrr}
1 & 1 & 2 \\
1 & 1 & -2 \\
1 & -1 & 0
\end{array} \right) , ~
\left( \begin{array}{ccc}
0 & 1 & 0 \\
1 & 0 & 0 \\
0 & 0 & 1
\end{array} \right) , ~
\left( \begin{array}{ccc}
1 \\
~ & -1 \\
~ & ~ & -1
\end{array} \right) \right\rangle
$$
(Weyl group of type $B_3$, Shephard \& Todd \#2a) of order
48, with Molien series
$$\Phi = \frac{1}{(1- \la^2) (1- \la^4) (1- \la^6)} ~.$$
However, the invariant of degree 4 is
$$\delta_4 = (x^2 - z^2 ) (y^2 - z^2) ~,$$
which cannot be obtained from the swe of any self-dual
code of length 4.
The ring of invariants here and the ring in \eqn{EqG19a} have the same
quotient field.
So there is no group whose ring of invariants is \eqn{EqG19a}.

\noindent{\bf Codes:}
\beql{EqG21a}
\begin{array}{ll} \hline
\mbox{Ring} & \mbox{Codes} \\ \hline
\mbox{\eqn{EqG15a}} & \mbox{$i_2$ \eqn{Eq16a}, $h_6$ \eqn{Eq17}} \\
\mbox{\eqn{EqG16a}} & \mbox{$i_2$ \eqn{Eq16a}, $h_6$ \eqn{Eq17}, $e_8 \otimes \FF_4$, $(e_7 e_5)^+$ (\S\ref{GU5}); $d'_{12}$} \\
\mbox{\eqn{EqG18a}} & \mbox{$i_2$ \eqn{Eq16a}, $h_6$ \eqn{Eq17}, $e_8 \otimes \FF_4$, $(e_7 e_5)^+$ (\S\ref{GU5})} \\
\mbox{\eqn{EqG19a}} & \mbox{$i_2$ \eqn{Eq16a}, $h_6$ \eqn{Eq17}, $e_8 \otimes \FF_4$; $(e_7 e_5)^+$ (\S\ref{GU5})}
\end{array}
\eeq
Here $d'_{12}$ is the code obtained from $d_{12}^+$ of Section~\ref{GU3} by multiplying the last four coordinates by $\om$.
\subsection{Family \Efour: Self-dual codes over $\FF_4$ with Euclidean inner product}\label{GL4}
\hsp
(This is inadequately treated in \cite{Me55}, where only even codes are considered.)
Neither the hwe nor the swe can be obtained directly from invariant theory, but must be obtained by collapsing the cwe.
Since $(v,v) =0 \Leftrightarrow \sum v_i^2 =0 \Leftrightarrow \sum v_i =0 \Leftrightarrow (v, 1^n) =0$, we may assume $1^n \in C$.
\paragraph{cwe of $C$}
$$G = \left\langle
\frac{1}{2}
\left( \begin{array}{rrrr}
1 & 1 & 1 & 1 \\
1 & 1 & -1 & -1 \\
1 & -1 & -1 & 1 \\
1 & -1 & 1 & -1
\end{array} \right) , ~
\left( \begin{array}{cccc}
0 & 1 & 0 & 0 \\
1 & 0 & 0 & 0 \\
0 & 0 & 0 & 1 \\
0 & 0 & 1 & 0
\end{array} \right) , ~
\left( \begin{array}{cccc}
1 & 0 & 0 & 0 \\
0 & 0 & 1 & 0 \\
0 & 0 & 0 & 1 \\
0 & 1 & 0 & 0
\end{array} \right)\right\rangle
$$
order 192
$$\Phi = \frac{1+ \la^{16}}{(1- \la^2 ) (1- \la^4) (1- \la^6 ) (1- \la^8 )}$$
\beql{EqG24a}
R: \frac{1,~ abcd (a^2 - b^2) (a^2 -c^2) \cdots (c^2 - d^2)}{\mbox{symmetric polynomials in $a^2, b^2, c^2, d^2$}}
\eeq
where
\begin{eqnarray*}
a = (+x-y-z-t) /2, && b= (-x+y-z-t)/2 ~, \\
c = (-x-y+z-t)/2 , && d= (-x-y-z+t)/2 ~.
\end{eqnarray*}
\paragraph{cwe of $C$, assuming $\overline{C}$ has same cwe as $C$:}
$$G = \mbox{previous group together with}
\left( \begin{array}{cccc}
1 & 0 & 0 & 0 \\
0 & 1 & 0 & 0 \\
0 & 0 & 0 & 1 \\
0 & 0 & 1 & 0
\end{array}
\right)
$$
(Weyl group of type $B_4$, Shephard \& Todd \#2a), order 384
$$
\Phi = \frac{1}{(1- \la^2) (1- \la^4) (1- \la^6) (1- \la^8)}
$$
\beql{EqG26a}
R: \mbox{symmetric polynomials in}~ a^2, b^2, c^2, d^2 ~.
\eeq
\paragraph{swe of $C$:}
(Set $t=z$ in the above cwe)
$$\Phi = \frac{1+ \la^8 + \la^{16}}{(1- \la^2 )(1- \la^4) (1- \la^6)}$$
\beql{EqG27a}
R: \frac{1,~ \{ (x^2 -z^2) (y^2 - z^2) \}^2,~ \{ (x^2 - z^2) (y^2 - z^2) \}^4}{x^2 + y^2 + 2z^2 ,~ x^4 + y^4 + 2z^4 + 12 xyz^2,~ z^2 (x-y)^2 (xy-z^2)} ~.
\eeq
\paragraph{hwe of $C$:}
(Set $t=z=y$ in the cwe)
$$\Phi = \frac{1+ \la^6}{(1- \la^2) (1- \la^4)}$$
\beql{EqG27b}
R: \frac{1,~ y^2 (x^2 - y^2)^2}{x^2 + 3y^2 ,~ y^2(x-y)^2}~.
\eeq
Rather surprisingly,
\eqn{EqG24a}, \eqn{EqG27a}, \eqn{EqG27b} appear to be new.

\noindent{\bf Codes:}
The following codes will be used:
\beql{EqG28aa}
\begin{array}{ll}
\multicolumn{2}{l}{i_2 = [11], ~~ cwe = x^2 + y^2 + z^2 + t^2 , ~ swe = x^2 + y^2 + 2z^2 , ~~ hwe = x^2 + 3y^2} \\ [+.1in]
\multicolumn{2}{l}{c_4 = \left[ \begin{array}{cccc}
1 & 1 & 1 & 1 \\
0 & 1 & \om & \oom
\end{array}
\right] , ~~ \mbox{a [4,2,3] Reed-Solomon code,}} \\ [+.1in]
\multicolumn{2}{l}{~~~~~~~~cwe = x^4 + y^4 + z^4 + t^4 + 12 xyzt, ~swe = x^4 + y^4 + 2z^4 + 12 xyz^2 ,} \\
\multicolumn{2}{l}{~~~~~~~~hwe = x^4 + 12 xy^3 + 3y^4 ~.} \\ [+.1in]
c_6 = \left[ \begin{array}{cccccc}
1 & 1 & 1 & 1 & 1 & 1 \\
0 & 0 & 0 & 1 & \om & \oom \\
1 & \om & \oom & 0 & 0 & 0
\end{array}
\right], \\ [+.1in]
\multicolumn{2}{l}{~~~~~~~~cwe = x^6 + \cdots \mbox{(4 terms)} + 6x^3 yzt + \cdots \mbox{(4 terms)} + 9x^2y^2z^2 + \cdots \mbox{(4 terms)} ,} \\
\multicolumn{2}{l}{~~~~~~~~hwe = x^6 + 6x^3 y^3 + 27 x^2 y^4 + 18xy^5 + 12y^6 ~.}
\end{array}
\eeq

\noindent{\bf Codes:}
\beql{EqG28a}
\begin{array}{ll} \hline
\mbox{Ring} & \mbox{Codes} \\ \hline
\mbox{\eqn{EqG24a}} & i_2, c_4 , c_6, e_8 \otimes \FF_4 ; ? \\
\mbox{\eqn{EqG26a}} & i_2, c_4, c_6, e_8 \otimes \FF_4 \\
\mbox{\eqn{EqG27a}} & i_2, c_4, c_6 ; e_8 \otimes \FF_4 \\
\mbox{\eqn{EqG27b}} & i_2, c_4; c_6
\end{array}
\eeq
\paragraph{Remark.}
The question mark in the first line of the table indicates that we do not have a code that produces the degree 16 polynomial in the numerator of \eqn{EqG24a}.
Such a code would necessarily be odd and
have the property that the cwe of $\overline{C}$ is not equal to that of $C$.
Presumably a random self-dual code would do, but we would prefer to find a code with some nice structure.
\subsection{Family \EfiveI: Additive self-dual codes over $\FF_4$ using trace inner product}\label{GL5a}
\paragraph{hwe of $C$:}
$$G = \left\langle \frac{1}{2}
\left( \begin{array}{rr}
1 & 3 \\ 1 & -1
\end{array} \right) \right\rangle ~, \quad {\rm order~2}
$$
$$
\Phi = \frac{1}{(1- \la) (1- \la^2)}
$$
\beql{EqG29a}
R: \frac{1}{x+y ,~ y(x-y)}
\eeq
\paragraph{cwe of $C$, $1^n \in C$:}
$$G = \left\langle M_4 = \frac{1}{2} \left( \begin{array}{rrrr}
1 & 1 & 1 & 1 \\
1 & 1 & -1 & -1 \\
1 & -1 & 1 & -1 \\
1 & -1 & -1 & 1
\end{array}
\right) ~, \quad
\af_4 = \left( \begin{array}{cccc}
0 & 1 & 0 & 0 \\
1 & 0 & 0 & 0 \\
0 & 0 & 0 & 1 \\
0 & 0 & 1 & 0
\end{array} \right) \right\rangle ,
{\rm order}~ 8
$$
$$\Phi = \frac{1+ \la^3}{(1- \la) (1- \la^2)^2 (1- \la^4)}$$
\beql{EqG30b}
R: \frac{1,~ BCD}{A, B^2 + C^2,~ D^2,~ B^2 C^2}
\eeq
where
\beql{EqG30a}
A = (x+y)/2, B= (x-y)/2, C= (z+t)/2, D= (z-t)/2 ~.
\eeq
\paragraph{swe of $C$, $1^n \in C$:}
(Set $t=z$ in cwe)
$$\Phi = \frac{1}{(1- \la) (1- \la^2) (1- \la^4)}$$
\beql{EqG31b}
R: \frac{1}{A,~ B^2 + C^2,~ B^2C^2}
\eeq
where
\beql{EqG31a}
A = (x+y)/2, B= (x-y) /2, C=z ~.
\eeq
\paragraph{cwe of $C$, $1^n \in S$:}
Note that $(1^n ,u) = wt(u) - n_1 (u) \equiv wt(u)$ $(\bmod~2)$
if and only if the number of 1's in $u$ is even.
So if $1^n \in S$, the cwe is invariant under ${\rm diag} \{1,-1,1,1\}$.
$$
G= \langle M_4, ~{\rm diag} \{1,-1,1,1\} \rangle ~, \quad {\rm order~6}
$$
$$
\Phi = \frac{1}{(1- \la)^2 (1- \la^2) (1- \la^3)}
$$
\beql{EqG32a}
R: \frac{1}{D,~ A+B+C,~ A^2 + B^2 + C^2,~ A^3 + B^3 +C^3}
\eeq
where $A$, $B, \ldots$ are as in \eqn{EqG30a}.
\paragraph{swe of $C$, $1^n \in S$:}
(Set $t=z$ in cwe)
$$\Phi = \frac{1}{(1- \la) (1- \la^2) (1- \la^3)}$$
\beql{EqG33b}
R: \frac{1}{\mbox{symmetric polynomials in $A$, $B$, $C$}}
\eeq
\paragraph{hwe of $S$:}
$$\Phi : \frac{1}{(1- \la) (1- \la^2)}$$
\beql{EqG331a}
R: \frac{1}{2y,~ - \frac{1}{2} (x^2 - y^2)}
\eeq
As a corollary, the weight of a vector in the shadow is congruent
to $n$ $(\bmod~2)$.
\paragraph{hwe of $W^{(1)} - W^{(3)}$:}
Again we use the terminology $W^{(i)}$, $i=0, \ldots, 3$, for the cosets of
$C_0$ in $C_0^\perp$ (as in Sect.~\ref{Shad})
$$G = \left\langle M_2 = \frac{1}{2}
\left( \begin{array}{rr}
1 & 3 \\ 1 & -1
\end{array} \right) , \quad \af_2 =
\left( \begin{array}{rr}
1 & 0 \\ 0 & -1
\end{array} \right) \right\rangle
$$
with character $\chi (M_2) =1$, $\chi (\af_2) = (-1)^n$
(Ker $\chi \cong S_3$)
$$\Phi = \frac{1}{(1- \la^2) (1- \la^3)}$$
\beql{EqG332a}
R: \frac{1}{x^2 + 3y^2 ,~ y (x^2 - y^2)} ~.
\eeq
\paragraph{cwe of $S$, $1^n \in C$:}
Belongs to image of \eqn{EqG30b} under the map that sends
$(x,y,z,t)$ to $(x,y,z,t) \beta_4 M_4$:
\beql{EqG333a}
R: \frac{1,~ ABD}{C,~ A^2 + B^2,~ D^2 ,~ A^2 B^2}
\eeq
\paragraph{cwe of $W^{(1)} - W^{(3)}$, $1^n \in C$:}
$G = \langle M_4, \af_4, \beta_4 \rangle$ with character
$\chi (M_4) =1$, $\chi ( \af_4) = \chi (\beta_4) = (-1)^n$, order 48
$$\Phi = \frac{1}{(1- \la) (1- \la^2) (1- \la^3) (1- \la^4)}$$
\beql{EqG334a}
R: \frac{1}{D,~ A^2 + B^2+C^2,~ ABC,~ A^4 + B^4 + C^4}
\eeq
where
$$
A = x+y , ~ B = x-y , ~C = z+t ,~  D = z-t .
$$
\paragraph{swe of $W^{(1)} - W^{(3)}$, $1^n \in C$:}
$$
\Phi = \frac{1}{(1- \la^2)(1- \la^3) (1- \la^4)}
$$
\beql{EqG334b}
R: \frac{1}{A^2 + B^2 +C^2,~ ABC,~ A^4 + B^4 +C^4}
\eeq
\paragraph{cwe of $S$, $1^n \in S$:}
Belongs to image of \eqn{EqG32a} under the map $x \to y$, $y \to x$,
$z \to t$, $t \to z$:
\beql{EqG335a}
R: \frac{1}{-D,~ A-B+C,~ A^2 + B^2 +C^2,~ A^3 - B^3 + C^3}
\eeq
\paragraph{swe of $S$, $1^n \in S$:}
Set $D=0$ in \eqn{EqG335a}.
\paragraph{hwe of $S$, $1^n \in S$:}
Same as \eqn{EqG331a}.
\paragraph{cwe of $W^{(1)} - W^{(3)}$, $1^n \in S$:}
$G = \langle M_4, \beta_4 = {\rm diag} \{1, -1, -1, -1\} \rangle$,
with character $\chi (M_4) =1$, $\chi (\beta_4) = (-1)^n$, order 12
$$\Phi = \frac{1}{(1- \la)^2 (1- \la^2) (1- \la^3)}$$
\beql{EqG336a}
R: \frac{1}{D,~ A-B-C,~ A^2 + B^2 +C^2 ,~ A^3 - B^3 - C^3}
\eeq
\paragraph{Remark.}
We may obtain $W^{(1)} - W^{(3)}$ by applying $\af_4$ to $W^{(0)} - W^{(2)}$, which in turn is obtained by applying $\beta_4$ to $W^{(0)} +W^{(2)}$.
\subsection{Family \EfiveII: Additive even self-dual codes over $\FF_4$ using trace inner product}\label{GL5b}
\paragraph{hwe of $C$:}
Same as family \Ethree, see \eqn{EqG15a}.
\paragraph{cwe of $C$, $1^n \in C$:}
$$
G = \langle M_4, \af_4, \beta_4 \rangle~, \quad
\mbox{order 48}
$$
$$\Phi = \frac{1+ \la^4}{(1- \la^2)^2 (1- \la^4) (1- \la^6)}$$
\beql{EqG35a}
R: \frac{1, ~ABCD}{D^2,~ A^2 + B^2 + C^2,~ A^4 + B^4 + C^4 ,~ A^6 + B^6 + D^6}
\eeq
\paragraph{swe of $C$, $1^n \in C$:}
(Set $t=z$ in cwe)
$$\Phi = \frac{1}{(1- \la^2)(1- \la^4) (1- \la^6)}$$
\beql{EqG36a}
R: \frac{1}{\mbox{symmetric polynomials in $A^2$, $B^2$, $C^2$}}
\eeq
\paragraph{cwe of $C$ (not assuming $1^n \in C$):}
$$G = \langle M_4 , \beta_4 \rangle, \quad \mbox{order 12}$$
\beql{EqG37a}
\Phi = \frac{1+ \la^2 + 2 \la^4}{(1- \la^2)^3 (1- \la^6)}
\eeq

\noindent{\bf Codes:}
The following codes will be used:
\begin{eqnarray*}
i_1 & = & [1], \quad cwe = swe = hwe = x+y \\
i'_1 & = & [ \om ], \quad cwe = swe = x+z \\
i''_1 & = & [ \oom] ,\quad cwe = x+t, ~swe = x+z \\
i_2 & = & [11, \om \om], \quad \mbox{see \eqn{EqG28aa}} \\
i'_2 & = & [11, \om \oom ], \quad cwe = x^2 + y^2 + 2zt , swe = x^2 + y^2 + 2z^2, hwe = x^3 + 3y^2 \\
c_3 & = & [111, \om \om 0 , \om 0 \om ] \\
c_4 & = & [1111, \om (\om 00)] \\
c'_4 & = & [\om \om \om \om, 1(100)]
\end{eqnarray*}
\beql{EqG37ca}
\begin{array}{ll} \hline
\mbox{Ring} & \mbox{Codes} \\ \hline
\mbox{\eqn{EqG29a}} & i_1, i_2 \\
\mbox{\eqn{EqG30b}} & i_1, i_2, i'_2 , c_4; c_3 \\
\mbox{\eqn{EqG31b}} & i_1, i_2, c_4 \\
\mbox{\eqn{EqG32a}} & i'_1 , i''_1 , i_2, c'_3 \\
\mbox{\eqn{EqG33b}} & i'_1, i_2, c'_3 \\
\mbox{\eqn{EqG331a}} & i_1, i_2 \\
\mbox{\eqn{EqG333a}} & i_1, i_2, i'_2, c_4; c_3 \\
\mbox{\eqn{EqG336a}} & i'_1, i''_1 , i_2, c_3 \\
\mbox{\eqn{EqG35a}} & i_2 , i'_2, c'_4 , h_6; c_4 \\
\mbox{\eqn{EqG36a}} & i_2, c'_4, h_6
\end{array}
\eeq
\subsection{Family \Esix: Codes over $\FF_q$, $q$ a square, with Hermitian inner product}\label{GL6}
\hsp
The case $q=4$ has been studied in Section~\ref{GL3}.
The next case is $q=9$, but as little attention has been paid
so far to codes over this field we shall not discuss the cwe or swe further.
It is possible to say a little about the Hamming
weight enumerator in the general case.
\paragraph{hwe of $C$} (See Theorem \ref{ThMS6}):
$$G = \left\langle \frac{1}{\sqrt{q}} \left( \begin{array}{cc}
1 & q-1 \\ 1 & -1
\end{array}
\right) ,~ \left( \begin{array}{rr}
-1 & 0 \\ 0 & -1
\end{array} \right) \right\rangle,
~{\rm order~4}
$$
$$\Phi = \frac{1}{(1- \la^2)^2}$$
\beql{EqG37cb}
R: \frac{1}{x^2 + (q-1) y^2,~ y(x-y)}
\eeq
(This is somewhat unsatisfactory, since $y(x-y)$ forces a vector of weight 1,
which is impossible in a self-dual code.)
\subsection{Family \Eseven: Codes over $\FF_q$ with Euclidean inner product}\label{GL7}
\hsp
The cases $q= 2,3$ and 4 have been studied in Sections~\ref{GL1a}, \ref{GL2}, \ref{GL4}.
As $q$ increases the results rapidly become more complicated.

We first discuss the case $q=5$ and then say a little about the
general case.

\paragraph{cwe of $C$, $q=5$:}
Let $\xi = e^{2 \pi i/5}$.
$$G = \left\langle \frac{1}{\sqrt{5}}
( \xi^{rs} )_{r,s = 0, \ldots, 4} , ~
{\rm diag} \{1,\xi , \xi^{-1} , \xi^{-1}, \xi \} \right\rangle~, \quad
{\rm order~240}
$$
\beql{EqG41a}
\Phi = \frac{\phi ( \lambda)}{(1- \lambda^4) (1- \lambda^6) (1- \lambda^{10})^2}
\eeq
where $\phi ( \lambda )$ is a polynomial of degree 26, with $\phi (1) =60$.
A good basis for this ring would therefore involve about 65 polynomials!
Such Behavior is typical of most groups --- see Huffman and Sloane\index{Huffman-Sloane theorem}\index{theorem!Huffman-Sloane} \cite{Me61}.\paragraph{swe of $C$, $q=5$}
$$G = \left\langle \left(
\begin{array}{ccc}
1 & 2 & 2 \\
1 & \xi + \xi^4 & \xi^2 + \xi^3 \\
1 & \xi^2 + \xi^3 & \xi + \xi ^4
\end{array}\right),
{\rm diag} \{1, \xi , \xi^4 \} ,
\left( \begin{array}{ccc}
1 & 0 & 0 \\
0 & 0 & 1 \\
0 & 1 & 0
\end{array} \right) \right\rangle ~,
$$
(the reflection group $[3,5]$, a three-dimensional representation of the icosahedral group, Shephard and Todd \#23), order 120
$$\Phi = \frac{1}{(1- \la^2) (1- \la^6) (1- \la^{10})}$$
\beql{EqG42a}
R: \frac{1}{\af ,~ \beta ,~ \gamma}
\eeq
where
\begin{eqnarray*}
\af & = & x^2 + 4yz \\
\beta & = & x^4 yz - x^2 y^2 z^2 - x(y^5 + z^5 ) +2y^3 z^3 \\
\gamma & = & 5x^6 y^2 z^2 - 4x^5 (y^5 + z^5) -10x^4 y^3 z^3 +
10x^3 (y^6 z+ yz^6) + 5x^2y^4z^4 \\
&&~~~~~~~~~- ~ 10x (y^7 z^2 + y^2 z^7) + 6y^5 z^5 + y^{10} + z^{10} ~.
\end{eqnarray*}
{\bf Codes:}
$[12]$, $[(100)(133)]$ and either
$$
d_5^{2+} = [(01234)(00000), (00000)(01234), 1111111111]
$$
or
\beql{EqG43a}
e_{10}^+ = [(00014)(00023), 1111111111] 
\eeq
for the invariant of degree 10.

In \cite{Me89} it was observed that these invariants were already known to Klein\index{Klein, F.} \cite{Me8922}, \cite{Me8923}.
This paper then went on to remark that ``it is worth mentioning that precisely the same invariants have recently been studied by Hirzebruch\index{Hirzebruch, F.} in connection with
cusps of the Hilbert
modular surface associated with $\QQ ( \sqrt{5} )$ --- see \cite{Me8919}, p.~306.
However, there does not seem to be any connection between this work
and ours''.
An elegant explanation for this was soon found by Hirzebruch \cite{Hirz86}.
The basic idea is to take a self-dual code over $\FF_5$ and to obtain from it (using a version of Construction A\index{Construction A} \cite{SPLAG}) a lattice over $\ZZ[ \sqrt{5}]$.
The theta series of this lattice is a Hilbert modular form which can be written down from the swe of the code.
This produces an isomorphism between the ring of swe's and the appropriate ring of Hilbert modular forms.\index{Hilbert modular form}
\index{modular form!Hilbert}
The monograph \cite{Ebe94} gives a comprehensive account of these connections.

Incidentally, we do not know if the cwe ring described by \eqn{EqG41a} collapses to \eqn{EqG42a}.
\paragraph{hwe of $C$, $q=5$:}
(\cite{Patt80}) (Set $z=y$ in swe)
$$\Phi = \frac{1+ \la^{10} + \la^{20}}{(1- \la^2) (1- \la^6)}$$
$$R: \frac{1,~ \oga ,~ \oga^2}{\oaf,~ \obe}$$
where
\begin{eqnarray*}
\oaf & = & x^2 + 4y^2 , \\
\obe & = & y^2 (x-y)^2 (x^2 + 2xy + 2y^2) , \\
\oga & = & y^4 (x-y)^4 (5x^2 +12 xy + 8y^2) ~.
\end{eqnarray*}
\paragraph{cwe of $C$, $q=5$, $1^n \in C$:}
(The group is now considerably larger, but the ring of invariants is no simpler)
$$G = \langle \mbox{previous group}, ~
{\rm diag} \{ 1, \xi , \xi^2, \xi^3 , \xi^4 \} \rangle
$$
$\cong \pm 5^{1+2} . Sp_2 (5)$, a Clifford group\index{Clifford group} \cite{SS4}, \cite{SS5},
\index{group!Clifford}
\cite{CCKS96}, \cite{SS22} (see also \cite{HoMu73}), order 30000
$$
\Phi = \frac{1+ 3 \la^{20} + 13 \la^{30} + 18 \la^{40} + 28 \la^{50} + 34 \la^{60} +
17 \la^{70} + 4 \la^{80} + 2 \la^{90}}{(1- \la^{10}) (1- \la^{20} )^2 (1- \la^{30})^2} ~.
$$
The sum of the coefficients in the numerator is 120, so again there is no possibility of giving a good basis.

The degree 10 invariant is the cwe of either of the codes of length 10 given in \eqn{EqG43a}.
\paragraph{cwe of $C$, general $q$:}
It is hard to say anything in general, but if $q$ is an odd prime $p$ we can at least describe the structure of the group $G$ under which the cwe is invariant.
$$G = \left\langle M = \frac{1}{\sqrt{p}}
(\xi^{rs} )_{r,s=0, \ldots, p-1} ,~
J = {\rm diag} \{1, \xi , \xi^4 , \xi^9 , \ldots \}, -I \right\rangle ~.
$$
If $1^n \in C$ then the cwe is invariant under the larger group $G^+ = \langle G,P \rangle$, where
$$P: x_j \to x_{j+1} , ~~\mbox{(subscripts $\bmod~p$)}$$
We use $\Xi (H)$ to denote the center of a group $H$.
\begin{theorem}\label{thG50}
(a) Suppose $p \equiv 1$ $(\bmod~4)$.
Then $G$ has structure $Z(2) \times SL_2 ( p)$ and center $\Xi (G) = \langle -I \rangle$.
$G^+$ has structure $Z(2) \times p^{1+2}$
$SL_2 (p)$ and $\Xi (G^+) = \langle -I, \xi I \rangle$.
(b) Suppose $p \equiv 3~ ( \bmod~4)$.
Then $G$ has structure $Z(4) \times SL_2 (p)$ and $\Xi (G) = \langle iI \rangle$.
$G^+$ has structure $Z(4) \times p^{1+2} SL_2 (p)$ and
$\Xi (G^+) = \langle iI, \xi I \rangle$.
In either case $G$ and $G^+$ are preserved by the Galois group $Gal ( \QQ [\sqrt{p}, \xi ] / \QQ)$.
\end{theorem}
\paragraph{Remarks.}
(i) The group $G$ was first studied in the present context
by Gleason \cite{Gle70}.
The groups $G$ and $G^{+}$ (also for composite
 odd $q$, and with the appropriate modification for even $q$ as well)
are a special case of the construction in \cite{Weil64}.
Weil obtains analogs of $G^+$, in which $\FF_q$ can be replaced by any
locally compact abelian group isomorphic to its Pontrjagin dual.\footnote{We
are grateful to N. D. Elkies for this comment.}

(ii)~The analogous results for $p=2$ are given in
Sections~\ref{GL1a} and \ref{GL1b}.
(iii)~In both cases (a) and (b)
$G^+$ is the full normalizer (with coefficients restricted to $\QQ [ \sqrt{p}, \xi ]$) of the extraspecial $p$-group $E= \langle P, Q \rangle$, where
$Q = {\rm diag} \{ 1, \xi , \xi^2, \xi^3, \ldots \}$
(cf. \cite{CCKS96}).
\paragraph{Proof.}
$G$ normalizes $E$, since $MPM^{-1} = Q^{-1}$, $MQM^{-1} =P$, $JPJ^{-1} = \xi^a PQ^{-2}$,
$JQJ^{-1} = \xi^b Q$ for appropriate integers $a$ and $b$.
(Note that $\xi I = PQP^{-1} Q^{-1} \in E$.)
Thus we have a surjective homomorphism $\phi$ from $G$ to $SL_2 (p): M \to \left( {0 \atop 1} ~{-1 \atop 0} \right)$,
$J \to \left( {1 \atop -2} ~ {0 \atop 1} \right)$.
In particular $G$ is transitive on $E/ \langle \xi I \rangle$.

Suppose $G \cap E$ is nontrivial.
If there were a noncentral element of $E$ in $G$ then by the transitivity of $G$ it would follow that $\xi^c P \in G$ and $\xi^d Q \in G$ for some $c,d$.
But then $\xi I \in G$.
This would force the length of $C$ to be a multiple of $p$,
which is false (since there is always a code of length 4).
Hence $G \cap E = \{I\}$.

$E$ is irreducible, so the centralizer of $E$ consists only of multiples of $I$.
It follows that $\ker \phi$ consists of multiples of elements of $E$.
But the fourth power of an element of $\ker \phi$ would be in $E$, and this must be $I$.
Thus $\ker \phi$ is either $\langle -I \rangle$ or $\langle iI \rangle$.
If $p \equiv 1 $ $(\bmod~4)$ then $i \not\in \QQ [\sqrt{p} , \xi ]$, so the first possibility obtains.
It remains to show that $iI \in G$ when $p \equiv 3$ $(\bmod~4)$.
The matrix $(MJ)^pM^2$ is readily verified to belong to $\ker \phi$.
But $\det ((MJ)^p M^2) = ( \det M)^{p+2}$.
Since $M^2$ maps $x_j$ to $x_{-j}$, $\det M^2 = -1$, so $\det M = \pm i$.
It follows that $(MJ)^p M^2$ is $\pm iI$.~~~$\bsq$
\begin{coro}\label{co1}
If $p \equiv 3$ $(\bmod~4)$ then a self-dual code over $\FF_p$ must have length
divisible by 4.
\end{coro}
\paragraph{Proof.}
$iI \in G$.~~~$\bsq$

The conclusion
of Corollary~\ref{co1} also
holds for self-dual codes over $\FF_q$, $q \equiv 3$ $(\bmod~4)$ \cite{Ple7}.
\paragraph{hwe of $C$, general $q$:}
Belongs to the ring \eqn{EqG37cb}.
If $q \equiv 3$ $(\bmod~4)$ we can say more (\cite{HaOu97}):
$$
G= \left\langle \frac{1}{\sqrt{q}}
\left( \begin{array}{rr}
1 & q-1 \\ 1 & -1
\end{array}\right) , \quad
\left( \begin{array}{cc}
i & 0 \\ 0 & i
\end{array} \right) \right \rangle , \quad {\rm order ~8}
$$
$$\Phi = \frac{1+ \la^4}{(1- \la^4)^2}$$
$$R: \frac{1, ~ x^2y^2 - 2xy^3 + y^4}{x^4 +4(q-1)xy^3 + (q-1) (q-3) y^4,~x^3 y+ (q-3) xy^3 - (q-2)y^4}
$$

\subsection{Family \EeightI: Self-dual codes over $\ZZ_4$}\label{GL8}
\paragraph{cwe of $C$:} (\cite{Klemm89})
$$G = \left\langle M_4 = \frac{1}{2} \left(
\begin{array}{rrrr}
1 & 1 & 1 & 1 \\
1 & i & -1 & -i \\
1 & -1 & 1 & -1 \\
1 & -i & -1 & i
\end{array} \right) , ~
\af_4 = {\rm diag} \{ 1, i,1,i\} \right\rangle ,
~\quad{\rm order~64}
$$
$$\Phi = \frac{1+ \la^{10}}{(1- \la) (1- \la^4)^2 (1- \la^8)}$$
\beql{EqG60a}
R: \frac{1,~ (BCD)^2 (B^4 - C^4)}{A,~ B^4 + C^4,~ D^4,~ B^4 C^4}
\eeq
where
\beql{EqG60b}
A = x+z, ~B= y+t , ~C= x-z , ~ D= y-t ~.
\eeq
\paragraph{swe of $C$:}
(Set $t=y$ in cwe)
$$\Phi = \frac{1}{(1- \la) (1 - \la^4) (1- \la^8)}$$
\beql{EqG61a}
R: \frac{1}{A,~ B^4 +C^4 ,~ B^4 C^4}
\eeq
\paragraph{hwe of $C$:}
(Set $t=z=y$ in cwe)
$$\Phi = \frac{1+ \la^8}{(1- \la) (1- \la^4)}$$
\beql{EqG62a}
R: \frac{1,~ y^4 (x-y)^4}{x+y,~ y(x-y) (x^2 + xy + 2y^2)}
\eeq
\paragraph{cwe, $1^n \in C$}
$$G = \left\langle M_4, \af_4, \left(
\begin{array}{cccc}
0 & 1 & 0 & 0 \\
0 & 0 & 1 & 0 \\
0 & 0 & 0 & 1 \\
1 & 0 & 0 & 0
\end{array}
\right) \right\rangle~,\quad{\rm order~1024}
$$
$$
\Phi = \frac{(1+ \la^{12} ) (1+ \la^{16} )}{(1- \la^4) (1- \la^8)^2 (1- \la^{16})}
$$
\beql{EqG63a}
R: \frac{(1,~ A^{12} + B^{12} + C^{12} + D^{12} ) \times (1,~\sigma_{16} )}{A^4 + B^4 + C^4 +D^4 , ~ A^8 + B^8 + C^8 + D^8 , ~ \sigma_8 , ~ A^4 B^4 C^4 D^4}
\eeq
where
\begin{eqnarray*}
\sigma_8 & = & A^4 D^4 + B^4 C^4 ~, \\
\sigma_{16} & = & (ABCD)^2 (A^4 B^4 + C^4 D^4 - A^4 C^4 - B^4 D^4 )
\end{eqnarray*}
\paragraph{swe of $C$, $\pm 1^n \in C$}
(Set $t=y$ in cwe)
$$\Phi = \frac{1+ \la^{12}}{(1-\la^4)(1- \la^8)^2}$$
\beql{EqG64a}
R: \frac{1,~ A^4 B^4 C^4}{A^4 + B^4 + C^4 , ~ A^8 ~ B^8 ~ C^8, ~ B^4 C^4}~.
\eeq
This ring may also be described as $R_0 \oplus B^4 C^4 R_0 \oplus B^8 C^8 R_0$,
where $R_0$ is the ring of symmetric polynomials in $A^4$, $B^4$, $C^4$.
\paragraph{hwe of $C$, $\pm 1^n \in C$:}
(Set $t=z=y$ in cwe)
$$\Phi = \frac{(1+ \la^8) (1+ \la^{12})}{(1- \la^4) (1- \la^8)}$$
\beql{EqG64b}
R: \frac{(1,~y^2 (x^2 + 3y^2) (x^2 - y^2)^2) \times (1,~ y^4 (x^2 - y^2)^4)}{(x^2 + 3y^2)^2,~ y^4(x-y)^4}
\eeq
\paragraph{cwe of $C$, $1^n \in S$:}
If $1^n \in S$, Part (i) of Theorem~\ref{GP6} implies that if a vector
$0^a 1^b 2^c d^d \in C$ then $b-d +2c \equiv \frac{1}{2} (b+d)+ 2c$ $(\bmod~4)$,
i.e. $b \equiv 3d$ $(\bmod~8)$, and so $\beta_4 = {\rm diag} \{1,\eta , 1, \eta^5\} \in G$.
$$G = \langle M_4, \beta_4 \rangle , ~~{\rm order} ~192$$
$$\Phi = \frac{1+ \la^{18}}{(1- \la) (1- \la^4) (1- \la^8) (1- \la^{12})}$$
\beql{EqG65a}
R:
\frac{1,~ B^2 C^2 D^2 (B^4 -C^4) (B^4 +D^4) (C^4 +D^4)}{A,~ B^4 + C^4 - D^4,~ B^8 + C^8 +D^8,~ B^{12} +C^{12} - D^{12}}
\eeq
\paragraph{swe and hwe of $C$, $\pm 1^n \in S$:}
Same as \eqn{EqG61a} and \eqn{EqG62a}, respectively.
\paragraph{cwe of $S$:}
the image of \eqn{EqG60a} under $A \to B$, $B \to \eta  C,C \to A$, $D \to \eta^3 D$
$$\Phi = \frac{(1+ \lambda^{10})}{(1- \lambda ) (1- \lambda^4)^2 (1- \lambda^8)}$$
\beql{EqG65b}
R: \frac{1,~A^2 C^2 D^2 ( -A^4 -C^4)}{B,~A^4 -C^4, ~-D^4, ~ -A^4 C^4}
\eeq
\paragraph{swe of $S$:}
the image of \eqn{EqG61a} under $A \to B$, $B \to \eta C$, $C \to A$
$$\Phi = \frac{1}{(1- \lambda) (1- \lambda^4) (1- \lambda ^8)}$$
\beql{EqG65c}
R: \frac{1}{B, ~ A^4 -C^4 , ~ -A^4 C^4}
\eeq
It follows that the norms of vectors in the shadow are congruent to $n \bmod~8$.
\paragraph{cwe of $S$, $1^n \in C$:}
the image of \eqn{EqG64a} under $A \to B$, $B \to \eta C$, $C \to A$, $D \to \eta^3 D$.
\paragraph{swe of $S$, $\pm 1^n \in C$:}
$$
\Phi = \frac{1+ \lambda^{12}}{(1- \lambda^4)(1- \lambda^8)^2}
$$
\beql{EqG65d}
R: \frac{1, ~ -A^4 B^4 C^4}{A^4 + B^4 - C^4 ,~ A^8+B^8 + C^8 , ~ - A^4 C^4}
\eeq
\paragraph{cwe of $S$, $1^n \in S$:}
\beql{EqG65e}
R: \frac{1,~A^2 C^2 D^2 (A^4 + C^4 ) (C^4 + D^4) (A^4 - D^4)}{B,~ A^4 -C^4 + D^4 , ~ A^8 + C^8 + D^8 , ~ A^{12} - C^{12} + D^{12}}
\eeq
\paragraph{swe of $S$, $\pm 1^n \in S$:}
same as \eqn{EqG65c}.
\paragraph{cwe of $W^{(1)} - W^{(3)}$:}
$$G = \langle M_4, \gamma_4 = {\rm diag} \{ 1, \eta , -1, \eta \} \rangle~,
\quad {\rm order~768} ~,$$
with character $\chi (M_4) =i^n$, $\chi ( \gamma_4) = \eta^n$
$$\Phi = \frac{1+ \la^{18}}{(1- \la) (1- \la^4) (1- \la^8 ) (1- \la^{12})}$$
\beql{EqG66a}
R: \frac{1, A^2 B^2 C^2 (A^4 + B^4 ) (A^4 + C^4 ) (B^4 - C^4)}{D,~ \mbox{symmetric polynomials in} ~A^4, - B^4, -C^4}
\eeq
\paragraph{swe of $W^{(1)} - W^{(3)}$:}
$$\Phi = \frac{1+ \la^{18}}{(1- \la^4) (1- \la^8) (1- \la^{12})}$$
\beql{EqG66b}
R: \mbox{omit $D$ from \eqn{EqG66a}}
\eeq
(This ring has also been studied in \cite{DHS97}.)
\paragraph{cwe of $W^{(1)} - W^{(3)}$ with $1^n \in C$:}
$G= \langle M_4, \beta_4, \gamma_4 \rangle$, order 6144, with character $\chi (M_4) = i^n$, $\chi (\beta_4)=1$,
$\chi( \gamma_4) = \eta^n$;
$\ker(\chi)$ has order 3072
$$\Phi = \frac{1+ \la^{32}}{(1- \la^4) (1- \la^8) (1- \la^{12}) (1- \la^{16})}$$
\beql{EqG66ca}
R: \frac{1,~ A^2 B^2 C^2 D^2 (A^4 + B^4) (A^4 +C^4) (A^4 -D^4) (B^4 - C^4) (B^4 + D^4) (C^4 +D^4)}{\mbox{symmetric polynomials in $A^4 , -B^4, -C^4, D^4$}}
\eeq
\paragraph{swe of $W^{(1)} - W^{(3)}$ with $1^n \in C$:}
$$\Phi = \frac{1}{(1- \lambda^4 ) (1- \la^8 ) (1- \la^{12})}$$
$$
R: \frac{1}{\mbox{symmetric polynomials in $A^4, ~-B^4 , ~ -C^4$}}
$$
\subsection{Family \EeightII: Type II self-dual codes over $\ZZ_4$}\label{GL9}
\paragraph{cwe of $C$, $1^n \in C$}
\cite{Me209}, \cite{BSBM97}
In view of the remarks following Theorem~\ref{GP6},
this is not a severe restriction.
$$G= \langle M_4, \beta_4, \gamma_4 \rangle ~,\quad
{\rm order~6144}
$$
$$
\Phi = \frac{(1+ \la^{16}) (1+ \la^{32})}{(1-\la^8 )^2 (1- \la^{16}) (1- \la^{24})}$$
\beql{EqG69a}
R: \frac{(1,~ f_{16}) \times (1,~ f_{32})}{A^8 + B^8 + C^8 + D^8 , ~ f_8 , ~ A^{16} + \cdots + D^{16} , ~ A^{24} + \cdots + D^{24}}
\eeq
where
\begin{eqnarray*}
f_8 & = & A^4 C^4 + C^4 D^4 + D^4 B^4 + B^4 A^4 - A^4 D^4 - B^4 C^4 ~, \\
f_{16} & = & (ABCD)^4 ~, \\
f_{32} & = & (ABCD)^2 (A^4 + C^4) (C^4 + D^4) (D^4 + B^4) (B^4 + A^4) (A^4 -D^4) (B^4 -C^4)
\end{eqnarray*}
\paragraph{swe of $C$, $\pm 1^n \in C$}
$$\Phi = \frac{1+ \la^{16}}{(1- \la^8)^2 (1- \la^{24})}$$
\beql{EqG71a}
R: \frac{1,~\theta_{16}}{\theta_8,~ h_8,~ \theta_{24}}
\eeq
where
\begin{eqnarray*}
\theta_8 & = &
x^8 + 28 x^6 z^2 + 70 x^4 z^4 + 28 x^2 z^6 + z^8 + 128y^8 ~, \\
\theta_{16} & = & \{ x^2 z^2 (x^2 + z^2)^2 -4y^8 \}
\{ (x^4 + 6x^2 z^2 + z^4)^2 - 64y^8 \} ~, \\
\theta_{24} & = & y^8 ( x^2 - z^2)^8 ~, \\
h_8 & = & \{xz (x^2 + z^2 )-2y^4 \}^2 ~.
\end{eqnarray*}
\paragraph{cwe of $C$, $1^n \in C$, Lee weights divisible by 4}
(\cite{Me209})
$$G = \left\langle M_4, \beta_4, \gamma_4 , \left(
\matrix{
0& 0 & 1 & 0 \cr
0& i & 0 & 0 \cr
1 & 0 & 0 & 0 \cr
0 & 0 & 0 & i \cr
}
\right) \right\rangle
$$
(Shephard \& Todd \#2a), order 49152
$$\Phi = \frac{1}{(1- \la^8 ) (1- \la^{16} )^2 (1- \la^{24} )}$$
\beql{EqG71aa}
R: \frac{1}{f_{16} ,~\mbox{symmetric polynomials in $A^8$, $B^8$, $C^8$, $D^8$}}
\eeq
\paragraph{swe of $C$, $\pm 1^n \in C$, Lee weights divisible by 4}
(\cite{Me209})
(Set $t=y$ in cwe)
$$\Phi = \frac{1}{(1- \la^8 ) (1- \la^{16} ) (1- \la^{24} )}$$
\beql{EqG71ab}
R: \frac{1}{\theta_8 ,~ \theta_{16},~ \theta_{24}}
\eeq
However, the following result shows that the extra condition on the Lee weights
may not be a good thing.
For it was shown in \cite{Me184} that
most interesting linear codes over $\ZZ_4$ do not have linear images
under the Gray map.\index{Gray map}
\begin{theorem}\label{thLee}
{\rm \cite{Me209}}
If $C$ is a self-dual code over $\ZZ_4$ with all Lee weights divisible by 4,
then the binary image of $C$ under the Gray map \eqn{EqLift2}
is linear.
\end{theorem}

\noindent
For the proof, see \cite{Me209}.

\paragraph{Codes.}
The following codes will be used:
$i_1$ and $\sD_4^\oplus$ are defined in Section~\ref{GU8},
and $o_8$ is the octacode \eqn{Eq19a}.
$\sJ_{10}$ is the self-dual code with generator matrix
$$
\left[
\matrix{
1&0&1&1&1&1&1&0&1&1 \cr
0&1&0&0&0&0&3&3&1&0 \cr
0&0&2&0&0&0&0&0&0&2 \cr
0&0&0&2&0&0&0&0&0&2 \cr
0&0&0&0&2&0&0&0&0&2 \cr
0&0&0&0&0&2&0&0&0&2 \cr
0&0&0&0&0&0&2&0&2&0 \cr
0&0&0&0&0&0&0&2&2&2 \cr
}
\right]
$$
and $| \sJ_{10}| = 4^2 2^6$ (\cite{Me168}).
$\sJ_{16}$ has generator matrix
$$
\left[
\matrix{
1&0&0&0&0&0&1&1&1&0&3&3&1&0&3&2 \cr
0&1&0&0&0&0&1&0&0&1&3&3&1&1&2&3 \cr
0&0&1&0&0&0&1&0&0&0&0&2&2&0&3&3 \cr
0&0&0&1&0&0&0&1&1&1&3&0&2&3&3&1 \cr
0&0&0&0&1&0&0&1&1&1&0&0&1&1&1&1 \cr
0&0&0&0&0&1&0&0&0&0&0&3&2&0&1&1 \cr
0&0&0&0&0&0&2&0&0&0&0&2&0&2&0&2 \cr
0&0&0&0&0&0&0&2&0&0&0&0&0&2&2&2 \cr
0&0&0&0&0&0&0&0&2&0&0&0&0&2&2&2 \cr
0&0&0&0&0&0&0&0&0&2&0&0&0&2&0&0 \cr
}
\right]
$$
and $|\sJ_{16} | = 4^6 2^4$.

$\sK_{4m}$ $(m \ge 1$, but note that $\sK_4 \cong \sD_4^\oplus$) is a self-dual
code introduced by Klemm \cite{Klemm89},
having generator matrix
\beql{EqG71b}
\left[
\matrix{
1&1&1& \ldots &1&1 \cr
0&2&0& \ldots &0 & 2 \cr
0 & 0 & 2 & \ldots & 0 & 2 \cr
. & . & . & \ldots & . & . \cr
0 & 0 & 0 & \ldots & 2 & 2 \cr
}
\right] ~;
\eeq
$| \sK_{4m} | = 4^1 2^{4m-2}$, $g=2^{4m-1} (4m)!$,
$cwe = (A^{4m} + B^{4m} + C^{4m} + D^{4m} )/2$ (see \eqn{EqG60b}).
\beql{EqG72}
\mbox{\begin{tabular}{ll} \hline
Ring & Codes \\ \hline
\eqn{EqG60a} & $i_1$, $\sD_4^\oplus$ in 2 versions \eqn{EqC23a}, $o_8$; $\sJ_{10}$ \\
\eqn{EqG61a} & $i_1 , \sD_4^\oplus$, $o_8$ \\
\eqn{EqG62a} & $i_1 , \sD_4^\oplus$; $o_8$ \\
\eqn{EqG63a} & $\sK_4$, $\sK_8$, $o_8$, $\sK_{16}$; $\sK_{12}, \sJ_{16}$ \\
\eqn{EqG64a} & $\sK_4$, $\sK_8$, $o_8$; $\sK_{12}$ \\
\eqn{EqG64b} & $\sK_4$, $o_8$; $\sK_8$, $\sK_{12}$ \\
\eqn{EqG69a} & $\sK_8$, $o_8$, $\sK_{16}$, $\sK_{24}$; $\sJ_{16}$,? \\
\eqn{EqG71a} & $\sK_8$, $o_8$, $\sK_{24}$; $\sJ_{16}$ \\
\eqn{EqG71ab} & $\sK_8$, $\sK_{16}$, $\sK_{24}$
\end{tabular}
}
\eeq
Again the question mark indicates that we do not have a satisfactory code to produce the desired polynomial.
\subsection{Family \Enine: Self-dual codes over $\ZZ_m$}\label{GL10}
\hsp
The Hamming weight enumerator of
a self-dual code over $\ZZ_m$ for general
$m$ has been considered in \cite{HaOu97}.
\section{Weight enumerators of maximally self-orthogonal codes}\label{SO}\index{maximally self-orthogonal code}
\hsp
In some cases it is possible to prove results analogous to those in Section~\ref{GL} for codes which are maximally self-orthogonal yet not self-dual,
the $[7,3,4]$ Hamming code $e_7$\index{Hamming code!$e_7$} with weight enumerator $p_7 = x^7 + 7x^3 y^4$ being a typical example.
A more trivial example is the zero code $z_1 = \{0\}$, with weight enumerator $p_1 =x$.

The following results are proved in \cite{Me35}.

For $n$ odd, let $C$ be an $[n, \frac{1}{2} (n-1)]$ self-orthogonal binary code.
Thus $C^\perp = C \cup (1+C)$.
The weight enumerator of $C$ belongs to the module
$R= p_1 \CC [x^2 + y^2, x^2 y^2 (x^2 - y^2)^2 ] \oplus p_7 \CC [ x^2 + y^2 , x^2 y^2 (x^2 - y^2) ^2 ]$,
which in the notation of the previous section would be described by
$$\Phi = \frac{\la + \la^7}{(1- \la^2)(1- \la^8)}$$
\beql{EqS01}
R: \frac{p_1,~ p_7}{x^2 + y^2,~ x^2y^2 (x^2 - y^2)^2}
\eeq
(compare \eqn{EqG6a}).
If in addition $C$ is doubly-even, the module is described by:

\noindent
(a)~if $n=8m-1$,
$$\Phi = \frac{\la^7 + \la^{23}}{(1- \la^8)(1- \la^{24})}$$
\beql{EqS02}
R: \frac{p_7 ,~ p_{23}}{x^8 + 14x^4 y^4 + y^8,~ x^4 y^4 (x^4 - y^4 )^4}
\eeq
(b)~if $n=8m+1$,
$$
\Phi = \frac{\la + \la^{17}}{(1- \la^8 ) (1- \la^{24})}
$$
\beql{EqS03}
R: \frac{p_1,~ p_{17}}{x^8 + 14x^4 y^4 + y^8,~ x^4 y^4 (x^4 - y^4)^4}
\eeq
(compare \eqn{EqG5b}).
Here $p_{17} = x^{17} + 17x^{13} y^4 + 187x^9 y^8 + 51 x^5 y^{12}$,
$p_{23} = x^{23} + 506 x^{15} y^8 + 1288 x^{11} y^{12} + 253 x^7 y^{16}$.

\paragraph{Codes.}  The code $g_{23}$ is the cyclic version of $g_{24}$ obtained by deleting any coordinate.\index{Golay code!cyclic ($g_{23}$)}\index{code!Golay}

\begin{center}
\begin{tabular}{ll} \hline
Ring & Codes \\ \hline
\eqn{EqS01} & $i_2 , e_8 ; z_1 , e_7$ \\
\eqn{EqS02} & $e_8 , g_{24} ; e_7 , g_{23}$ \\
\eqn{EqS03} & $e_8 , g_{24} ; z_1 , (d_{10} e_7)^+$
\end{tabular}
\end{center}

There are analogous results for ternary codes:
see \cite{Me80}.


\section{Upper bounds}\label{BD}\index{upper bounds}
\hsp
Of course, we are interested
not just in codes per se, but
also in good (or, at the very least, interesting) codes, that is, codes
with large minimal distance (Hamming, Lee, or Euclidean, as appropriate).
In order to know if a particular code is good, it is necessary to know how good
comparable
codes could be; that is, for a given length and dimension, what is the
optimal minimal distance?  For general codes, this question was studied in
Chapters~xx (Levenshtein), yy (Brouwer) and zz (Litsyn);
we are, of course, interested in self-dual codes.  As one
might imagine, the constraint of self-duality usually leads to stronger
bounds.

We will concentrate most of our attention on binary codes (family \Eone),
pointing out analogues to other families as they arise.

Essentially all of the bounds we will be discussing are special cases of
the linear programming\index{linear programming bound} (or LP) bound (Section 2.5 of Chapter~yy (Brouwer)); that is, they rely
\index{bound!linear programming}\index{LP bound}
on the fact that both the weight enumerator of the code and the weight
enumerator of its dual are nonnegative.  For a self-dual code, these weight
enumerators are, of course, equal.  
So for Type II self-dual binary
codes, for instance, we have the following:
\begin{theorem}\label{thBD1}
If there exists a Type II self-dual binary code of length $n$ and minimal
distance $d$, then there exists a homogeneous polynomial $W(x,y)$ with
nonnegative (integer) coefficients such that
\begin{eqnarray*}
2^{n/2} W(x+y,x-y) &=& W(x,y)\\
W(1,y) &=& 1+O(y^d)\\
W(x,iy) &=& W(x,y).
\end{eqnarray*}
\end{theorem}
These conditions assert that the code is self-dual, that it has
minimal distance $d$, and that it is of Type II, respectively.

The analogues for other classes of codes should be clear; in each case, the
appropriate enumerator (Hamming, symmetrized, complete) is nonnegative,
invariant under the appropriate transformations (see Section \ref{GL}), and is zero
on all terms of low weight.  In some cases, we can add further constraints
from shadow\index{shadow} theory (Section \ref{Shad}), since the weight enumerator of the
shadow of the code is also nonnegative.  For instance:
\begin{theorem}\label{thBD2}
If there exists a Type I self-dual binary code of length $n$ and minimal
distance $d$, then there exist homogeneous polynomials $W(x,y)$ and
$S(x,y)$ with nonnegative (integer) coefficients such that
\begin{eqnarray*}
W(x,y) & = & 2^{-n/2} W(x+y , x-y) \\
W(1,y) &=& 1+O(y^d)\\
S(x,y) &=& 2^{-n/2} W(x+y,i(x-y)).
\end{eqnarray*}
\end{theorem}
Again there
are analogues for each family for which shadows are
well-defined (\Eone, \Efive, \Eeight).

\paragraph{Remark.}  For a code $C$ from family \Esix\ (linear over $F_q$, $q$ a square, with
Hermitian inner product), it can be shown that the polynomial
$$
S(x,y)=q^{-n/2} W((\sqrt{q}-1) x + (\sqrt{q}+1) y,~ y-x)
$$
has nonnegative (but not necessarily integral) coefficients; note that this agrees with
the shadow enumerator for $q=4$.  This can be used to strengthen the LP
bound in those cases.  The known proof that this is nonnegative involves
constructing a quantum code\index{quantum code}\index{code!quantum}
$Q$ from $C$ (\cite{RainsNB}); $S(x,y)$ is then the
shadow enumerator of $Q$ (\cite{RainsQWE}, proved nonnegative in \cite{RainsPI}).
There is surely a more direct proof.

One way to apply the linear programming bound is to ignore the constraint
that the coefficients of $W(x,y)$ be nonnegative, and simply ask that the
low order coefficients be as specified.  This gives a surprisingly good
bound for Type II binary codes.
Recall from Theorem~\ref{ThMS3c} that for $C$ of Type II, $W(x,y)$
lies in the ring
$$
R={\Bbb C}[x^8+14 x^4 y^4+y^8, x^4 y^4(x^4-y^4)^4],
$$
and
if $C$ has length
$n$, $W(x,y)$ has degree $n$.  The subspace of
$R$ of degree $n$ has dimension $D =[{n\over 24}]+1$.  This lets us set
the first $D$ coefficients of $W(x,y)$ arbitrarily; in particular, there
exists a unique element $W^*(x,y)$ of $R$ such that $W^*(1,y)=1+O(y^{4D})$.
This is known as the {\it extremal} enumerator, since $W^*$ has the largest
\index{extremal weight enumerator}
minimal distance of any Type II self-dual enumerator.  It follows
immediately that the minimal distance of any Type II code of length $n$ is
bounded above by the minimal distance of $W^*$.

\begin{theorem}\label{thBD3}
{\rm \cite{Me30}}
The first nonzero coefficient of $W^*(1,y)$ occurs
precisely at degree $4D$; in
particular, the minimal distance of a Type II self-dual binary code of
length $n$ is at most $4[n/24]+4$.
\end{theorem}

In fact it is possible to use the B\"{u}rmann-Lagrange theorem\index{B\"{u}rmann-Lagrange theorem}\index{theorem!B\"{u}rmann-Lagrange} (Theorem~\ref{thBD6}) to derive an explicit formula for the number of words of weight $4 D$ in the extremal enumerator.
Let $\mu = [n/24]$, so that $D= \mu+1$.
Then we have
\begin{theorem}\label{thMS13}
{\rm (Mallows and Sloane \cite{Me30}.)}\index{Mallows and Sloane bound}\index{bound!Mallows and Sloane}
$A_{4 \mu + 4}^\ast$, the number of codewords of minimal nonzero
weight $4D= 4 \mu +4$ in the extremal weight enumerator, is given by:
\beql{EqMS54}
{\binom{n}{5}} {\binom{5 \mu -2}{\mu -1}} \Bigl/ {\binom{4 \mu +4}{5}} ,
\quad\mbox{if} \quad n = 24 \mu ~,
\eeq
\beql{EqMS55}
\frac{1}{4} n(n-1) (n-2) (n-4) \frac{(5 \mu )!}{\mu ! (4 \mu + 4)!} , \quad\mbox{if}\quad n = 24 \mu +8 ~,
\eeq
\beql{EqMS56}
\frac{3}{2} n(n-2) \frac{(5 \mu +2)!}{\mu !(4 \mu +4)!} , \quad\mbox{if} \quad
n=24 \mu +16 ~,
\eeq
and is never zero.
\end{theorem}

For the proof, see \cite{Me30} or \cite{MS}, Chapter 19.
There is a similar formula for Type I binary codes --- see \cite{MS}, Chapter~19, Problem (12).

Results similar to Theorem~\ref{thBD3} hold for other families:

\begin{theorem}\label{thBD4}
The minimal distance of a Type I binary self-dual code is at most $2[n/8]+2$.
The minimal distance of a Type II binary self-dual code is at most $4[n/24]+4$.
The minimal distance of a self-dual code from family \Etwo\ is at most $3[n/12]+3$.
The minimal distance of a self-dual code from family \Ethree\ is at most $2[n/6]+2$.
The minimal distance of a Type II self-dual code from family \Efive\ is at most
$2[n/6]+2$.
The minimal distance of a self-dual code from families \Efour,
\Efive, \Esix\ or \Eseven\ is at most
$[n/2]+1$.
\end{theorem}

Note that the last bound is simply the Singleton bound,\index{Singleton bound}\index{bound!Singleton}
obtained from the
ring ${\Bbb C}[x^2+(q-1)y^2,y(x-y)]$ of \eqn{EqG37cb}.
As we have already remarked in Section~\ref{GL6}, this is
not the correct ring (that is, the smallest ring containing all Hamming
enumerators of self-dual codes).  In some cases ($q=4$ or $q=5$), we know a
smaller ring; however, since the ring is no longer free, it is
much more difficult to use.  In particular, it is no longer the case that
we may
set the leading coefficients arbitrarily.  This
leads to the extremal enumerator not being unique, making it difficult to
determine its first nonzero coefficient.  Similarly, any attempt to
make an analogous argument for families \Eeight\ or \Enine\ will have the problem that, in
those cases, we are primarily interested in Lee weight
or Euclidean norm,
forcing us to work with the symmetrized weight enumerator.  This is, of
course, much more difficult to deal with than the Hamming enumerator.
A partial solution to this problem is provided by Theorem~\ref{thBD9} below.

In each case it can be shown (cf. \cite{Me41}) that the bounds
of Theorems~\ref{thBD3} and \ref{thBD4}
can be met for at most
finitely many $n$:
in fact, the next coefficient $(A_{4\mu+8}^\ast )$
{\it after} the leading nonzero coefficient in the extremal enumerator
becomes negative for sufficiently large $n$.
\index{negative coefficients exist}
Furthermore, for any constant $\alpha$,
the minimal distance can be within $\alpha$ of the bound only finitely often.
For Type II binary codes,
for instance,
it was shown in \cite{Me41} that
the $A_{4n+8}^\ast$ term first goes negative when $n$ is around 3720. 
Ma and Zhu \cite{MZ97} and Zhang \cite{Zha97}
have recently determined precisely when the $A_{4n+8}^\ast$ term first goes
negative, and have obtained similar results for several other families.
The following result incorporates the work of several authors.
\begin{theorem}\label{thBDa}
{\rm \cite{Zha97}}
Let
$C$ be a self-dual code of length $n$ from one of the families \EoneI, \EoneII,
\Etwo, \Ethree;
and let $c=2,4,3,2$, respectively, and $\mu = [n/8]$,
$[n/24]$, $[n/12]$, $[n/6]$.
Then the coefficient $A_{c(\mu+2)}^\ast$ in the extremal Hamming weight enumerator
is negative if and only if:
$$
\begin{array}{ll}
\mbox{$(\EoneI)$:} &
n=8i\,(i \ge 4),~ 8i+2\,(i \ge 5),~ 8i+4\,(i \ge 6) ,~ 8i+6 \,(i \ge 7); \\ [+.1in]
\mbox{$(\EoneII)$:} & n= 24i \,( i \ge 154) ,~ 24i + 8 \,(i \ge 159),~ 24i + 16 \,(i \ge 164); \\ [+.1in]
\mbox{$(\Etwo)$:} & n=12i \,(i \ge 70),~ 12i + 4\,(i \ge 75) ,~ 12i + 8 \,(i \ge 78); \\ [+.1in]
\mbox{$(\Ethree)$:} & n= 6i \,( i \ge 17),~ 6i + 2\,(i \ge 20) ,~ 6i+4\,( i \ge 22).
\end{array}
$$
In particular, the first time $A_{4 \mu +8}^\ast$ goes negative for Type II
codes is at $24 \times 154 = 3696$.
\end{theorem}

Of course other coefficients in the extremal weight enumerator may go negative
before this.
In the case of ternary self-dual codes, for example, family \Etwo, the extremal
Hamming weight enumerator contains a negative coefficient for lengths 72, 96, 120 and all $n \ge 144$.

The best asymptotic bound presently known for Type II codes is the following.
\begin{theorem}\label{thBD4b}
{\rm (Krasikov and Litsyn \cite{KrLit97}.)}\index{Krasikov
and Litsyn bound}\index{bound!Krasikov and Litsyn}
The minimal distance $d$ of a Type II binary code of length $n$ satisfies
$$d \le 0.166315 \ldots n + o(n) , \quad n \to \infty ~.$$
The constant in this expression is the real root of $8x^5 - 24 x^4 + 40x^3 - 30x^2 + 10x -1$.
\end{theorem}
The proof uses a variant of the linear programming bound.

For Type I binary codes,
the bound of Theorem~\ref{thBD4} is especially weak.
Ward\index{Ward, H. N.} \cite{Wa76} has shown
that the minimal distance can be $2[n/8]+2$ precisely when $n$ is one of
2, 4, 6, 8, 12, 14, 22 or 24.  This suggests that the bound can be greatly
strengthened, which is indeed the case.  Conway and Sloane\index{Conway and Sloane bound}\index{bound!Conway and Sloane} \cite{Me158}
showed that $d \le 2 [(n+6)/10]$ for $n > 72$, and Ward (\cite{Wa92}, see also Chapter ``Ward'') established $d \le n/6 + O( \log n )$.
It turns out, in fact, that the ``correct'' bound is $4[n/24]+4$ (except when
$n+2$ is a multiple of 24), just as for Type II codes.  The key to proving
this fact is the observation that we have not yet used the shadow\index{shadow}
enumerator.

\begin{theorem}\label{thBD5}
{\rm (Rains\index{Rains bound}\index{bound!Rains} \cite{RainsSB}.)}
Suppose $C$ is a $[n,n/2,d]$ self-dual binary code.  Then $d\le 4[n/24]+4$,
except when $n \equiv 22$ $(\bmod~24)$, when $d\le 4[n/24]+6$.  If $n$ is a multiple of 24,
any code meeting the bound is of Type II.
If $n \equiv 22$ $(\bmod~24)$, any code meeting the bound can be
obtained by shortening a Type II code of length $n+2$ that also meets the bound.
\end{theorem}

\paragraph{Proof (sketch).}
From \eqn{EqG6a}, $W(x,y)$ lies in the ring ${\Bbb
C}[x^2+y^2,x^2y^2(x^2-y^2)^2]$; consequently we can write
\begin{eqnarray*}
W(1,y)&=&\sum_j a_j y^{2j}\\
      &=&\sum_i c_i (1+y^2)^{n/2-4i} (y^2(1-y^2)^2)^i.
\end{eqnarray*}
Applying the shadow transform, we have
\begin{eqnarray*}
S(1,y)&=&\sum_j b_j y^{2j+t}\\
      &=&\sum_i c_i (2y)^{n/2-4i} (-(1-y^4)/2)^i,
\end{eqnarray*}
where $t=((n/2)\bmod 4)$.  Suppose $C$ had minimal distance $4[n/24]+6$.
This fact determines $c_i$ for $0\le i\le 2[n/24]+2$, and in particular
$c_{2[n/24]+2}$.  On the other hand, we can also express $c_{2[n/24]+2}$ as
a linear combination of the $b_j$ for small $j$.  It turns out that these
two expressions for $c_{2[n/24]+2}$ are incompatible; in particular, we
find that a certain nonnegative linear combination of the $b_j$ is
negative.

Rather than give the (somewhat messy) details of the proof, we will simply
show how one can compute the coefficients in these linear combinations.
This uses the B\"{u}rmann-Lagrange theorem:\index{B\"{u}rmann-Lagrange theorem}\index{theorem!B\"{u}rmann-Lagrange}

\begin{theorem}\label{thBD6}
{\rm (B\"urmann-Lagrange.)}
Let $f(x)$ and $g(x)$ be formal power series, with $g(0)=0$ and $g'(0)\ne
0$. If coefficients $\kappa_{ij}$ are defined by
$$
x^j f(x)=\sum_{0\le i} \kappa_{ij} g(x)^i,
$$
then
$$
\kappa_{ij}={1\over i}[\mbox{coeff. of $x^{i-1}$ in $[j x^{j-1} f(x)+x^j
f'(x)] \left({x\over g(x)}\right)^i$}].
$$
\end{theorem}
For proof and generalizations, see
\cite[p.~133]{MS1414},
\cite{MS533}, \cite{MS1137}, \cite{MS1138}, \cite{MS1139}.

For instance, to compute $c_{2[n/24]+2}$, we note that
$$
\sum_i c_i (1+y^2)^{n/2-4i} (y^2(1-y^2)^2)^i
=
1+O(y^{4[n/24]+6}).
$$
Dividing both sides by $(1+y^2)^{n/2}$ and substituting $y=\sqrt{Y}$, we get:
$$
\sum_i c_i \left({Y(1-Y)^2\over (1+Y)^4}\right)^i
=
(1+Y)^{-n/2} + O(Y^{2[n/24]+3}).
$$
We can then apply B\"urmann-Lagrange, with
$$
f(Y)=(1+Y)^{n/2}, \quad
g(Y)=Y (1-Y)^2 (1+Y)^{-4}
$$
to obtain
\begin{eqnarray*}
c_i &=& {1\over i}
[\mbox{coeff. of $Y^{i-1}$ in $[{d\over dY} (1+Y)^{-n/2}]
\left((1+Y)^4 (1-Y)^{-2}\right)^i$}]\\
&=& {-n\over 2i}
[\mbox{coeff. of $Y^{i-1}$ in $(1+Y)^{-n/2-1+4i} (1-Y)^{-2i}$}]\\
&=& {-n\over 2i}
[\mbox{coeff. of $Y^{i-1}$ in $(1+Y)^{-n/2-1+6i} (1-Y^2)^{-2i}$}].
\end{eqnarray*}
In particular, for $i=2[n/24]+2$,
$$
c_{2[n/24]+2} 
=
{-n\over 4[n/24]+4}
[\mbox{coeff. of $Y^{2[n/24]+1}$ in $(1+Y)^{-n/2+12[n/24]+11} (1-Y^2)^{-4[n/24]-4}$}].
$$
It follows that $c_{2[n/24]+2}\le 0$, with equality only when $n \equiv 22$ $(\bmod ~24)$,
since all coefficients of any power series of the form
$(1+Y)^a(1-Y^2)^{-b}$ are positive whenever $a$, $b>0$.

Similarly, we find that the coefficients of the expansion of
$c_{2[n/24]+2}$ in terms of the $b_j$ are positive.  This proves the
bound, except when $n \equiv 22$ $(\bmod ~24)$; the proof that the bound holds in that
case and that a code meeting the bound is even if $n \equiv 0$ $(\bmod~24)$
is left to the reader.~~~$\bsq$

This bound agrees with the full linear programming bound for $n\le 200$,
and, most likely, for much larger $n$.  However, it is likely that again
it can only be attained for finitely many $n$.

There is also an analogue of this bound for Type I codes from family \Efive.

\begin{theorem}\label{thBD7}
If $C$ is an additive self-dual code of length $n$ and minimal distance
$d$ from family \Efive,
then $d\le 2[n/6]+2$, except when $n \equiv 5 ~(\bmod~6))$, when $d\le 2[n/6]+3$.  If
$n$ is a multiple of 6, then any code meeting the bound is even.
\end{theorem}

We will call a code {\em extremal}\index{extremal code}\index{code!extremal} if it meets the strongest of the applicable
bounds from Theorems~\ref{thBD4}, \ref{thBD5}, and \ref{thBD7}.  For Type II binary codes, ternary
codes, and linear codes over $GF(4)$ this agrees with the historical
usage.  For Type I binary codes, however, ``extremal'' has generally
been used to mean a code meeting the much weaker bound of Theorem~\ref{thBD4};
in the light of Theorem~\ref{thBD5}, it seems appropriate to change the definition.

Concerning codes over $\ZZ_4$, Bonnecaze, Sol\'{e}, Bachoc and Mourrain \cite{BSBM97}\index{Bonnecaze, Sol\'{e}, Bachoc and Mourrain bound}\index{bound!Bonnecaze, Sol\'{e}, Bachoc and Mourrain} show:
\begin{theorem}\label{thBD9}
Suppose $C$ is a Type II self-dual code over $\ZZ_4$ of length $n$.
Then the minimal Euclidean norm of $C$ is at most
\beql{EqDB9a}
8 \left[\frac{n}{24} \right]+8 ~.
\eeq
\end{theorem}

The proof uses $C$ to define an even unimodular $n$-dimensional lattice $\La (C) = \{\frac{1}{2} u \in \RR^n : u ~(\bmod~4) \in C\}$, and examines its theta series.

As usual, one can derive an analogue for Type I codes:
\begin{theorem}\label{thBD9a}
{\rm \cite{RaSl97}}
Suppose $C$ is a Type I self-dual code over $\ZZ_4$ of length $n$.
The minimal Euclidean norm of $C$ is at most
\beql{EqBD9b}
8 \left[ \frac{n}{24} \right] +8 ~,
\eeq
except when $n \equiv 23$ $(\bmod~24)$, in which case the bound is
\beql{EqBD9c}
8 \left[ \frac{n}{24} \right] + 12 ~.
\eeq
If equality holds in \eqn{EqBD9c} then $C$ is a shortened version of a Type II code of length $n+1$.
\end{theorem}

We say that codes meeting either of these bounds are {\em norm-extremal}.\index{norm-extremal}
\index{code!norm-extremal}
For Type II codes this agrees with the definition given in \cite{BSBM97}.

There should be an analogous concept of {\em Lee-extremal}, but at present we do not know what this is.
\index{code!Lee-extremal}
Of course, the bounds \eqn{EqBD9b} and \eqn{EqBD9c} also apply to Lee weight.
But this is not a satisfactory bound, since it is not even tight at length 24,
where the highest attainable Lee weight is 12 rather than 16 (see Table~\ref{TXZ4}).

The fact that, from Theorem~\ref{thBD5}, an extremal binary code of length a
multiple of 24 must be doubly-even suggests that these codes are likely to
be particularly nice.  Indeed, we have the following
result, which is a consequence of the Assmus-Mattson theorem\index{Assmus-Mattson theorem}\index{theorem!Assmus-Mattson} (see \cite[Chap.~6]{MS},
Theorem 11.14 of Chapter~1, Section~5 of Chapter xx (Tonchev)).

\begin{theorem}\label{thBD8}
Let $C$ be an extremal binary code of length $24m$.  Then
the codewords of $C$ of any given weight form a 5-design.
\end{theorem}

Similarly, the supports of the minimal codewords of an extremal ternary
code of length $12m$ form a 5-design.
For codewords of larger weight, the natural incidence structure is {\it
almost} a 5-design, except that it may have repeated blocks.  Similarly,
for an extremal additive code over $\FF_4$ of length $6m$, the supports
with multiplicities of the codewords of any fixed weight form a 5-design
with repeated blocks.
Harada \cite{Hara99} has shown that the $\ZZ_4$-lift of the Golay code
$g_{24}$ also yields 5-designs.
More generally, one can show that the words of any fixed symmetrized type,
in any of the 13 Lee-optimal self-dual codes of length 24 over $\ZZ_4$,
form a colored 5-design, possibly with repeated blocks \cite{BRSZ4}.
See also \cite{Oze96}.

\section{Lower bounds}\label{Lbounds}\index{lower bounds}
\hsp
There are two ways to obtain lower bounds on the optimum minimal distance
of a code of length $n$.  The first way, naturally, is simply to construct
a good code.
Just as for general linear codes, there is also a nonconstructive lower bound,
analogous to the Gilbert-Varshamov bound\index{Gilbert-Varshamov bound}\index{bound!Gilbert-Varshamov} (cf. Theorems~3.1, 3.4, 3.5 of Chapter~1).

We first consider the case of self-dual binary codes (family \Eone).

\begin{theorem}\label{GV1}
{\rm \cite{Thomp1}, \cite{Me24}}
Let $n$ be any positive even integer.  Let $d_{GV}$ be the largest
\index{good self-dual codes exist}
integer such that
\beql{GV1a}
\sum_{0<i<d\atop 2|i} {n\choose i} < 2^{n/2-1}+1.
\eeq
Then there exists a self-dual binary code of length $n$ and minimal
distance at least $d_{GV}$.
\end{theorem}

\paragraph{Proof}
If we can show that the expected number of nonzero vectors of weight less
than $d_{GV}$ in a {\it random} self-dual code of length $n$ is less
than 1, it will immediately follow that there exists {\it some} self-dual
code of length $n$ with no such vectors.

Let us therefore compute the average weight enumerator\index{average weight enumerator} of the set of
\index{weight enumerator!average}
self-dual codes.  Consider the group $G$ of binary matrices that preserve
the quadratic form $I$.  On the vector space of even weight vectors, modulo
the all 1's vector, the quadratic form becomes symplectic, and the group
acts as the full symplectic group.  In particular, it is therefore
transitive on nonzero vectors of even weight, modulo $1^n$.  It follows
that the expected number of vectors of weight $2i$ in a random code
must be proportional to ${n\choose 2i}$, except for $i=0$ or $i=n/2$.
Thus the average weight enumerator has the form:
\begin{eqnarray*}
\overline{W}(x,y) &=& a x^n+b \sum_{1\le i\le {n/2}-1} {n\choose 2i} x^{n-2i}
y^{2i} + c y^n\\
&=&
a x^n+c y^n + b ({1\over 2}(x+y)^n+{1\over 2}(x-y)^n-x^n-y^n).
\end{eqnarray*}
Since every self-dual binary code contains the 0 vector and the all 1's
vector, $\overline{W}(1,0)=\overline{W}(0,1)=1$; since every self-dual
code contains a total of $2^{n/2}$ vectors, $\overline{W}(1,1)=2^{n/2}$.
Solving for $a$, $b$, and $c$, we find:
$$
\overline{W}(x,y) = x^n + y^n + {1\over 2^{n/2-1}+1}
\sum_{1\le i\le {n/2}-1} {n\choose 2i}
x^{n-2i} y^{2i} ~.
$$
Thus the average number of nonzero vectors of weight less than $d$ is
$$
{1\over 2^{n/2-1}+1}
\sum_{0<i<d\atop 2|i} {n\choose i}.~~~\hfill \bsq
$$

\begin{coro}\label{GV2}
{\rm \cite{Thomp1}, \cite{Me24}}
There exists an infinite sequence of self-dual $[n_i,n_i/2,d_i]$ binary
codes, such that $n_i$ tends to infinity, and
$$
\liminf_{i\to\infty} {d_i\over n_i} \ge \delta,
$$
where $\delta\sim .11002786$ is the unique solution less than $1\over 2$ of
$$
H_2(\delta)=-\delta\log_2(\delta)-(1-\delta)\log_2(1-\delta)={1\over 2}.
$$
\end{coro}

\paragraph{Proof.}
Take the logarithm of both sides of \eqn{GV1a}, divide by $n$, and let $n$
tend to infinity.  The resulting inequality is
$$
H_2(\delta)\le {1\over 2},
$$
as desired.~~~$\bsq$

Similar results hold if one restricts ones attention to codes of Type II:

\begin{theorem}\label{GV3}
{\rm \cite{Thomp1}, \cite{Me24}}
Let $n$ be any positive multiple of 8.  Let $d_{GV}$ be the largest
integer such that
\beql{GV3a}
\sum_{0<i<d\atop 4|i} {n\choose i} < 2^{n/2-2}+1
\eeq
Then there exists a doubly-even self-dual binary code of length $n$ and
minimal distance at least $d_{GV}$.
\end{theorem}

\paragraph{Proof.}
Again we compute the average weight enumerator.  The key observation is
that the function ${1\over 2}wt(v)$ induces a quadratic form on the space
of even weight vectors modulo the all 1's vector.  The group of matrices
that preserve this quadratic form is transitive on the kernel of this
quadratic form; that is, vectors of weight divisible by 4, modulo $1^n$.
This allows us to write down the average weight enumerator:
$$
\overline{W}_{II}(x,y)
=
x^n+y^n+{1\over 2^{n/2-2}+1} \sum_{0<i<n/4} {n\choose 4i} x^{n-4i}
y^{4i}.
\eqno{\bsq}
$$

Asymptotically, this agrees with Corollary~\ref{GV2} (as well as 
the Gilbert-Varshamov bound\index{Gilbert-Varshamov bound}\index{bound!Gilbert-Varshamov}).
For finite $n$, it is actually (slightly) stronger!  That is,
the constraint that the code be Type II makes it easier to find good codes.

Similar arguments prove:

\begin{theorem}\label{GV4}
In each family from the list \EoneI, \EoneII, \Etwo, \Ethree, \Efour,
\EfiveI, \EfiveII, \Esix\ and \Eseven\ there exists a sequence of
self-dual codes with length tending to infinity satisfying
$$
\liminf_{i\to\infty} {d_i\over n_i} \ge \delta,
$$
where
$$
H_q(\delta)=\delta\log_q(q-1)-\delta\log_q(\delta)-(1-\delta)\log_q(1-\delta)
= {1\over 2}.
$$
\end{theorem}

The result for families \Esix\ and \Eseven\ was first given by Pless and Pierce \cite{Ple16}.

Similar results hold for self-dual codes over ${\Bbb Z}_4$:

\begin{theorem}\label{GV5}
There exists a family of Type II self-dual codes over ${\Bbb Z}_4$, with
length tending to infinity, such that
$$
\liminf_{i\to\infty} {l_i\over 2 n_i} \ge \delta,
$$
where $l_i$ is the minimal Lee weight of the $i$th code and
$\delta = H_2^{-1} (1/2)$, as before.
\end{theorem}

\begin{theorem}\label{GV6}
There exists a family of Type II self-dual codes over ${\Bbb Z}_4$,
with length tending to infinity, such that
$$
\liminf_{i\to\infty} {N_i\over n_i} \ge .34737283\ldots,
$$
where $N_i$ is the minimal Euclidean norm of the $i$th code.
\end{theorem}
\section{Enumeration of self-dual codes}\label{GU}\index{enumeration!self-dual codes}
\index{self-dual codes!enumeration of}
\subsection{Gluing theory}\label{GU1}\index{gluing theory}
\hsp
Gluing is a technique for building up self-dual codes from smaller codes,
and is especially useful when one is attempting to classify
all self-dual codes of a given length.
Typically one finds that there are many codes with low minimal distance and only a few with high minimal distance.
Gluing theory is good at finding all the codes of low distance.

The first formal description of gluing theory appeared in \cite{cp}.
It has also been used in
\cite{Me65}, \cite{SPLAG}, \cite{Me168}, \cite{Me76}, \cite{Me89}, etc.

\begin{figure}[htb]
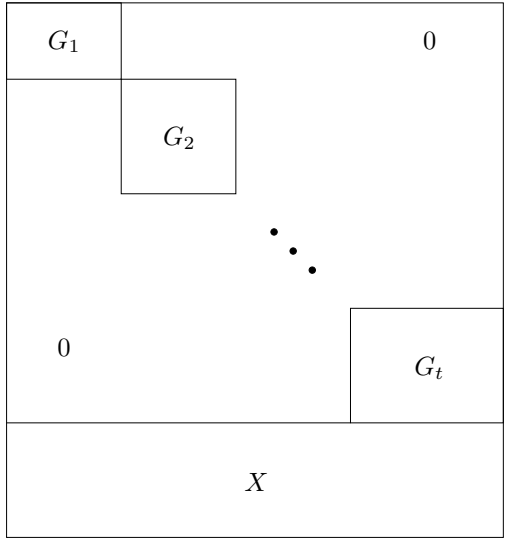

\begin{center}
\input box.pstex_t
\end{center}

\caption{Generator matrix $G$ for a code formed by gluing components $C_1, \ldots, C_t$ together. $G_i$ is a generator matrix for $C_i$, and $X$ denotes the rest of the generator matrix for $C$.}
\label{FigFG1}
\end{figure}

The theory applies to codes from any of the families that we have discussed in this chapter.
Let $C_1 , \ldots , C_t$ be self-orthogonal codes of lengths $n_1, \ldots, n_t$ with generator matrices $G_1, \ldots, G_t$.
If $C$ is a self-dual code with the generator matrix shown in Fig.~\ref{FigFG1} then we say that $C$ is formed by {\em gluing} the {\em components}
\index{component code} \index{code!component}
\index{$+$, gluing notation}
$C_1, \ldots, C_t$ together, and we write
\beql{Eqglue}
C = (C_1 C_2 \ldots C_t)^+
\eeq
to indicate this process.
(Whenever possible the subcodes are chosen so that every minimal weight
codeword of $C$ belongs to one of the $C_i$.)
The codewords in $C$ which contain a nonzero linear combination of the rows of the matrix $X$
are called {\em glue words},\index{glue word}
since these hold the components together.
A glue word has the form
\beql{EqGU1}
u = u_1 u_2 \ldots u_t ~,
\eeq
where each glue element $u_i$ has length $n_i$.
Since $C$ is self-dual, $u_i$ is in $C_i^\perp$.

Let us choose coset representatives $a_0 = 0$, $a_1, \ldots, a_{s-1}$ for $C_i$ in $C_i^\perp$, where
$s = | C_i^\perp | / | C_i |$, so that
$$
C_i^\perp = \bigcup_{j=0}^{s-1} (a_j + C_i ) ~.
$$
Then we can assume that each $u_i$ in \eqn{EqGU1} is one of $a_0, \ldots, a_{s-1}$.

As illustrations we give the two indecomposable
binary Type I
self-dual codes of length 18 (see Tables~\ref{T2a} and \ref{T2d}),
using the components from the list in Section~\ref{GU3}.
The first code is formed by gluing three copies of
the component $d_6$ together:
\beql{EqGU2}
\begin{array}{|l@{}l@{}l@{}l@{}l@{}l|l@{}l@{}l@{}l@{}l@{}l|l@{}l@{}l@{}l@{}l@{}l|} \hline
1&1&1&1&~&~&~&~&~&~&~&~&~&~&~&~&~&~ \\
~&~&1&1&1&1&~&~&~&~&~&~&~&~&~&~&~&~ \\ \hline
~&~&~&~&~&~&1&1&1&1&~&~&~&~&~&~&~&~ \\
~&~&~&~&~&~&~&~&1&1&1&1&~&~&~&~&~&~ \\ \hline
~&~&~&~&~&~&~&~&~&~&~&~&1&1&1&1&~&~ \\
~&~&~&~&~&~&~&~&~&~&~&~&~&~&1&1&1&1 \\ \hline
0&1&0&1&0&1&0&0&0&0&1&1&0&1&0&1&1&0 \\
0&1&0&1&1&0&0&1&0&1&0&1&0&0&0&0&1&1 \\
0&0&0&0&1&1&0&0&0&0&1&1&0&0&0&0&1&1 \\ \hline
\end{array}
\eeq
The three glue vectors shown are $abc$, $cab$ and $bbb$.

The second code is formed by gluing together $d_{10}$, $e_7$ and a ``free'' 
(or empty) component $f_1$:
\beql{EqGU3}
\begin{array}{|l@{}l@{}l@{}l@{}l@{}l@{}l@{}l@{}l@{}l|l@{}l@{}l@{}l@{}l@{}l@{}l|l|} \hline
1&1&1&1&~&~&~&~&~&~&~&~&~&~&~&~&~&~ \\
~&~&1&1&1&1&~&~&~&~&~&~&~&~&~&~&~&~ \\
~&~&~&~&1&1&1&1&~&~&~&~&~&~&~&~&~&~ \\
~&~&~&~&~&~&1&1&1&1&~&~&~&~&~&~&~&~ \\ \hline
~&~&~&~&~&~&~&~&~&~&1&1&1&0&1&0&0&~ \\
~&~&~&~&~&~&~&~&~&~&0&1&1&1&0&1&0&~ \\
~&~&~&~&~&~&~&~&~&~&0&0&1&1&1&0&1&~ \\ \hline
0&1&0&1&0&1&0&1&0&1&0&0&0&0&0&0&0&1 \\
0&1&0&1&0&1&0&1&1&0&1&1&1&1&1&1&1&0 \\ \hline
\end{array}
\eeq
The two glue vectors shown are $a0A$ and $cd0$.

Of course a self-{\em dual} code has no (nonzero) glue.
If
a self-orthogonal code $C$ has a component $B$, say, which is self-dual, then $C$ is a direct sum
$C=B \oplus C'$, where $C'$ is again self-orthogonal.

It may happen that there is a glue word in which only one $u_i$ is nonzero,
in which case we say that the component $C_i$ has {\em self-glue},\index{self-glue} and that $u$ is a {\em self-glue vector}.
So if $C$ has a single component $C_1$ (say) with self-glue, we write $C = C_1^+$ (compare \eqn{Eqglue}).

A basic result of gluing theory is the following.
\begin{theorem}\label{thGU1}
If a self-dual code $C$ is formed by gluing together two codes $C_1$ and $C_2$
in such a way that there is no self-glue, then the quotient groups $C_1^\perp / C_1$ and $C_2^\perp / C_2$ are
isomorphic.
\end{theorem}

We omit the easy proof.
The isomorphism is given by $u_1 + C_1 \to u_2 + C_2$ whenever there is a glue vector $u_1 u_2$.
\subsection{Automorphism groups of glued codes}\label{GU2}
\hsp
One advantage of the gluing method is that it makes it much easier to find the automorphism group of a self-dual code $C$.
We will denote the group by $G(C)$ rather than $Aut (C)$ in this section.
It is essential that every automorphism of $C$ takes the set of component codes $C_1, \ldots, C_t$ to itself.
We will always choose the components so that this is true.

This being the case, any automorphism in $G(C)$ will effect some permutation of the $C_i$, so that $G(C)$ will have a normal subgroup $G_{01}$ consisting of just those elements for which this permutation is trivial.
The group of permutations of the components that are realized in this way we call
$G_2 (C)$ --- it is isomorphic to the quotient group $G(C) / G_{01}$.

Let $G_0 (C)$ be the normal subgroup of $G_{01}$ consisting of those automorphisms which, for every $i$, send each glue element $u_i$ into a vector in the same
coset $u_i + C_i$, i.e. which fix the glue elements modulo the components.
Then $G_{01} / G_0 (C)$ is isomorphic to a group acting on the glue elements
of each component:
we call this group $G_1 (C)$.
Thus the full group $G(C)$ is compounded of the groups $G_0 (C)$, $G_1(C)$ and $G_2(C)$, and has order
\beql{EqGU4}
|G(C)| = |G_0(C) | |G_1(C) || G_2 (C) | ~.
\eeq
Also $G_0 (C) $ is the direct product of the groups $G_0 (C_i)$.
But in general $G_1 (C)$ is only a subgroup of the direct product of
the $G_1 (C_i)$, and therefore must be computed directly for each $C$.
\subsection{Family \Eone: Enumeration of binary self-dual codes}\label{GU3}
\hsp
The enumeration of binary self-dual codes of length $n \le 32$ has been carried out in a series of papers:
Pless \cite{Ple10} for $n \le 20$;
Conway (unpublished) for Type II of length 24;
Pless and Sloane \cite{Me39} for $n=22$, 24;
Conway and Pless \cite{cp} for $n=26$ to 30 and Type II of length 32 (see also
Pless \cite{Ple12}).
Some errors in the last two references were corrected in Conway, Pless and Sloane \cite{Me162}.
The results are summarized in Table~\ref{T2z}.
\begin{table}
\caption{Numbers of self-dual codes of length $n$. (a)~Indecomposable Type II, (b)~total Type II, (c)~indecomposable self-dual, (d)~total self-dual.}
$$
\begin{array}{c|ccccccccc}
n & 0 & 2 & 4 & 6 & 8 & 10 & 12 & 14 & 16 \\ \hline
a & 1 & - & - & - & 1 & - & - & - & 1 \\
b & 1 & - & - & - & 1 & - & - & - & 2 \\
c & 1 & 1 & 0 & 0 & 1 & 0 & 1 & 1 & 2 \\
d & 1 & 1 & 1 & 1 & 2 & 2 & 3 & 4 & 7
\end{array}
$$
$$
\begin{array}{c|cccccccc}
n & 18 & 20 & 22 & 24 & 26 & 28 & 30 & 32 \\ \hline
a & - & - & - & 7 & - & - & - & 74 \\
b & - & - & - & 9 & - & - & - & 85 \\
c & 2 & 6 & 8 & 26 & 45 & 148 & 457 \\
d & 9 & 16 & 25 & 55 & 103 & 261 & 731
\end{array}
$$
\label{T2z}
\end{table}

In this section we describe these codes, drawing heavily from the tables in \cite{Me162}.

Since (from \eqn{EqAGa}) there are at least 17493 inequivalent Type II codes of length 40, length 32 is probably a good place to stop.

Although the Type I codes of length 32 have not been classified,
it is shown in \cite{Me158} that there are precisely three inequivalent $[32,16,8]$ extremal Type I codes.

The following self-orthogonal codes will be used as components.

$d_4: [1111]$, glue: $a=0011$, $b=0101$, $c=0110$, $|G_0| =4$, $G_1 = S(3)$ on $\{a,b,c\}$, $|G_1| = 6$.

$d_{2n} ( n \ge 3)$:
\beql{EqGU8}
[111100 \ldots, 00111100 \ldots, \ldots, \ldots 001111 ]~,
\eeq
glue: $a=0101 \ldots 01$, $b=0000 \ldots 11$,
$c=0101 \ldots 10$, $|G_0| = 2^{n-1} n!$,
$|G_1| =2$ (swap $a$ and $c$)

$e_7: [(1110100)]$, glue: $a=1111111$,
$G_0 = L_3 (2)$, $|G_0| = 168$, $|G_1| =1$.

$e_8$ is the $[8,4,4]$ Hamming code, see Section~\ref{Wee2}.

$f_n :$ If some coordinate positions contain very few codewords, it is often best to regard these places as containing the {\em free} (or {\em empty}) {\em component} $f_n = \{0^n \}$.
In this case we label the coordinate positions by $A$, $B$, $C, \ldots$, and use $ABD$ for example to denote the glue word
$110100 \ldots~$.
Also $|G_0 |=1$.

The above components are important in view of the following decomposition theorem for binary codes with low minimal distance.
\begin{theorem}\label{thGU2}
(a)~If a self-orthogonal code $C$ has minimal distance 2 then $C= i_2^k \oplus C'$, where $C'$ has minimal distance at least 4.
(b)~If a self-orthogonal code $C$ is generated by words of weight 4 then $C$ is a direct sum of copies of the codes $d_{2m}$ $(m \ge 2)$, $e_7$ and $e_8$.
\end{theorem}
\paragraph{Proof.}
(a)~Suppose $C$ contains a word of weight 2, say $u= 1100 \ldots~$.
Then any other word $v \in C$ must meet $u$ evenly, so begins $00 \ldots$ or $11 \ldots$~.
Hence $C= B \oplus C'$ where $B= [11]$.
(b)~A set of mutually self-orthogonal words of weight 4 whose supports are linked is easily seen to be either a $d_{2m}$ for some $m \ge 2$, or an $e_7$ or $e_8$.~~~$\bsq$
\paragraph{Remarks.}
(1)~Suppose $C$ is a self-dual code with minimal distance 4, and let $C'$ be the subcode generated by words of weight 4.
Then $C'$ is as described in part (b) of the theorem, and $C$ can be regarded as being obtained by gluing $C'$ to some other subcode $C''$ (the latter may be the free component $f_n$).

(2)~Generalizing Part (a) of the theorem, it is easy to show that any self-orthogonal code over a field $\FF_q$ with length $n > 2$ and minimal distance 2 is decomposable (\cite{Me65}, Theorem~3).

The following are some additional components that will be used in Table~\ref{T2a}.

The code $g_{24-m}$
$( m = 0, 2, 3, 4, 6, 8 )$ is obtained by taking the words of the Golay code $g_{24}$ that vanish on $m$ digits (and
then deleting those digits).
For the [16,5,8] first-order Reed-Muller code\index{Reed-Muller code} $g_{16}$ 
\index{code!Reed-Muller}
the 8 digits must be a special
octad, while for $g_{18}$ they must be an umbral hexad
(see \cite{SPLAG} for terminology).
For $0 \le m \le 6$, $g_{24-m}$ is a $[24-m , 12-m , 8]$ code.

The [24,11,8] {\em half-Golay code}\index{half-Golay code} $h_{24}$ consists of the Golay codewords that intersect a given tetrad evenly.
\index{code!half-Golay}\index{Golay code!half-}

The {\em odd Golay code}\index{odd Golay code}
\index{Golay code!odd} \index{code!odd Golay}
\index{Golay code}\index{code!Golay, q.v.}
$h_{24}^+$ is the
$[24,12,6]$ Type I code generated by $h_{24}$ and an appropriate vector of weight 6.
Alternatively, the odd Golay code may be obtained as follows.
Let $v \in \FF_{2}^{\,24}$ be a fixed vector of weight 4,
say $v= 1^4 0^{24}$.
Then $h_{24}^+ = \{ u \in g_{24} : wt (u \cap v)~{\rm even} \} \bigcup \{ u+v : u \in g_{24}, wt (u \cap v) ~{\rm odd} \}$, with generator matrix
$$
\begin{array}{|cccc|cccc|cccc|cccc|cccc|cccc|} \hline
1&1&1&1&1&1&1&1&0&0&0&0&0&0&0&0&0&0&0&0&0&0&0&0 \\ 
1&1&1&1&0&0&0&0&1&1&1&1&0&0&0&0&0&0&0&0&0&0&0&0 \\ 
1&1&1&1&0&0&0&0&0&0&0&0&1&1&1&1&0&0&0&0&0&0&0&0 \\ 
1&1&1&1&0&0&0&0&0&0&0&0&0&0&0&0&1&1&1&1&0&0&0&0 \\ 
1&1&1&1&0&0&0&0&0&0&0&0&0&0&0&0&0&0&0&0&1&1&1&1 \\ \hline
1&1&0&0&0&0&0&0&0&0&0&0&1&1&0&0&1&0&1&0&1&0&1&0 \\ 
0&0&0&0&1&1&0&0&0&0&0&0&1&0&1&0&1&1&0&0&1&0&1&0 \\ 
0&0&0&0&0&0&0&0&1&1&0&0&1&0&1&0&1&0&1&0&1&1&0&0 \\ \hline
1&0&1&0&0&0&0&0&0&0&0&0&1&0&1&0&1&0&0&1&1&0&0&1 \\ 
0&0&0&0&1&0&1&0&0&0&0&0&1&0&0&1&1&0&1&0&1&0&0&1 \\ 
0&0&0&0&0&0&0&0&1&0&1&0&1&0&0&1&1&0&0&1&1&0&1&0 \\ \hline
1&0&0&0&1&0&0&0&1&0&0&0&1&0&0&0&1&0&0&0&1&0&0&0 \\ \hline
\end{array}
$$
This code has weight enumerator
$x^{24} + 64 x^{18} y^6 + 375 x^{16} y^8 + 960 x^{14} y^{10} + 1296 x^{12} y^{12}
+960 x^{10} y^{14} + 375 x^8 y^{16} + 64 x^8 y^{18} + y^{24}$,
and $Aut (h_{24}^+)$ is the ``sextet group'' $2^6\mathord{:}3.S(6)$, of order
$2^{10}.3^3 .5 = 138240$
(\cite{Me39}, \cite[p.~309]{SPLAG}).

The first 11 rows of the above matrix
generate $h_{24}$; if the last row is replaced by
$$
\begin{array}{|cccc|cccc|cccc|cccc|cccc|cccc|} \hline
0&1&1&1&1&0&0&0&1&0&0&0&1&0&0&0&1&0&0&0&1&0&0&0 \\ \hline
\end{array}
$$
we get the Golay code $g_{24}$ itself; and if 
the last row is replaced by
$$
\begin{array}{|cccc|cccc|cccc|cccc|cccc|cccc|} \hline
1&1&1&1&0&0&0&0&0&0&0&0&0&0&0&0&0&0&0&0&0&0&0&0 \\ \hline
\end{array}
$$
we get $(d^6_4)^+$.

Under the action of $Aut ( g_{24} )$ there are two distinct ways to select tetrads
$t = \{c,d,e,f \}$, $u = \{ a,b,e,f \}$, $v = \{a,b,c,d \}$ so that $t+ u+v = 0$,
depending on whether $\{ a,b, \ldots ,f \}$ is a special hexad or an umbral hexad (see Fig.~\ref{FG2}).
Correspondingly there are two [24,10,8] {\em quarter-Golay codes}\index{quarter-Golay code}
\index{Golay code!quarter}\index{code!quarter-Golay}
$q_{24}^+$, $q_{24}^-$,
consisting of the codewords of $g_{24}$ that intersect all of $t,u,v$ evenly.

We refer to \cite{cp} and \cite{Me162} for a description of the glue vectors for these codes.
\begin{figure}[htb]
\caption{Two choices for a hexad (special or umbral), used to define the two $[24,10,8]$ quarter-Golay codes $q_{24}^+$ and $q_{24}^-$.}
$$
\begin{array}{|cc|cc|cc|} \hline
~&~&~&~&~&~ \\
~&~&~&~&a&b \\ \hline
~&~&~&~&c&d \\
~&~&~&~&e&f \\ \hline
\end{array} \quad 
\begin{array}{|cc|cc|cc|} \hline
a&b&c&d&e&f \\
~&~&~&~&~&~ \\ \hline
~&~&~&~&~&~ \\
~&~&~&~&~&~ \\ \hline
\end{array}
$$
\label{FG2}
\end{figure}

Our first table (Table~\ref{T2a}) lists all indecomposable binary self-dual codes of length $n \le 22$,
together with the indecomposable Type II codes of length 24,
using the $+$ notation of \eqn{Eqglue}.
For these codes (and for most of those in the following tables) there is only one way to glue the specified components together without introducing
additional minimal-weight words.
We have therefore omitted the glue words from the table.
(However, more information about these codes, including the glue words, will be given in Table~\ref{T2d}.)
\begin{table}[htb]
\caption{Indecomposable binary self-dual codes of length $n \le 24$ (\S ~indicates a Type II code. For length 24 only Type II codes are listed).}
$$
\begin{array}{rl} \hline
\multicolumn{1}{c}{\mbox{Length $n$}} & \mbox{Components} \\ \hline
2 & i_2 \\
8 & e_8^{\S} \\
12 & d_{12}^+ \\
14 & e_7^{2+} \\
16 & d_{16}^{\S} , d_8^{2+} \\
18 & (d_{10} e_7 f_1)^+, d_6^{3+} \\
20 & d_{20}^+ , (d_{12} d_8)^+ , (d_8^2 d_4)^+ , (e_7^2 d_6)^+ , (d_6^3 f_2)^+, d_4^{5+} \\
22 & g_{22}^+ , (d_{14} e_7 f_1)^+, (d_{10}^2 f_2)^+ , (d_{10} d_6^2)^+ , (d_8 e_7 d_6 f_1)^+, \\
& (d_8 d_6^2 f_2)^+ , (d_6^2 d_4^2 f_2)^+, (d_4^4 f_3^2)^+ \\
24^{\S} & g_{24}, d_{24}^+ , d_{12}^{2+} , (d_{10} e_7^2)^+, d_8^{3+} , d_6^{4+} , d_4^{6+}
\end{array}
$$
\label{T2a}
\end{table}

The next table (Tables~\ref{T2b} and \ref{T2b2}) gives the full list of all 85 (decomposable or indecomposable) Type II codes of length 32.
This table is taken from \cite{Me162}, and is a corrected version of the table in \cite{cp}.
The codes are labeled from C1 to C85 in the first column (using the same order as in \cite{cp} and \cite{Me162}).
The second column gives the components (omitting the superscripts ``$+$'' to save space).

The third and fourth columns give the orders of the groups $G_1 (C)$ and $G_2 (C)$, and the fifth column gives the order of the full group,
using \eqn{EqGU4},
where $|G_0 (C)|$ is the product of the orders of the $G_0 (C_i)$ for the components.
The latter are given in Table~\ref{T2c}.
The next column gives $A_4$, the number of codewords of weight 4.
The weight enumerator of the code is then (from Theorem~\ref{ThMS3c})
$$(x^8 + 14x^4 y^4 + y^8 )^4 -
(56 - A_4) x^4 y^4 (x^4 - y^4)^4 (x^8 +14 x^4 y^4 + y^8) ~.$$
The last four columns give the number of self-dual codes (the ``children'', cf. Chapter~xx (Pless)) of lengths 30, 28, 26, 24 that arise from the code.

To save space, we have omitted the glue vectors from Tables~\ref{T2b} and \ref{T2b2}.
In many cases they are uniquely determined by the components, and in any case
they can be found in full in \cite{cp}, with corrections in \cite{Me162}.

\begin{table}[htb]
\caption{Doubly-even self-dual (or Type II) binary codes of length 32 (Part 1)}

\begin{center}
\begin{tabular}{llrrlrrrrr} \hline
Code & Components & \multicolumn{1}{c}{$|G_1|$} & \multicolumn{1}{c}{$|G_2|$} & $|G|$ & \multicolumn{1}{c}{$A_4$} & \multicolumn{1}{c}{$n_{30}$} &
\multicolumn{1}{c}{$n_{28}$} & \multicolumn{1}{c}{$n_{26}$} &
\multicolumn{1}{c}{$n_{24}$} \\ \hline
C1 & $d_{32}$ & 1 & 1 & $2^{30} 3^{6} 5^{3} 7^{2} 11.13$ & 120 & 2 & 2 & 1 & 1 \\
C2 & $d_{24} e_{8}$ & 1 & 1 & $2^{27} 3^{6} 5^{2} 7^{2} 11$ & 80 & 4 & 3 & 2 & 2 \\
C3 & $d_{20} d_{12}$ & 1 & 1 & $2^{26} 3^{6} 5^{3} 7$ & 60 & 5 & 4 & 2 & 2 \\
C4 & $d_{18} e_{7}^{2}$ & 1 & 2 & $2^{22} 3^{6} 5.7^{3}$ & 50 & 5 & 3 & 2 & 1 \\
C5 & $d_{16}^{2}$ & 1 & 2 & $2^{29} 3^{4} 5^{2} 7^{2}$ & 56 & 3 & 2 & 1 & 1 \\
C6 & $d_{16} e_{8}^{2}$ & 1 & 2 & $2^{27} 3^{4} 5.7^{3}$ & 56 & 5 & 3 & 2 & 2 \\
C7 & $d_{16} d_{8}^{2}$ & 1 & 2 & $2^{27} 3^{4} 5.7$ & 40 & 6 & 4 & 2 & 2 \\
C8 & $d_{14} d_{10} e_{7} f_{1}$ & 1 & 1 & $2^{20} 3^{4} 5^{2} 7^{2}$ & 38 & 11 & 5 & 3 & 2 \\
C9 & $d_{14} d_{6}^{3}$ & 1 & 6 & $2^{20} 3^{6} 5.7$ & 30 & 6 & 4 & 2 & 1 \\
C$10$ & $d_{12}^{2} e_{8}$ & 1 & 2 & $2^{25} 3^{5} 5^{2} 7$ & 44 & 5 & 3 & 2 & 2 \\
C11 & $d_{12}^{2} d_{8}$ & 1 & 2 & $2^{25} 3^{5} 5^{2}$ & 36 & 6 & 4 & 2 & 2 \\
C12 & $d_{12} d_{8}^{2} d_{4}$ & 1 & 2 & $2^{24} 3^{4} 5$ & 28 & 11 & 7 & 2 & 2 \\
C13 & $d_{12} e_{7}^{2} d_{6}$ & 1 & 2 & $2^{19} 3^{5} 5.7^{2}$ & 32 & 9 & 5 & 3 & 1 \\
C14 & $d_{12} d_{6}^{3} f_{2}$ & 1 & 6 & $2^{19} 3^{6} 5$ & 24 & 9 & 4 & 2 & 1 \\
C15 & $d_{12} d_{4}^{5}$ & 1 & 120 & $2^{22} 3^{3} 5^{2}$ & 20 & 5 & 3 & 1 & 1 \\
C16 & $d_{10}^{3} f_{2}$ & 1 & 6 & $2^{22} 3^{4} 5^{3}$ & 30 & 5 & 2 & 1 & 1 \\
C17 & $d_{10}^{2} d_{6}^{2}$ & 1 & 4 & $2^{22} 3^{4} 5^{2}$ & 26 & 7 & 4 & 2 & 1 \\
C18 & $d_{10} e_{8} e_{7}^{2}$ & 1 & 2 & $2^{20} 3^{4} 5.7^{3}$ & 38 & 8 & 4 & 3 & 2 \\
C19 & $d_{10} d_{8} e_{7} d_{6} f_{1}$ & 1 & 1 & $2^{19} 3^{4} 5.7$ & 26 & 17 & 7 & 4 & 2 \\
C20 & $d_{10} d_{8} d_{6}^{2} f_{2}$ & 1 & 2 & $2^{20} 3^{4} 5$ & 22 & 15 & 6 & 3 & 2 \\
C21 & $d_{10} d_{6}^{2} d_{4}^{2} f_{2}$ & 1 & 4 & $2^{19} 3^{3} 5$ & 18 & 16 & 6 & 2 & 1 \\
C22 & $d_{10} d_{4}^{4} f_{6}$ & 6 & 24 & $2^{19} 3^{3} 5$ & 14 & 9 & 3 & 1 & 1 \\
C23 & $d_{10} g_{22}$ & 2 & 1 & $2^{15} 3^{3} 5^{2} 7.11$ & 10 & 4 & 2 & 1 & 1 \\
C24 & $e_{8}^{4}$ & 1 & 24 & $2^{27} 3^{5} 7^{4}$ & 56 & 2 & 1 & 1 & 1 \\
C25 & $e_{8} d_{8}^{3}$ & 1 & 6 & $2^{25} 3^{5} 7$ & 32 & 5 & 3 & 2 & 2 \\
C26 & $e_{8} d_{6}^{4}$ & 1 & 24 & $2^{21} 3^{6} 7$ & 26 & 5 & 3 & 2 & 1 \\
C27 & $e_{8} d_{4}^{6}$ & 3 & 720 & $2^{22} 3^{4} 5.7$ & 20 & 4 & 2 & 1 & 1 \\
C28 & $e_{8} g_{24}$ & 1 & 1 & $2^{16} 3^{4} 5.7^{2} 11.23$ & 14 & 3 & 1 & 1 & 1 \\
C29 & $d_{8}^{4}$ & 1 & 24 & $2^{27} 3^{5}$ & 24 & 3 & 2 & 1 & 1 \\
C30 & $d_{8}^{4}$ & 1 & 8 & $2^{27} 3^{4}$ & 24 & 4 & 2 & 1 & 1 \\
C31 & $d_{8}^{3} d_{4}^{2}$ & 1 & 6 & $2^{23} 3^{4}$ & 20 & 6 & 3 & 1 & 1 \\
C32 & $d_{8}^{2} e_{7}^{2} f_{2}$ & 1 & 4 & $2^{20} 3^{4} 7^{2}$ & 26 & 10 & 3 & 2 & 1 \\
C33 & $d_{8}^{2} d_{6}^{2} f_{4}$ & 1 & 4 & $2^{20} 3^{4}$ & 18 & 14 & 4 & 2 & 1 \\
C34 & $d_{8}^{2} d_{4}^{4}$ & 1 & 16 & $2^{24} 3^{2}$ & 16 & 8 & 4 & 1 & 1 \\
C35 & $d_{8} e_{7} d_{6}^{2} d_{4} f_{1}$ & 1 & 2 & $2^{18} 3^{4} 7$ & 20 & 18 & 7 & 3 & 1 \\
C36 & $d_{8} d_{6}^{4}$ & 1 & 8 & $2^{21} 3^{5}$ & 18 & 7 & 4 & 2 & 1 \\
C37 & $d_{8} d_{6}^{2} d_{6} d_{4} f_{2}$ & 1 & 2 & $2^{18} 3^{4}$ & 16 & 22 & 9 & 3 & 1 \\
C38 & $d_{8} d_{6}^{2} d_{4}^{2} f_{4}$ & 1 & 4 & $2^{18} 3^{3}$ & 14 & 20 & 7 & 2 & 1 \\
C39 & $d_{8} d_{6} d_{4}^{3} f_{6}$ & 2 & 6 & $2^{17} 3^{3}$ & 12 & 17 & 6 & 2 & 1 \\
C40 & $d_{8} d_{4}^{6}$ & 1 & 48 & $2^{22} 3^{2}$ & 12 & 7 & 4 & 1 & 1 \\
C41 & $d_{8} d_{4}^{4} f_{8}$ & 2 & 24 & $2^{18} 3^{2}$ & 10 & 12 & 4 & 1 & 1 \\
C42 & $d_{8} d_{4}^{2} g_{16}$ & 36 & 2 & $2^{17} 3^{3}$ & 8 & 9 & 3 & 1 & 1 \\
C43 & $d_{8} h_{24}$ & 1 & 1 & $2^{16} 3^{4} 5$ & 6 & 5 & 2 & 1 & 1 \\
C44 & $e_{7}^{4} d_{4}$ & 1 & 24 & $2^{17} 3^{5} 7^{4}$ & 29 & 4 & 2 & 1 & 0 \\
\end{tabular}
\end{center}
\label{T2b}
\end{table}
\begin{table}[htb]
\caption{Doubly-even self-dual (or Type II) binary codes of length 32 (Part 2)}
\begin{center}
\begin{tabular}{llrrlrrrrr} \hline
Code & Components & \multicolumn{1}{c}{$|G_1|$} & \multicolumn{1}{c}{$|G_2|$} & $|G|$ & \multicolumn{1}{c}{$A_4$} & \multicolumn{1}{c}{$n_{30}$} &
\multicolumn{1}{c}{$n_{28}$} & \multicolumn{1}{c}{$n_{26}$} &
\multicolumn{1}{c}{$n_{24}$} \\ \hline
C45 & $e_{7}^{2} d_{6}^{3}$ & 1 & 6 & $2^{16} 3^{6} 7^{2}$ & 23 & 6 & 3 & 2 & 0 \\
C46 & $e_{7} d_{6}^{3} d_{4} f_{3}$ & 1 & 6 & $2^{15} 3^{5} 7$ & 17 & 13 & 4 & 2 & 0 \\
C47 & $e_{7} d_{6} d_{4}^{4} f_{3}$ & 1 & 24 & $2^{17} 3^{3} 7$ & 14 & 12 & 4 & 2 & 0 \\
C48 & $e_{7} d_{4}^{4} f_{9}$ & 18 & 24 & $2^{15} 3^{4} 7$ & 11 & 7 & 2 & 1 & 0 \\
C49 & $e_{7} d_{4} g_{21}$ & 6 & 1 & $2^{12} 3^{4} 5.7^{2}$ & 8 & 6 & 2 & 1 & 0 \\
C50 & $d_{6}^{5} f_{2}$ & 1 & 10 & $2^{16} 3^{5} 5$ & 15 & 6 & 2 & 1 & 0 \\
C51 & $d_{6}^{4} d_{4}^{2}$ & 1 & 48 & $2^{20} 3^{5}$ & 14 & 6 & 3 & 1 & 0 \\
C52 & $d_{6}^{4} d_{4} f_{2}^{2}$ & 1 & 8 & $2^{17} 3^{4}$ & 13 & 13 & 4 & 1 & 0 \\
C$53$ & $d_{6}^{4} f_{8}$ & 2 & 24 & $2^{16} 3^{5}$ & 12 & 8 & 2 & 1 & 0 \\
C54 & $d_{6}^{3} d_{4}^{3} f_{2}$ & 1 & 6 & $2^{16} 3^{4}$ & 12 & 12 & 5 & 1 & 0 \\
C55 & $d_{6}^{3} d_{4}^{2} f_{6}$ & 1 & 6 & $2^{14} 3^{4}$ & 11 & 15 & 3 & 1 & 0 \\
C56 & $d_{6}^{2} d_{4}^{4} f_{2}^{2}$ & 1 & 8 & $2^{17} 3^{2}$ & 10 & 16 & 4 & 1 & 0 \\
C57 & $d_{6}^{2} d_{4}^{4} f_{4}$ & 2 & 16 & $2^{19} 3^{2}$ & 10 & 13 & 4 & 1 & 0 \\
C58 & $d_{6}^{2} d_{4}^{3} f_{2} f_{6}$ & 1 & 12 & $2^{14} 3^{3}$ & 9 & 18 & 4 & 1 & 0 \\
C59 & $d_{6}^{2} d_{4}^{2} f_{12}$ & 12 & 4 & $2^{14} 3^{3}$ & 8 & 14 & 3 & 1 & 0 \\
C60 & $d_{6}^{2} g_{20}$ & 4 & 2 & $2^{15} 3^{3} 5$ & 6 & 6 & 2 & 1 & 0 \\
C$61$ & $d_{6} d_{4}^{2} d_{4}^{3} f_{3}^{2}$ & 1 & 12 & $2^{15} 3^{2}$ & 8 & 19 & 6 & 1 & 0 \\
C62 & $d_{6} d_{4}^{4} f_{10}$ & 2 & 8 & $2^{15} 3$ & 7 & 21 & 4 & 1 & 0 \\
C63 & $d_{6} d_{4}^{3} f_{14}$ & 8 & 6 & $2^{13} 3^{2}$ & 6 & 18 & 4 & 1 & 0 \\
C64 & $d_{6} d_{4}^{2} g_{18}$ & 36 & 2 & $2^{10} 3^{4}$ & 5 & 12 & 3 & 1 & 0 \\
C$65$ & $d_{6} d_{4} g_{16} f_{6}$ & 72 & 1 & $2^{12} 3^{3}$ & 4 & 14 & 4 & 1 & 0 \\
C$66$ & $d_{6} f_{13}^{2}$ & 5616 & 2 & $2^{8} 3^{4} 13$ & 3 & 6 & 2 & 1 & 0 \\
C67 & $d_{4}^{8}$ & 6 & 1344 & $2^{23} 3^{2} 7$ & 8 & 2 & 1 & 0 & 0 \\
C68 & $d_{4}^{8}$ & 1 & 1152 & $2^{23} 3^{2}$ & 8 & 3 & 1 & 0 & 0 \\
C69 & $d_{4}^{8}$ & 1 & 336 & $2^{20} 3.7$ & 8 & 2 & 1 & 0 & 0 \\
C$70$ & $d_{4}^{6} f_{8}$ & 4 & 48 & $2^{18} 3$ & 6 & 7 & 2 & 0 & 0 \\
C$71$ & $d_{4}^{6} f_{8}$ & 1 & 48 & $2^{16} 3$ & 6 & 9 & 2 & 0 & 0 \\
C72 & $d_{4}^{4} d_{4} f_{12}$ & 6 & 24 & $2^{14} 3^{2}$ & 5 & 10 & 2 & 0 & 0 \\
C73 & $d_{4}^{5} f_{12}$ & 1 & 60 & $2^{12} 3.5$ & 5 & 6 & 1 & 0 & 0 \\
C74 & $d_{4}^{4} g_{16}$ & 8 & 24 & $2^{18} 3$ & 4 & 7 & 2 & 0 & 0 \\
C$75$ & $d_{4}^{4} f_{16}$ & 8 & 8 & $2^{14}$ & 4 & 14 & 2 & 0 & 0 \\
C76 & $d_{4}^{3} g_{18} f_{2}$ & 8 & 6 & $2^{10} 3^{2}$ & 3 & 13 & 2 & 0 & 0 \\
C$77$ & $d_{4}^{2} q_{24}^+$ & 6 & 2 & $2^{15} 3^{2}$ & 2 & 6 & 1 & 0 & 0 \\
C$78$ & $d_{4}^{2} q_{24}^-$ & 3 & 2 & $2^{10} 3^{2}$ & 2 & 8 & 1 & 0 & 0 \\
C$79$ & $d_{4} f_{4}^{6}$ & 16 & 72 & $2^{11} 3^{2}$ & 2 & 8 & 1 & 0 & 0 \\
C$80$ & $d_{4} f_{7}^{4}$ & 168 & 8 & $2^{8} 3.7$ & 1 & 8 & 2 & 0 & 0 \\
C81 & $q_{32}$ & 1 & 1 & $2^{5} 3.5.31$ & 0 & 1 & 0 & 0 & 0 \\
C82 & $r_{32}$ & 1 & 1 & $2^{15} 3^{2} 5.7.31$ & 0 & 1 & 0 & 0 & 0 \\
C83 & $g_{16}^{2}$ & 20160 & 2 & $2^{15} 3^{2} 5.7$ & 0 & 2 & 0 & 0 & 0 \\
C84 & $f_{4}^{8}$ & 256 & 336 & $2^{12} 3.7$ & 0 & 2 & 0 & 0 & 0 \\
C85 & $f_{2}^{16}$ & 2 & 11520 & $2^{9} 3^{2} 5$ & 0 & 3 & 0 & 0 & 0 \\
\end{tabular}
\end{center}
\label{T2b2}
\end{table}
\begin{table}[htb]
\caption{The groups $G_0$ for the components mentioned in Tables~\ref{T2a}, \ref{T2b} and \ref{T2b2}.}
\begin{center}
\begin{tabular}{ccc} \hline
Component & $G_0$ & $|G_0 |$ \\ \hline
$d_{2m}$ & $2^{m-1}  .  S(m)$ & $2^{m-1} m!$ \\
$e_{7}$ & $L_3 (2)$ & 168 \\
$e_{8}$ & $GA_3 (2)$ & 1344 \\
$f_{n}$ & 1 & 1 \\
$g_{16}$ & $2^{4}$ & 16 \\
$g_{18}$ & $Z(3)$ & 3 \\
$g_{20}$ & $M_{20}$ & $2^6 3 .5$ \\
$g_{21}$ & $M_{21}$ & $2^6 3^2 5.7$ \\
$g_{22}$ & $M_{22}$ & $2^7 3^2 5.7.11$ \\
$g_{24}$ & $M_{24}$ & $2^{10} 3^3 5.7.11.23$ \\
$h_{24}$ & $2^6 : 3 S(6)$ & $2^9 3^3 5$ \\
$q_{24}^+$ & $2^6  .  (S(3)  \times 2^2 )$ & $2^9 3$ \\
$q_{24}^-$ & $2^2 \times S(4)$ & $2^5 3$ \\
\end{tabular}
\end{center}
\label{T2c}
\end{table}

The enumeration in Tables~\ref{T2b} and \ref{T2b2} has been subjected to many checks, including the verification of the mass formula
$$
\sum \frac{1}{|Aut (C)|} =
\frac{391266122896364123}{532283035423762022400}
$$
(in agreement with \eqn{EqAGa}).
\paragraph{Remark.}\label{RemGU1}
There are just five Type II codes of length 32 with minimal distance 8:
the quadratic residue code $C81 = q_{32}$, generated by
$$(1001000110110111100010101110000)1 ~;$$
the second-order Reed-Muller code $C82 = r_{32}$, generated by
$$(1110010000010000001100000000000)1 ~;$$
and the three codes $C83 = g_{16}^{2+}$, $C84 = f_4^{8+}$ and $C85 = f_2^{16+}$.
Explicit generator matrices for the last three are shown in Fig.~\ref{FG3}.
\begin{figure}[htb]
\caption{Generator matrices for the $[32,16,8]$ Type II codes $C83 = g_{16}^{2+}$, $C84 = f_4^{8+}$ and $C85 = f_2^{16+}$.}

{\footnotesize
$$
\begin{array}{l}
11101000111010000000000000000000 \\
10110100101101000000000000000000 \\
10011010100110100000000000000000 \\
10001101100011010000000000000000 \\
00000000000000001101100011011000 \\
00000000000000001010110010101100 \\
00000000000000001001011010010110 \\
00000000000000001000101110001011 \\
11011000110110001101100000000000 \\
10101100101011001010110000000000 \\
10010110100101101001011000000000 \\
10001011100010111000101100000000 \\
00000000111010001110100011101000 \\
00000000101101001011010010110100 \\
00000000100110101001101010011010 \\
00000000100011011000110110001101 \\
\end{array}
$$
$$
\begin{array}{l}
11101000000000001110100011101000 \\
10110100000000001011010010110100 \\
10011010000000001001101010011010 \\
10001101000000001000110110001101 \\
00000000111010001110100010110100 \\
00000000101101001011010010011010 \\
00000000100110101001101010001101 \\
00000000100011011000110111000110 \\
11011000110110001101100000000000 \\
10101100101011001010110000000000 \\
10010110100101101001011000000000 \\
10001011100010111000101100000000 \\
11011000101100010000000011011000 \\
10101100110110000000000010101100 \\
10010110101011000000000010010110 \\
10001011100101100000000010001011 \\
\end{array}
$$
$$
\begin{array}{l}
10000000000000001111100010001000 \\
01000000000000001111010001000100 \\
00100000000000001111001000100010 \\
00010000000000001111000100010001 \\
00001000000000001000111110001000 \\
00000100000000000100111101000100 \\
00000010000000000010111100100010 \\
00000001000000000001111100010001 \\
00000000100000001000100011111000 \\
00000000010000000100010011110100 \\
00000000001000000010001011110010 \\
00000000000100000001000111110001 \\
00000000000010001000100010001111 \\
00000000000001000100010001001111 \\
00000000000000100010001000101111 \\
00000000000000010001000100011111 \\
\end{array}
$$
}
\label{FG3}
\end{figure}
\paragraph{Subtraction.}\index{subtraction}
Suppose for concreteness that $C$ is a Type I code of length 26 with doubly even subcode $C_0$.
Then we obtain a Type II code $B$ (say) of length 32 by gluing $C_0$ to $d_6$,
as follows.
Write $C_0^\perp = C_0 \cup C_1 \cup C_2 \cup C_3$,
as in Section~\ref{Shad}, where $C= C_0 \cup C_2$, the shadow of $C$ is
$C_1 \cup C_3$, and $C_i = u_i +C_0$ for $i=1,2,3$.
Then $B$ is generated by
\beql{EqGU6}
\begin{array}{|c|c|} \hline
C_0 & ~~ \\ \hline
~~ & d_6 \\ \hline
u_1 & a \\ \hline
u_2 & b \\ \hline
u_3 & c \\ \hline
\end{array}
\eeq

This is a special case of the following construction.
Let $C$, $D$ be any strictly Type I codes, of lengths $n_1$ and $n_2$, respectively, with
$C_0^\perp = \cup_{i=0}^3 C_i$, $D_0^\perp = \cup_{i=0}^3 D_i$.
Then $B= \cup_{i=0}^3 C_i \times D_i$ is self-dual if $n_1 + n_2 \equiv 0$ $ ( \bmod~4)$, and is Type II if $n_1 + n_2 \equiv 0$ $( \bmod ~8)$.
The weight enumerator of $B$ is then
$$\sum_{i=0}^3 W_{C_i} (x,y) W_{D_i} (x,y) ~.$$
Several constructions in the literature (\cite{BrPl91}, Theorems~1 and 2;
\cite{DGH97a}, Theorem~3.1, for example) are special cases of this construction.
In \eqn{EqGU6} we have $D= i_2^3$.

In this way any Type I code of length 26 leads to a unique (up to equivalence)
Type II code of length 32.

Conversely, all Type I codes of length 26 can be obtained by choosing a $d_6$
inside a Type II code of length 32 and inverting the above process.

More generally, suppose $B$ is a Type II code of length $n$.
We choose a copy of $D = i_2^m$ so that $D_0 = d_{2m} \subset B$.
Then we obtain a Type I code of length $n-2m$ by taking the vectors $v$ such that $vw \in B$ for some $w \in D$.
We call this process {\em subtraction}.\index{subtraction}
Every Type I code of length $n-2m$ can be obtained in this way by starting with a
unique Type II code and subtracting an appropriate $d_{2m}$.
Of course any Type II code of length $n-2m$ is a direct summand of some Type II code of any greater length.

Table~\ref{T2d} shows all (decomposable or indecomposable) codes of lengths
$n \le 22$ with minimal distance $d \ge 4$, as obtained by subtracting
suitable codes $d_{2m}$ from one of the codes in Tables~\ref{T2b} and \ref{T2b2}.
The second column indicates the parent code in Tables~\ref{T2b} and \ref{T2b2} and the $d_{2m}$ to be subtracted.
The next two columns gives the components, with a \S~ to indicate a Type II code,
and the name (if any) given to this code in \cite{Ple10} or \cite{Me39}.
The remaining columns give the orders of the glue groups $G_1$ and $G_2$,
the weight distribution, and generators for the glue.

Table~\ref{T2e} gives the self-dual codes (both Type I and Type II) of length 24 and minimal distance $d \ge 4$.

A complete list of all Type I or Type II self-dual codes of lengths $n \le 24$ can be obtained by forming direct sums of the codes in Tables~\ref{T2d} and \ref{T2e} in all possible ways with the codes
$i_2^m$ $(m=0,1, \ldots )$.
\begin{table}[htb]
\caption{}
\begin{center}
\begin{tabular}{rllccrrrrrrl}
\multicolumn{12}{c}{Binary self-dual codes with $n \le 22$, $d \ge 4$} \\ \hline
\multicolumn{1}{c}{$n$} & \multicolumn{1}{c}{Code} & Compts. & Name & \multicolumn{1}{c}{$|G_1 |$} & \multicolumn{1}{c}{$|G_2 |$} & $A_4$ & $A_6$ & $A_8$ & $A_{10}$ & $A_{12}$ & Generators for glue \\ \hline
0 & C1$(d_{32} )$ & $i_0$ & - & 1 & 1 &  &  &  &  &  & - \\
8 & C2$(d_{24} )$ & $e_8$ & $A_8$ & 1 & 1 & 14 & 0 & 1 &  &  & - \\
12 & C3$(d_{20} )$ & $d_{12}$ & $B_{12}$ & 1 & 1 & 15 & 32 & 15 & 0 & 1 & $a$ \\
14 & C4$(d_{18} )$ & $e_7^2$ & $D_{14}$ & 1 & 2 & 14 & 49 & 49 & 14 & 0 & $d d$ \\
16 & C5$(d_{16} )$ & $d_{16}$ & $E_{16}$ & 1 & 1 & 28 & 0 & 198 & 0 & 28 & $a$ \\
 & C6($d_{16}$) & $e_8^2$ & $A_8 \oplus A_8$ & 1 & 2 & 28 & 0 & 198 & 0 & 28 & - \\
 & C7($d_{16}$) & $d_8^2$ & $F_{16}$ & 1 & 2 & 12 & 64 & 102 & 64 & 12 & $(ab)$ \\
18 & C8($d_{14}$) & $d_{10} e_7 f_1$ & $I_{18}$ & 1 & 1 & 17 & 51 & 187 & 187 & 51 & $aoA , cd$- \\
 & C9($d_{14}$) & $d_6^3$ & $H_{18}$ & 1 & 6 & 9 & 75 & 171 & 171 & 75 & $(abc), bbb$ \\
20 & C3($d_{12}$) & $d_{20}$ & $J_{20}$ & 1 & 1 & 45 & 0 & 210 & 512 & 210 & $a$ \\
 & C10($d_{12}$) & $d_{12} e_8$ & $A_8 \oplus B_{12}$ & 1 & 1 & 29 & 32 & 226 & 448 & 226 & $a$- \\
 & C11($d_{12}$) & $d_{12} d_8$ & $K_{20}$ & 1 & 1 & 21 & 48 & 234 & 416 & 234 & $(ab)$ \\
 & C12($d_{12}$) & $d_8^2 d_4$ & $S_{20}$ & 1 & 2 & 13 & 64 & 242 & 384 & 242 & $(ab)x , bby$ \\
 & C13($d_{12}$) & $e_7^2 d_6$ & $L_{20}$ & 1 & 2 & 17 & 56 & 238 & 400 & 238 & $doa, d d b$ \\
 & C14($d_{12}$) & $d_6^3 f_2$ & $R_{20}$ & 1 & 6 & 9 & 72 & 246 & 368 & 246 & $aaaA , cccB,  (abc)$- \\
 & C15($d_{12}$) & $d_4^5$ & $M_{20}$ & 1 & 120 & 5 & 80 & 250 & 352 & 250 & $(ooxyx)$ \\
22 & C8($d_{10}$) & $d_{14} e_7 f_1$ & $N_{22}$ & 1 & 1 & 28 & 49 & 246 & 700 & 700 & $aoA , b d A$ \\
 & C16($d_{10}$) & $d_{10}^2 f_2$ & $P_{22}$ & 1 & 2 & 20 & 57 & 270 & 676 & 676 & $(ao) \ast , cc$- \\
 & C17($d_{10}$) & $d_{10} d_6^2$ & $Q_{22}$ & 1 & 2 & 16 & 61 & 282 & 664 & 664 & $aoc , oaa , bbb$ \\
 & C18($d_{10}$) & $e_8 e_7^2$ & $ E_8 \oplus D_{14}$ & 1 & 2 & 28 & 49 & 246 & 700 & 700 & -$d d$ \\
 & C19($d_{10}$) & $d_8 e_7 d_6 f_1$ & $R_{22}$ & 1 & 1 & 16 & 61 & 282 & 664 & 664 & $odbA , boaA , aob$- \\
 & C20($d_{10}$) & $d_8 d_6^2 f_2$ & $S_{22}$ & 1 & 2 & 12 & 65 & 294 & 652 & 652 & $baoA , aoo AB , abb$-, \\
&&&&&&&&&&&$occ$- \\
 & C21($d_{10}$) & $d_6^2 d_4^2 f_2$ & $T_{22}$ & 1 & 4 & 8 & 69 & 306 & 640 & 640 & $aoxoA, ooyyAB$, \\
&&&&&&&&&&&$aayo$-, $bozx$-, $obxz$- \\
 & C22($d_{10}$) & $d_4^4 f_6$ & $U_{22}$ & 6 & 24 & 4 & 73 & 318 & 628 & 628 & $oxyz BC , ozxy AC$, \\
&&&&&&&&&&&$ooxx AE , oyoy AD$, \\
&&&&&&&&&&&$ozzo AF , xxxx$-, \\
&&&&&&&&&&&$yyyy$- \\
 & C23($d_{10}$) & $g_{22}$ & $G_{22}$ & 2 & 1 & 0 & 77 & 330 & 616 & 616 & the all-ones vector
\end{tabular}
\end{center}
\label{T2d}
\end{table}
\begin{table}[htb]
\caption{}
\begin{center}
\begin{tabular}{llllllll}
\multicolumn{8}{c}{Binary self-dual codes with length 24 and $d \ge 4$} \\ \hline
Code & Components & Name & $d$ & Code & Components & Name & $d$ \\ \hline
C2($e_8$) & $d_{24}$ \S & $E_{24}$ & 4 & C32($d_8$) & $d_8 e_7^2 f_2$ & $J_{24}$ & 4 \\
C6($e_8$) & $d_{16} e_8$ \S & --- & 4 & C33($d_8$) & $d_8 d_6^2 f_4$ & $R_{24}$ & 4 \\
C7($d_8$) & $d_{16} d_8$ & $H_{24}$ & 4 & C34($d_8$) & $d_8 d_4^4$ & $T_{24}$ & 4 \\
C10($e_8$) & $d_{12}^2$ \S & $A_{24}$ & 4 & C35($d_8$) & $e_7 d_6^2 d_4 f_1$ & $P_{24}$ & 4 \\
C11($d_8$) & $d_{12}^2$ & --- & 4 & C26($e_8$) & $d_6^4$ \S & $D_{24}$ & 4 \\
C12($d_8$) & $d_{12} d_8 d_4$ & $I_{24}$ & 4 & C36($d_8$) & $d_6^4$ & $Q_{24}$ & 4 \\
C18($e_8$) & $d_{10} e_7^2$ \S & $B_{24}$ & 4 & C37($d_8$) & $d_6^2 d_6 d_4 f_2$ & $S_{24}$ & 4 \\
C19($d_8$) & $d_{10} e_7 d_6 f_1$ & $K_{24}$ & 4 & C38($d_8$) & $d_6^2 d_4^2 f_2$ & $U_{24}$ & 4 \\
C20($d_8$) & $d_{10} d_6^2 f_2$ & $N_{24}$ & 4 & C39($d_8$) & $d_6 d_4^3 f_6$ & $W_{24}$ & 4 \\
C24($e_8$) & $e_8^3$ \S & --- & 4 & C27($e_8$) & $d_4^6$ \S & $F_{24}$ & 4 \\
C25($d_8$) & $e_8 d_8^2$ & --- & 4 & C40($d_8$) & $d_4^6$ & $V_{24}$ & 4 \\
C25($e_8$) & $d_8^3$ \S & $C_{24}$ & 4 & C41($d_8$) & $d_4^4 f_8$ & $X_{24}$ & 4 \\
C29($d_8$) & $d_8^3$ & $L_{24}$ & 4 & C42($d_8$) & $d_4^2 g_{16}$ & $Y_{24}$ & 4 \\
C30($d_8$) & $d_8^3$ & $M_{24}$ & 4 & C43($d_8$) & $h_{24}$ & $Z_{24}$ & 6 \\
C31($d_8$) & $d_8^2 d_4^2$ & $O_{24}$ & 4 & C28($e_8$) & $g_{24}$ \S & $G_{24}$ & 8 \\ 
\end{tabular}
\end{center}
\label{T2e}
\end{table}

There are over 1000 self-dual codes of lengths 26--30 (see Table~\ref{T2z}, \cite{cp}, \cite{Me162}).
The highest minimal distance is 6, and there are respectively 1, 3 and 13 codes
with $d=6$ of lengths 26, 28 and 30.
\subsection{Family \Etwo: Enumeration of ternary self-dual codes}\label{GU4}
\hsp
Ternary self-dual codes of lengths $n \le 20$ (and the maximal self-orthogonal codes of lengths $n \le 19$, $n \not\equiv 0$ $(\bmod~4)$ have been
enumerated by Pless \cite{Ple8} and Mallows, Pless and Sloane \cite{Me51} for $n \le 12$, Conway, Pless and Sloane \cite{Me65} for $n \le 16$, and Pless,
Sloane and Ward \cite{Me69} for $n \le 20$.
Leon, Pless and Sloane \cite{Me76} give a partial enumeration of the self-dual codes of length 24,
making use of the complete list of Hadamard
matrices of order 24,
and show that there are precisely two codes with minimal
distance 9 (cf. Table~\ref{TX3a} below).

We will make use of the following components.

$e_3$: [111], glue: $\pm a$, $a=120$.
If the coordinates are labeled 1, 2, 3 then $G_0$
is generated by $(1,2,3)$ and $(1,2)$ ${\rm diag} \{-1, -1, -1\}$ and has order 6;
$|G_1| =2$.

$t_4$ is the $[4,2,3]_3$ tetracode, and $g_{12}$ is the $[12,6,6]_3$ ternary
Golay code, see Section~\ref{Wee2}.

$g_{10}$ is the $[10,4,6]_3$ code consisting of the vectors $u$ such that $00 u \in g_{12}$.
If $x$ and $y$ are chosen so that $11 x \in g_{12}$, $12y \in g_{12}$, then the glue words for $g_{10}$ can be taken to be $\pm x$, $\pm y$, $\pm x \pm y$.
$|G_0| = 360$, $|G_1| =8$.

$p_{13}$:
Let $Q_0$, $Q_1 , \ldots, Q_{12}$ be the points of a projective plane\index{projective plane} of order 3,
labeled so that the 13 lines are represented by the cyclic shifts $t_0$, $t_1, \ldots, t_{12}$ of the vector $t_0$ given by
$$\begin{array}{ccccccccccccc}
Q_0 & Q_1 & Q_2 & Q_3 & Q_4 & Q_5 & Q_6 & Q_7 & Q_8 & Q_9 & Q_{10} & Q_{11} & Q_{12} \\
1 & 1 & 0 & 1 & 0 & 0 & 0 & 0 & 0 & 1 & 0 & 0 & 0
\end{array}
$$
(\cite{MS}, p.~695, \cite{Me146}).
The vectors $t_0 , \ldots, t_{12}$ generate a $[13,7,4]_3$ code $p_{13}^\perp$.
The dual is $p_{13}$, a $[13,6,6]_3$ self-orthogonal code consisting of the vectors $\sum_{i=0}^{12} a_i t_i$ with $a_i \in \FF_3$ and $\sum a_i =0$, and having weight distribution
$A_0 =1$, $A_6 = 156$, $A_9 = 494$, $A_{12} = 78$.
$G_0 ( p_{13}) \simeq PGL_3 (3)$, of order 5616, $|G_1 (p_{13} ) | =2$.
The glue words are $\pm t_0$.

The indecomposable self-dual codes of lengths $n \le 16$ are shown in Table~\ref{T3a}.
\index{matrix!Hadamard}
$H_8$ denotes a suitably normalized version of the Hadamard matrix\index{Hadamard matrix} of order 8.
\begin{table}[htb]
\caption{Indecomposable ternary self-dual codes of lengths $n \le 20$}
$$
\begin{array}{c|c|c|c|c|c|l} \hline
n & \mbox{Components} & |G_0 | & |G_1 | & |G_2| & d & \mbox{glue} \\ \hline
4 & t_4 & 48 & 1 & 1 & 3 & - \\ 
8 & - & ~ & ~ & ~ & ~ & ~ \\ 
12 & e_3^{4+} & 6^4 & 2 & 24 & 3 & aaa0, 0 \bar{a} aa \\ 
~ & g_{12} & 190080 & 1 & 1 & 6 & - \\ 
16 & (e_3^4 f_4)^+ & 6^4 .1 & 8 & 24 & 3 & (a000) (2111) \\ 
~ & (e_3^2 g_{10} )^+ & 6^2 .360 & 4 & 2 & 3 & a0x , 0ay \\ 
~ & (e_3 p_{13})^+ & 6.5616 & 2 & 1 & 3 & at_0 \\ 
~ & f_8^{2+} & 1^2 & 2^7 .168 & 2 & 6 & [I | H_8] \\ 
20 & \mbox{17 codes} & ~ & ~ & ~ & ~ & ~ \\
~ & \mbox{-- see \cite{Me69}} & ~ & ~ & ~ & ~ & ~ \\
\end{array}
$$
\label{T3a}
\end{table}

The analogue of Theorem~\ref{thGU2} is:
any self-orthogonal ternary code generated by words of weight 3 is a direct
sum of copies of $e_3$ and $t_4$.
A technique for classifying self-orthogonal codes generated by words
of weight 6 (using ``center sets'')\index{center set} is given in \cite{Me69}.
\subsection{Family \Ethree: Enumeration of Hermitian self-dual codes over $\FF_4$}\label{GU5}
\hsp
These have been classified for lengths $n \le 16$ \cite{Me65}
--- see Table~\ref{T4a}.

We will make use of the following components.

$d_{2n}$ $(n \ge 2)$:
generated by \eqn{EqGU8}.
There are 16 cosets of $d_{2n}$ in $d_{2n}^\perp$, and as glue words we choose
$0$, $\om^\nu a$, $\om^\nu b$, $\om^\nu c$, $\om^\nu d$, $\om^\nu e$,
$\nu \in \{0,1,2 \}$, where
\begin{eqnarray*}
a & = & 1010 \ldots 1010 \\
b & = & 0000 \ldots 0011 \\
c & = & 1010 \ldots 1001 \\
d & = & 1010 \ldots 10 \om \oom \\
e & = & 1010 \ldots 10 \oom \om
\end{eqnarray*}
Also $|G_0| = 2^{n-1} n!$, $|G_1| = 36$ $(n=2)$, or 12 $(n \ge 3$).

$e_5: [ \om \oom \oom \om 0, 0 \om \oom \oom \om ]$,
glue: $\om^\nu 1$, $\nu \in \{0,1,2\}$, $G_0 = A(5)$,\index{group!alternating $A(n)$} of order 60,
$|G_1| =6$.

$h_6$ is the hexacode, $e_7 \otimes \FF_4$, $e_8 \otimes \FF_4$ are $\FF_4$-versions of the Hamming codes in Sections~\ref{Wee2},
and $1_n$ is the $[n,1,n]_4$ repetition code.
\begin{table}[htb]
\caption{Indecomposable Hermitian self-dual codes over $\FF_4$ of lengths $n \le 16$}
$$
\begin{array}{c|c|c|c|c|c|c} \hline
n & \mbox{Components} & |G_0 | & |G_1 | & |G_2 | & d & \mbox{Glue} \\ \hline
2 & i_2 & 12 & 1 & 1 & 2 & - \\
4 & - & - & - & - & - & - \\
6 & h_6 & 2160 & 1 & 1 & 4 & - \\
8 & e_8 & 8064 & 1 & 1 & 4 & - \\
10 & d_{10}^+ & 2^4 . 5! & 6 & 1 & 4 & d \\
~ & e_5^{2+} & 60^2 & 6 & 2 & 4 & 11 \\
12 & d_{12}^+ & 2^5 . 6! & 6 & 1 & 4 & a \\
~ & (e_7 e_5)^+ & 60.168 & 6 & 1 & 4 & 11 \\
~ & d_6^{2+} & 24^2 & 6 & 2 & 4 & (bd) \\
~ & d_4^{3+} & 4^3 & 54 & 6 & 4 & (0de) \\
14 & d_{14}^+ & 2^6 . 7 ! & 6 & 1 & 4 & d \\
~ & e_7^{2+} & 168^2 & 6 & 2 & 4 & 11 \\
~ & (d_8 e_5 f_1)^+ & 2^3 .4! 60 & 6 & 1 & 4 & d01, e10 \\
~ & (e_5^2 d_4)^+ & 4.60^2 & 18 & 2 & 4 & 01d, 10e \\
~ & (d_8 d_6)^+ & 2^3 4! 2^2 3! & 6 & 1 & 4 & ab, bd \\
~ & (d_6^2 f_2 )^+ & (2^2 .3!)^2 & 6 & 2 & 4 & (d0) 11, bb \om \oom \\
~ & (d_6 d_4^2)^+ & 4^2 2^2 3! & 18 & 2 & 4 & bbb, a0d, cd0 \\
~ & (d_4^3 f_2)^+ & 4^3 & 6 & 6 & 4 & aa011, 0aa \om \om , \dot{b} \ddot{b} 0 \om \oom , 0 \dot{b} \ddot{b} \oom 1 \\
~ & (d_4^2 1_6)^+ & 4^2 & 108 & 2 & 4 & b00011 \om \om, a00 \oom \oom 11 0, \\
&&&&&&0b 01 \om 0 1 \om , 0a 0 11 \oom \oom 0 \\
~ & q_{14} & 6552 & 1 & 1 & 6 & -  \\
16 & \mbox{31 codes} & ~ & ~ & ~ & ~ & \\
~ & \mbox{(see \cite{Me65})} & ~ & ~ & ~ & ~ & ~
\end{array}
$$
\label{T4a}
\end{table}
\paragraph{Remarks.}
(1) The group orders differ slightly form those in \cite{Me65}, since now we are allowing conjugation in the group.

(2)~The dots and double-dots in the glue column indicate multiplication by
$\omega$ or $\omega^2$, respectively.

(3)~The unique distance 6 code at length 14, $q_{14}$, is the $[14, 7,6]_4$ extended quadratic residue code generated by
$$1(1 \om \oom \om \om \oom \oom \oom\oom \om \om \oom \om ) ~.$$

(4)~The analogue of Theorem~\ref{thGU2} is:
(a) any self-orthogonal code with minimal distance 2 has $i_2$ as a direct summand;
(b) any self-orthogonal code generated by words of weight 4 is a direct sum of copies of $d_4$, $d_6$, $d_8 , \ldots$, $e_5$, $h_6$, $e_7$ and $e_8$.
\subsection{Family \Efour: Enumeration of Euclidean self-dual codes over $\FF_4$}\label{GU6}
\hsp
Although even codes of length up to 14 were classified in \cite{Me55}, the odd codes do not seem to have been classified.
\subsection{Family \Efive: Enumeration of trace self-dual additive codes over $\FF_4$}\label{GU7}
\hsp
These have been classified up to length 7
(and Type II code up to length 8)
in \cite{Me223}, \cite{Hoh96}.

The analogue of Theorem~\ref{thGU2} is the following.
Let $d_n$ be the code of length $n$ generated by all even-weight binary
vectors $(n \ge 2)$, and let $i_2 = [11, \om \om ]$.
Then any trace self-orthogonal additive code over $\FF_4$ generated by words
of weight 2 is a direct sum of copies of $i_2$, $d_2$, $d_3$, $d_4, \ldots$.

$d_n^+$ (mentioned in Table~\ref{TX5a}) is the code of length $n$,
containing $2^n$ words, generated by $d_n$ and $\om \om \ldots \om$.
\subsection{Family \Eeight: Enumeration of self-dual codes over $\ZZ_4$}\label{GU8}
\hsp
These have been classified for lengths up to 16 in the following papers:
 Conway and Sloane \cite{Me168} for $n \le 9$, Fields, Gaborit, Leon and Pless \cite{FGLP} for $n \le 15$, and Pless, Leon and Fields \cite{PLF96} for Type II codes of length 16.

In this section we will present enough component codes to state the analogue of Theorem~\ref{thGU2}.

The smallest self-dual code is $i_1 = \{0,2\}$.
If a self-orthogonal code $C$ contains a vector of the form $2^1 0^{n-1}$ then
$C= i_1 \oplus C'$ is
decomposable.
The next-simplest possible vectors are ``tetrads'',\index{tetrad} of type $\pm 1^4 0^{n-4}$.
We list a number of self-orthogonal codes that are generated by tetrads;
$t$ denotes the total number of tetrads in the code.

The first four codes have the property that the associated binary code
$C^{(1)}$ is the self-dual code $d_{2m}$ of \eqn{EqGU8}.

$\sD_{2m}$
$(m \ge 2)$ is generated by the tetrads
$11130 \ldots 0$, $0011130 \ldots 0 ,~ \ldots ,~ 0 \ldots 01113$;
$| \sD_{2m} |  =  4^{m-1}$,
$| Aut ( \sD_{2m}) | = 2 . 4!~ ( m = 2 )$ or $2^2 . 2^m ~ ( m > 2 )$, $t = 2 (m-1)$.
$\sD_{2m}^\perp / \sD_{2m}$ is a group of type $4^2$ with generators $v_1  =  0101 \ldots 01$, $v_2  =  00 \ldots 0011$.

$\sD_{2m}^{\rm O}$ $( m \ge 2 )$ is generated by $\sD_{2m}$ and the tetrad $1300 \ldots 0011$ (or equivalently the vector $2020 \ldots 20$);
$| \sD_{2m}^{\rm O} |  =  4^{m-1} 2$, $|Aut ( \sD_{2m}^{\rm O})| = 2^2  . 8$ $ ( m=2 )$ or
$2 . 2^{m-1} . 2m ~( m  > 2 )$,
$t=2m$.
$( \sD_{2m}^{\rm O} )^\perp / \sD_{2m}^{\rm O}$ is a cyclic group of order 4 generated by $v_1$ (if $m$ is odd), or a 4-group generated by $v_1$ and $2 v_2$ (if $m$ is even).

$\sD_{2m}^+$
$( m \ge 2$, but note that $\sD_4^+  \simeq \sD_4^{\rm O}$) is generated by $\sD_{2m}$ and $2 v_2$;
$| \sD_{2m}^+ | =  4^{m-1} 2$,
$|Aut ( \sD_{2m}^+ )| = 2^m . 2^{m+1}$, $t = 4(m-1)$.
$( \sD_{2m}^+ )^\perp / \sD_{2m}^+$ is a 4-group generated by $2 v_1$ and $v_2$.

$\sD_{2m}^\oplus$
$(m \ge 2 )$ is the self-dual code generated by $\sD_{2m}^{\rm O}$ and $\sD_{2m}^+$;
$| \sD_{2m}^\oplus |  =  4^{m-1} 2^2$, $|Aut ( \sD_{2m}^\oplus )| = 2^3  . 4!$
$(m=2)$ or $2^m . 2^m . 2m$ $(m > 2)$,
$t=4m$.
For use in \eqn{EqG72} we note that there are two permutation-inequivalent versions of $\sD_4^\oplus$, with generator matrices
\beql{EqC23a}
{\rm (a)} ~ \left[ \matrix{
1 & 1 & 1 & 1 \cr
0 & 2 & 0 & 2 \cr
0 & 0 & 2 & 2 \cr
} \right]
, \quad {\rm (b)}~
\left[ \matrix{
1 & 3 & 3 & 3 \cr
0 & 2 & 0 & 2 \cr
0 & 0 & 2 & 2 \cr
}\right] ~.
\eeq
$\sD_4^\oplus$ (in either version) has
$swe  =  x^4 + 6x^2 z^2 + z^4 + 8y^4$.

$\sE_7$ is generated by 1003110, 1010031, 1101003;
$| \sE_7 |  = 4^3$, $|Aut ( \sE_7)| = 2 . 4!$, $t=8$.
$\sE_7^\perp / \sE_7$ is a cyclic group of order 4 generated by 3111111.

$\sE_7^+$ is the self-dual code generated by $\sE_7$ and 2222222
(or equivalently by all cyclic shifts of 3110100);
$| \sE_7^+ | =  4^3 2$,
$|Aut ( \sE_7^+)| = 2 . 168$, $t=14$,
$swe  = x^7 + z^7 + 14y^4 (x^3 + z^3 )$
$+ 7x^3 z^3 ( x+z) + 42 xy^4 z (x+z)$.
For both $\sE_7$ and $\sE_7^+$ the associated binary code $C^{(1)}$ is the Hamming code $e_7$.

$\sE_8$ is the self-dual code generated by $0u$, $u \in \sE_7$ and 30001011.
An equivalent generator matrix has already been given in \eqn{EqLift1}.
$| \sE_8 | =  4^4$, $g = 8 . 2 . 4! = 384$, $t=16$, $swe = x^8 + 16y^8 + z^8 + 16y^4 ( x^4 + z^4 ) + 14 x^4 z^4 +$ $48 x y^4 z ( x^2 + z^2 )$ $+ 96 x^2 y^4 z^2$.
\begin{theorem}\label{thGUZ4}
Any self-orthogonal code over $\ZZ_4$ generated by vectors of the form
$\pm 1^4 0^{n-4}$ is equivalent to a direct sum of copies of the codes
$$\sD_{2m},~ \sD_{2m}^{\rm O} ,~ \sD_{2m}^+ ,~ \sD_{2m}^\oplus (m=1,2, \ldots ),~ \sE_7 ,~ \sE_7^+ ,~ \sE_8 ~.$$
\end{theorem}

The (somewhat complicated) inclusions between the codes mentioned in the theorem can be seen in Fig.~1 of \cite{Me168}.
\section{Extremal and optimal self-dual codes}\label{EX}
\hsp
Recall from Section~\ref{BD} that we have defined a self-dual code from any
of the families \Eone\ through \Eseven\ to be {\em extremal}\index{extremal code} if it meets the strongest
\index{code!extremal}
of the applicable bounds from Theorems~\ref{thBD4}, \ref{thBD5} and \ref{thBD7}, that is,
if its minimal distance $d$ is equal to

(\EoneI)
$4 \left[ \frac{n}{24} \right] +4 + \ep$, where $\ep = -2$ if $n=2$, 4 or 6, $\ep =2$ if $n \equiv 22$ $(\bmod~24)$, and $\ep =0$ otherwise,

(\EoneII) $4 \left[ \frac{n}{24} \right] + 4$,

(\Etwo) $3 \left[ \frac{n}{12} \right] + 3$,

(\Ethree) $2 \left[ \frac{n}{6} \right] +2$,

(\Efour) $\left[ \frac{n}{2} \right] +1$,

(\EfiveI) $2 \left[ \frac{n}{6} \right] +2 + \ep'$, where $\ep' = -1$ if $n=1$, $\ep ' =1$ if $n \equiv 5$ $(\bmod~6)$, and $\ep' =0$ otherwise,

(\EfiveII) $2 \left[ \frac{n}{6} \right] +2$,

(\Esix), (\Eseven) $\left[ \frac{n}{2} \right] +1$.

We also defined a code over $\ZZ_4$ to be norm-extremal if its minimal norm is

(\Eeight) $8 \left[ \frac{n}{24} \right] + 8 + \ep ''$

\noindent
where $\ep '' =4$ if $n \equiv 23$ $(\bmod~4)$, $\ep '' =0$ otherwise.

It is very likely (although we do not have a proof) that the above bounds
for families \Eone\ through \Eseven\ are the highest minimal distance that is
permitted by the pure linear programming bound applied to the Hamming weight
enumerator and (when relevant) the shadow enumerator.

In contrast, we call a code {\em optimal} if it has the highest
\index{code!optimal}\index{optimal code}
minimal distance of any self-dual code of that length.
An extremal code is automatically optimal.

In this section we will summarize what is presently known about extremal and optimal codes in the families we are considering.
Earlier summaries of extremal codes and lattices have appeared in
Chapter~7 of \cite{SPLAG}, \cite{huff3.5}.
In the tables we have tried to list all known codes
with the specified minimal distance (a period indicating that the list is
complete), or else to indicate how many extremal codes are known.
Whenever possible we have attempted to name at least one extremal code.
\subsection{Family \Eone: Binary codes}\label{EX1}
\hsp
Type I codes meeting the $d \le 2 [n/8] +2$ bound of Theorem~\ref{thBD4}
(the old definition of extremal) were completely
classified by Ward \cite{Wa76} (finishing the work begun in
\cite{Me30}, \cite{Ple10}, \cite{Me39}):
such codes exist if and only if $n$ is $2~(i_2)$, $4~( i_2^2)$, $6~(i_2^3)$,
$8~(e_8)$, $12 ~( d_{12}^+ )$, $14~ (e_7^{2+} )$, $22 ~(g_{22}^+)$ or $24 ~(g_{24} )$ --- compare Tables~\ref{T2a} and \ref{T2d}.
In each case the code is unique.

However, there are many more Type I codes that are extremal in the new sense, and they have not yet been fully classified.
It is known (Theorem~\ref{thBDa}) that extremal Type II codes do not exist for lengths $\ge 3952$ and
presumably a similar bound applies to extremal Type I codes.

Table~\ref{TX2a} shows the highest possible
minimal distance for binary self-dual codes of lengths $n \le 72$.
This is based on earlier tables in Fig.~19.2 of \cite{MS}, \cite{Me158} and \cite{DGH97a}.
In the table $d_I$ (resp. $d_{II}$) denotes the highest minimal distance of any
strictly
Type~I (resp. Type~II) self-dual code.
\subsection*{Remarks on Table~\ref{TX2a}}
\hsp
The fourth column of the table
gives the known codes having
the indicated minimal
distance.
As mentioned above, a
period indicates that the lists of codes is complete.
(The enumeration for lengths $n \le 32$ has already been discussed
in Section~\ref{GU3}.)
When $n$ is a multiple of 8 a semicolon separates the Type~I and Type~II codes.
\begin{table}[htb]
\caption{Highest minimal distance of binary self-dual codes}
\begin{center}
\begin{tabular}{rccl} \hline
$n~$ & $d_I$ & $d_{II}$ & Codes \\ \hline
2 & 2 &  & $i_2$. \\ 
4 & 2 &  & $i_2^2$. \\ 
6 & 2 &  & $i_2^3$. \\
8 & 2 & 4 & $i_2^4$; $e_8$. \\
10 & 2 &  & $i_2^5$, $e_8 i_2$. \\
12 & 4 &  & $d_{12}^+$. \\
14 & 4 &  & $e_7^{2+}$. \\
16 & 4 & 4 & $d_8^{2+}$; $d_{16}^+$, $e_8^2$. \\
18 & 4 &  & $d_6^{3+}$, $( d_{10} e_7 f_1 )^+$. \\
20 & 4 &  & 7 codes (Table~\ref{T2a}). \\
22 & 6 &  & $g_{22}$. \\
24 & 6 & 8 & $h_{24}^+$; $g_{24}$. \\
26 & 6 &  & $f_{13}^{ 2}$ \cite{cp}. \\
28 & 6 &  & 3 codes \cite{cp}. \\
30 & 6 &  & 13 codes \cite{cp}, \cite{Me162}. \\
32 & 8 & 8 & 3 codes \cite{Me158}; 5 codes (Table~\ref{T2b}). \\
34 & 6 &  & $\ge 200$ \\
36 & 8 &  & $\ge 2$ \\
38 & 8 & ~ & $\ge 3$ ~~~\cite{Me158}, \cite{HaKi95a} \\
40 & 8 & 8 & $\ge 22$; $\ge 1000$ (see text for references) \\
42 & 8 & ~ & $\ge 30$ \cite{DGH97a} \\
44 & 8 & ~ & $\ge 108$ \cite{DGH97a} \\
46 & 10 & ~ & $\ge 1$ \cite{Me158} \\
48 & 10 & 12 & $\ge 7$; $\ge 1$ $(XQ_{47})$ \\
50 & 10 & ~ & $\ge 6$ \\
52 & 10 & ~ & $\ge 499$ \cite{ht2} \\
54 & 10 & ~ & $\ge 54$ \\
56 & 10 or 12 & 12 & ?; $\ge 166$ \\
58 & 10 & ~ & $\ge 80$ \cite{DGH97a} \\
60 & 12 & ~ & $\ge 5$ \\
62 & 10 or 12 & ~ & ? \\
64 & 12 & 12 & $\ge 5$; $\ge 3270$ \cite{DGH97a} \\
66 & 12 & ~ & $\ge 3$ \\
68 & 12 & ~ & $\ge 65$ \\
70 & 12 or 14 & ~ & ? ~ \cite{Hara96a}, \cite{SchSch} \\
72 & 12 or 14 & 12 or 16 & ?; ?
\end{tabular}
\end{center}
\label{TX2a}
\end{table}

In the years since the manuscript of
\cite{Me158} was first circulated, a large number of sequels have been
written, supplying additional examples of self-dual codes in the range of Table~\ref{TX2a}.
The bibliography includes all the manuscripts known to us, even though
inevitably not all of them will be published.
It was not possible to mention all these references in the table, so instead
we list them here.
This list also includes a number of older papers.
Readers interested in extremal self-dual codes, especially of Type I, in the range of the table should therefore consult the following:
\cite{BrPl91}, \cite{Bus2},
\cite{Bus3}, \cite{Buy97}, \cite{Buy97a},
\cite{BB97},
\cite{by}, \cite{Dou95},
\cite{DGH97a}, \cite{DoHa96}, \cite{DoHa97},
\cite{DoHa97b},
\cite{DHo97},
\cite{gb92}, \cite{gh0}, \cite{gh0.1}, \cite{gh1},
\cite{gh2}, \cite{GuHa97a}, \cite{GHK97},
\cite{Hara96}, \cite{Hara96a}, \cite{Hara97a}, \cite{Hara97d},
\cite{HaKi95},
\cite{HaKi95a},
\cite{Ton15}, \cite{huff0}, \cite{ht},
\cite{ht2}, \cite{hy}, \cite{KaT90},
\cite{Kim94}, \cite{Oze98},
\cite{Pasq81},
\cite{Ple14},
\cite{Ple94}, \cite{Ton44},
\cite{rus1}, \cite{rus2}, \cite{rus3},
\cite{Ton6}, \cite{Ton11}, \cite{Ton16},
\cite{Tsai91}, \cite{Tsai92},
\cite{Tsai97},
\cite{yorgov1}, \cite{yorgov2}, \cite{yorgov3}, \cite{yr},
\cite{yy}, \cite{yz}.

Note that if we don't distinguish between Type I and Type II codes, but just ask what is the highest minimal distance of a binary self-dual code, then the answer is known
for all $n \le 60$.

The symbol $XQ_m$ in any of these tables indicates an extended quadratic residue code\index{quadratic residue code} of length $m+1$.
\index{code!double circulant}\index{code!quadratic residue}
Both quadratic residue codes and double circulant codes provide many examples
of good self-dual codes
(cf. Section~12 of Chapter~1, Chapter~xx (Ward), Chapter~yy (Pless),
\cite[Chapter~16]{MS}).
There are two basic types of binary double circulant codes,\index{double circulant code} having generator
matrices of the form
\beql{EqD1}
\begin{array}{|ccccc|ccccc|} \hline
1 & ~ & ~ & ~ & ~ & 0 & 1 & 1 & 1 & 1 \\
~ & 1 & ~ & ~ & ~ & 1 & ~ & ~ & ~ & ~ \\
~ & ~ & 1 & ~ & ~ & 1 & ~ & R & ~ & ~ \\
~ & ~ & ~ & 1 & ~ & 1 & ~ & ~ & ~ & ~ \\
~ & ~ & ~ & ~ & 1 & 1 & ~ & ~ & ~ & ~ \\
\hline
\end{array}
\eeq
or
\beql{EqD2}
\begin{array}{|ccccc|ccccc|} \hline
1 & ~ & ~ & ~ & ~ & ~ & ~ & ~ & ~ & ~ \\
~ & 1 & ~ & ~ & ~ & ~ & ~ & ~ & ~ & ~ \\
~ & ~ & 1 & ~ & ~ & ~ & ~ & R & ~ & ~ \\
~ & ~ & ~ & 1 & ~ & ~ & ~ & ~ & ~ & ~ \\
~ & ~ & ~ & ~ & 1 & ~ & ~ & ~ & ~ & ~ \\
\hline
\end{array}~,
\eeq
where $R$ is a circulant matrix with first row $r$ (say).
\eqn{EqD1} is used only when the length is a multiple of 4.
Such codes and their generalizations to other fields have
been studied by many authors, including
\cite{Bee1},
\cite{Bhar80}, \cite{gh0}--\cite{GHK97}, \cite{huff3.5},
\cite{hy},
\cite{Karlin},
\cite{MS878},
\cite{MS881},
\cite[Chap.~16]{MS},
\cite{Poli}, \cite{rus1},
\cite{TonR4},
\cite{VeRi85},
\cite{yorgov1}--\cite{yorgov3}.
Table~\ref{TaDC}, based on \cite{Me158} and \cite{Moo3},
gives a selection of double circulant binary codes.
Code H86 (from \cite{DGH97a})
is the shortest Type I self-dual code presently known with $d=16$.
The first column gives the name of the codes, following \cite{Me158}, and the last column gives $r$, the initial row of $R$,
in hexadecimal.
The codes marked $(\ast)$ are not necessarily optimal.
The minimal distance of the last two codes in the table was determined by
Moore \cite{Moo3a}, \cite{Moo3}.
For these two codes $r$ has 1's at the squares modulo 43 and 67,
respectively.
Moore remarked that the analogous code of length 168 {\em might} also
have been extremal. However, Aaron Gulliver (personal communication, 
Nov. 1997) has shown that the minimal distance of this code
is at most 28.
\begin{table}
\caption{Double circulant binary codes}
\label{TaDC}
\begin{center}
\begin{tabular}{ccccccc} \hline
Name & $n$ & $k$ & $d$ & Type & Form & $r$ (hexadecimal) \\ \hline
$g_{22}$ & 22 & 11 & 6 & I & \eqn{EqD2} & 97 \\
$g_{24}$ & 24 & 12 & 8 & II & \eqn{EqD1} & B7 \\
$A_{26} = f_{13}^{2+}$ & 26 & 13 & 6 & I & \eqn{EqD2} & 5F7 \\
$A_{28} = $D1 & 28 & 14 & 6 & I & \eqn{EqD1} & 8D \\
D2 & 34 & 17 & 6 & I & \eqn{EqD2} & 1ECE \\
D3 & 36 & 18 & 8 & I & \eqn{EqD1} & 2C6B \\
D4 & 38 & 19 & 8 & I & \eqn{EqD2} & 5793 \\
D5 & 40 & 20 & 8 & II & \eqn{EqD2} & 57EB \\
D6 & 40 & 20 & 8 & I & \eqn{EqD2} & 11E35 \\
D7 & 40 & 20 & 8 & I & \eqn{EqD2} & B393 \\
D8 & 44 & 22 & 8 & I & \eqn{EqD1} & 5E6B5 \\
D9 & 50 & 25 & 10 & I & \eqn{EqD2} & 31C4D \\
D10 & 52 & 26 & 10 & I & \eqn{EqD1} & 57F69D \\
D11 & 56 & 28 & 12 & II & \eqn{EqD1} & ADF1FF \\
D12 & 58 & 29 & 10 & I & \eqn{EqD2} & D5A89B \\
D12a & 58 & 29 & 10 & I & \eqn{EqD2} & 2DD1D3 \\
D13 & 60 & 30 & 12 & I & \eqn{EqD1} & 3EF6B77 \\
D14 & 64 & 32 & 12 & II & \eqn{EqD1} & 427BD0B \\
D15 & 64 & 32 & 12 & I & \eqn{EqD2} & 2EF3DD75 \\
D16 & 66 & 33 & 12 & I & \eqn{EqD2} & B2D97D9 \\
D17 & 68 & 34 & 12 & I & \eqn{EqD2} & 1F5C885F \\
D18$^\ast$ & 72 & 36 & 12 & I & \eqn{EqD2} & 2B8795E5 \\
D19$^\ast$ & 74 & 37 & 12 & I & \eqn{EqD2} & 1439372C7 \\
D20$^\ast$ & 82 & 41 & 12 & I & \eqn{EqD2} & A464B919B \\
H86 & 86 & 43 & 16 & I & \eqn{EqD2} & 7F7101712E2 \\
M88 & 88 & 44 & 16 & II & \eqn{EqD1} & ~ \\
M136 & 136 & 68 & 24 & II & \eqn{EqD1} & ~ \\
\end{tabular}
\end{center}
\end{table}

We see from Table~\ref{TX2a} that there are extremal Type I codes (in the new sense) that are not also Type II codes at lengths
$$
2,4,6,12,14,16,18,20,22,32,36,38,40,42,44,46,60,64,66,68
$$
that such codes do not exist at length
\beql{EqGU99}
8, 10, 24, 26, 28, 30, 34, 48, 50, 52, 54, 58
\eeq
and that their existence at lengths
\beql{EqGU100}
56, 62, 70, 72
\eeq
is at present an open question.
The nonexistence of the Type I codes of lengths in \eqn{EqGU99} is established by imposing
the extra condition that the shadow enumerator must have integral coefficients.

Concerning extremal Type II codes, with $d=4 [n/24]+4$, these exist for the following values of $n$:
$$8, 16, 24, 32, 40, 48, 56, 64, 80, 88, 104, 136$$
but their existence at lengths 72 and 96 and all greater lengths is open.
For lengths 8, 24, 32, 48, 80 and 104 we can use extended quadratic
residue codes, and for lengths 40, 56, 64, 88, 136 we can use double circulant
codes (see Table~\ref{TaDC}).

Only one $[48,24,12]$ code is presently known, $XQ_{47}$,
which is generated by {\bf 1} and
$$
1(01111011110010101110010011011000101011000010000)
$$
(with 1's at the nonzero squares modulo 47).
Huffman \cite{huff0} has shown that any Type II $[48, 24, 12]$ code with a nontrivial automorphism of odd order is equivalent to $XQ_{47}$.
Houghten, Lam and Thiel (cf. \cite{HLT94})
are attempting to establish by
direct search that $XQ_{47}$ is unique.
As Table~\ref{TX2a} shows, if $n \ge 40$ is congruent to 8 or 16 ($\bmod~24$)
there are often large numbers of extremal codes.
It is easy to find $[72, 36, 12]$ Type II codes, for example $XQ_{71}$;
\cite{DGH97a} shows that there are at least 33 inequivalent codes with these parameters.

Concerning the existence of self-dual codes with a specified minimal distance, the following results were established in \cite{Me158}.
Self-dual codes with minimal distance

$d \ge 6$ exist precisely for $n \ge 22$;

$d \ge 8$ exist precisely for $n = 24, 32$, and $n \ge 36$;

$d \ge 10$ exist precisely for $n \ge 46$;

$d \ge 12$ exist\footnote{The existence of a $[70,35,12]$ was not known when
\cite{Me158} was written, but such a code was later found by Scharlau and Schomaker \cite{SchSch}.} for $n=48$, 56, 60 and $n \ge 64$;
perhaps for $n=62$; and do not exist for all other values of $n$.
(As pointed out in \cite{Me158}, the $[58,29,12]$ self-dual code claimed
in \cite{Bhar80} is an error.)

Dougherty, Gulliver and Harada \cite{DGH97a}, extending work in \cite{Me158}, show that codes with

$d \ge 14$ exist for $n \ge 78$;
perhaps for $n=70$, 72, 74, 76;
and do not exist for all other values of $n$;

$d \ge 16$ exist for $n=80$, 86, 88, 96, 100--104, 112 and $n \ge 120$ (and
possibly for other values of $n$).
\subsection{Family \Etwo: Ternary codes}\label{EX2}
\hsp
Table~\ref{TX3a} shows the highest possible
minimal distance for ternary self-dual codes of lengths $n \le 72$.
\begin{table}[htb]
\caption{Highest minimal distance of ternary self-dual codes}
\begin{center}
\begin{tabular}{rcl} \hline
$n$ & \multicolumn{1}{c}{$d$} & Codes \\ \hline
4 & 3 & $t_4$. \\
8 & 3 & $t_4^2$. \\
12 & 6 & $g_{12}$. \\
16 & 6 & $f_8^{2+}$. \\
20 & 6 & 6 codes \cite{Me69}. \\
24 & 9 & $XQ_{23}$, $S(24)$ \cite{Me76}. \\
28 & 9 & $\ge 32$ \cite{ChSch87}, \cite{Hara97b}, \cite{KP}, \cite{huff3} \\
32 & 9 & $\ge 239$ \cite{huff3} \\
36 & 12 & $\ge 1$ $(S(36))$ \\
40 & 12 & $\ge 20$ \cite{Wa85}, \cite{Dawso84}, \cite{Hara97b}, \cite{huff3} \\
44 & 12 & $\ge 8$ \cite{Hara97b} \\
48 & 15 & $\ge 2$ $(XQ_{47}, S(48))$ \\
52 & \multicolumn{1}{c}{12 or 15} & ? \\
56 & 15 & $\ge 1$ \\
60 & 18 & $\ge 2$ $(XQ_{59}, S(60))$ \\
64 & 18 & $\ge 1$ \cite{Bee1}, \cite{Dawso84} \\
68 & \multicolumn{1}{c}{15 or 18} & ? \\
72 & 18 & $\ge 1$ $(XQ_{71} \cite{CoSe84})$
\end{tabular}
\end{center}
\label{TX3a}
\end{table}
\subsection*{Remarks on Table~\ref{TX3a}}
\hsp
For the entries at lengths $n \le 24$, see the discussion in Section~\ref{GU4}.

Extremal codes exist at lengths 4, 8, 12, 16, 20, 24, 28, 32, 36, 40, 44, 48, 56, 60 and 64.
Extremal codes do not exist at lengths 72, 96, 120 and all $n \ge 144$,
because then the extremal Hamming weight enumerator contains a negative
coefficient.
The existence of extremal codes in the remaining cases ($n= 52$, 68, 76, $\ldots$, 140) is undecided.

In Table~\ref{TX3a}, $XQ_n$ denotes an extended quadratic residue code of length $n+1$, and $S(n)$ denotes a Pless double circulant
(or ``symmetry'') code of length $n$ (see Section~5 of Chapter ``coding-constructions'',
\index{Pless symmetry code}
\index{code!Pless symmetry}
\index{code!symmetry}
\cite{Blake74}, \cite{Blake75}, \cite{Me51}, \cite{Ple9}, \cite{Ple11}).
A $[28,14,9]_3$ code was discovered by Cheng and R. Scharlau \cite{ChSch87}.
Another such code was given by Kschischang and Pasupathy \cite{KP},
namely the negacyclic code generated by the polynomial
$(x^2 + x-1) (x^6 -x^4 + x^3 + x^2 -1) (x^6 - x^5 -x-1)$, i.e. by the vectors
$$(2002021222020010000000000000)_- ~,$$
where the subscript $-$ indicates that the code is negacyclic.
Huffman \cite{huff3} shows that there are at least 14 inequivalent
$[28,14,9]_3$ codes with nontrivial automorphisms of odd order.

Ward \cite{Wa85} and Dawson \cite{Dawso84} independently discovered that $[40,20,12]_3$ codes can be constructed using generator matrices of the form
$[I_{20} H_{20} ]$, where $H_{20}$ is a Hadamard matrix of order 20.
There are three distinct Hadamard matrices of this order,
and Dawson shows that all three produce $[40,20,12]_3$ codes.
Harada \cite{Hara97b} shows that these three codes are inequivalent.
Dawson also shows that the same construction using the Paley-Hadamard matrix of order 32 leads to a $[64,32,18]_3$ self-dual code.
A $[64,32,18]_3$ code $B_{24}$ (equivalent to Dawson's) had been
constructed earlier by Beenker \cite{Bee1}.

The codes of length 32, 44, 52, 56 and 68
can be obtained by ``subtracting'' (see Section~\ref{GU3}) a copy of $t_4$ from a
code of length 4 greater.

Other constructions for ternary self-dual codes can be found in Harada \cite{Hara97b} and Ozeki \cite{Oze89b}.
\subsection{Family \Ethree: Hermitian self-dual codes over $\FF_4$}\label{EX3}
\hsp
Table~\ref{TX4a} shows the highest possible
minimal distance for Hermitian self-dual
codes over $\FF_4$ of lengths $n \le 32$.
\begin{table}[htb]
\caption{Highest minimal distance of Hermitian self-dual codes over $\FF_4$}
\begin{center}
\begin{tabular}{rrl} \hline
$n$ & \multicolumn{1}{c}{$d$} & Codes \\ \hline
2 & 2 & $i_2$. \\
4 & 2 & $i_2^2$. \\
6 & 4 & $h_6$. \\
8 & 4 & $e_8$. \\
10 & 4 & $d_{10}^+ , e_5^{2+}$. \\
12 & 4 & 5 codes (Table~\ref{T4a}). \\
14 & 6 & $q_{14}$. \\
16 & 6 & 4 codes \cite{Me65}. \\
18 & 8 & $S_{18}$ \cite{huff4}. \\
20 & 8 & 2 codes \cite{huff4}. \\
22 & 8 & $\ge 38$ codes \cite{huff1}, \cite{huff2} \\
24 & 8 & $\ge 1$ code \\
26 & \multicolumn{1}{c}{8 or 10} & ? \\
28 & 10 & $\ge 3$ codes \cite{huff1}, \cite{huff2} \\
30 & 12 & $XQ_{29}$ \cite{Me55} \\
32 & ? & ?
\end{tabular}
\end{center}
\label{TX4a}
\end{table}
\subsection*{Remarks on Table~\ref{TX4a}}
\hsp
A period in the ``Codes'' column indicates that the list is complete.

For the entries at lengths $n \le 16$, see the discussion in Section~\ref{GU5}.

Extremal codes exist at lengths 2, 4, 6, 8, 10, 14, 16, 18, 20, 22, 28 and 30.
They do not exist at lengths 12, 24, 102, 108, 114, 120, 122 and all $n \ge 126$ (the larger $n$
being eliminated by the presence of negative coefficients in the extremal Hamming
weight enumerator).
The remaining lengths $(26, 32, 34, \ldots)$ are undecided.

The $[18,9,8]_4$ code $S_{18}$ generated by
$$1(1 \om \oom \om \om \om \oom \oom \oom \oom \oom \oom \om \om \om \oom \om )$$
has a number of interesting properties
(see \cite{Me55}, \cite{Me146}, \cite{Me147},
\cite{Ple14a}).
It has automorphism group $3 \times (PSL_2 (16).4)$,
of order 48960 \cite{Me147} and is the unique $[18,9,8]_4$ code \cite{huff4}.

The long-standing question of the existence of a $[24,12,10]_4$ code was
settled in the negative by Lam and Pless \cite{lam} (see also
\cite{huff.5}).  The code $g_{24}\otimes\FF_4$ is an example of a
$[24,12,8]_4$ code.
\subsection{Family \Efive: Additive self-dual codes over $\FF_4$}\label{EX4}
\hsp
Table~\ref{TX5a}, taken from \cite{Me223}, shows the highest
possible minimal distance for additive codes over $\FF_4$
of lengths $n \le 30$
that are self-dual with respect to the trace inner product.
\begin{table}[htb]
\caption{Highest minimal distance of additive self-dual codes over $\FF_4$}
\begin{center}
\begin{tabular}{rcl@{~~~~~~~~~~}rcl} \hline
$n$ & $d$ & Codes & $n$ & $d$ & Codes \\ \hline
1 & 1 & $i_1$. & 16 & 6 & $\ge 4$ codes \cite{Me65} \\
2 & 2 & $i_2$. & 17 & 7 \\
3 & 2 & $d_3^+$. & 18 & 8 & $S_{18}$ \\
4 & 2 & 3 codes. & 19 & 7 \\
5 & 3 & $h_5$. & 20 & 8 & $\ge 2$ codes \cite{huff4} \\
6 & 4 & $h_6$. & 21 & 8 & $c_{21}$ \\
7 & 3 & ~ & 22 & 8 & $\ge 38$ codes \cite{huff1} \\
8 & 4 & $e_8$ & 23 & 8--9 & $c_{23}$ \\
9 & 4 & $c_9$ & 24 & 8--10 & $g_{24}\otimes\FF_4$ \\
10 & 4 & $d_{10}^+$ , $e_5^{2+}$ & 25 & 8--9 & $c_{25}$ \\
11 & 5 & ~ & 26 & 8--10 \\
12 & 6 & $z_{12}$. & 27 & 9--10 \\
13 & 5 & ~ & 28 & 10 \\
14 & 6 & $q_{14}$ & 29 & 11 \\
15 & 6 & $c_{15}$ & 30 & 12 & $XQ_{29}$
\end{tabular}
\end{center}
\label{TX5a}
\end{table}
\subsection*{Remarks on Table~\ref{TX5a}}
\hsp
A period in the ``Codes'' column indicates that the list is complete.

Extremal Type I codes exist at lengths 1--6, 8--12, 14--18,
20--22 and 28--30, and do not exist at lengths 7, 13 and 25.
Lengths 19, 23, 24, 26, 27 are undecided.

Many of the entries are
copied from the table of Hermitian self-dual
codes, Table~\ref{TX4a}.
The codes $d_n^+$ are defined in Section~\ref{GU7},
$h_6$ is the hexacode, and $h_5$ is the $[5,2.5,3]_{4+}$ shortened
hexacode, generated by $(01 \om \om 1)$, with weight enumerator
$x^5 + 10 x^2 y^3 + 15 xy^4+ 6y^5$ and $|Aut (h_5)| =120$.
Also, $c_9$, $c_{15}$,
$c_{21}$, $c_{23}$, $c_{25}$ are cyclic codes
with generators shown in Table~\ref{TX5b}.
If no name is given, the code can be obtained by shortening a code of length
one greater.
\begin{table}[htb]
\caption{Generators for cyclic additive codes over $\FF_4$}
$$
\begin{array}{ll}
c_{9} & (\om 10100101) \\
c_{15} & ( \om 11010100101011) \\
c_{21} & ( \oom \oom 1 \om 00111101011011000),
(101110010111001011100) \\
c_{23} & ( \om 0101111000000001111010) \\
c_{25} & (111010 \om 010111000000000000)
\end{array}
$$
\label{TX5b}
\end{table}
\subsection{Family \Eeight: Self-dual codes over $\ZZ_4$}\label{EXZ4}
\hsp
Table~\ref{TXZ4} gives the highest possible Hamming distance,
Lee distance and Euclidean norm for codes over $\ZZ_4$
of lengths $n \le 24$.  This is based on
\cite{Me168}, \cite{DHS97}, \cite{FGLP}, \cite{huff5}, \cite{PLF96}
and \cite{rainsZ4}.
The columns headed \# give the number of extremal codes.
\begin{table}[htb]
\caption{Highest Hamming distance $(d_H)$, Lee distance $(d_L)$ and Euclidean norm (Norm) of self-dual codes over $\ZZ_4$}
$$
\begin{array}{c|ccc|ccc|ccc} \hline
\mbox{Length} & \multicolumn{3}{|c}{\mbox{Hamming}} & \multicolumn{3}{|c}{\mbox{Lee}} & \multicolumn{3}{|c}{\mbox{Norm}} \\
n & d_H & \mbox{code} & \# & d_L & \mbox{code} & \# & \mbox{Norm} & \mbox{code} & \# \\ \hline
1 & 1 & i_1 & 1 & 2 & i_1 & 1 & 4 & i_1 & 1 \\
2 & 1 & i_1^2 & 1 & 2 & i_1^2 & 1 & 4 & i_1^2 & 1 \\
3 & 1 & i_1^3 & 1 & 2 & i_1^3 & 1 & 4 & i_1^3 & 1 \\
4 & 2 & D_4^\oplus & 1 & 4 & D_4^\oplus & 1 & 4 & i_1^4 & 2 \\
5 & 1 & D_4^\oplus i_1 & 2 & 2 & D_4^\oplus i_1 & 2 & 4 & i_1^5 & 2 \\
6 & 2 & D_6^\oplus & 1 & 4 & D_6^\oplus & 1 & 4 & i_1^6 & 3 \\
7 & 3 & E_7^+ & 1 & 4 & E_7^+ & 1 & 4 & i_1^7 & 4 \\
8 & 4 & o_8 & 2 & 6 & o_8 & 1 & 8 & o_8 & 1 \\
9 & 1 & o_8 i_1& 11 & 2 & o_8 i_1 & 11& 4 & i_1^9 & 11 \\
10 & 2 & D_4^\oplus D_6^\oplus & 5 & 4 & D_4^\oplus D_6^\oplus & 5 & 4 & i_1^{10} & 16 \\
11 & 2 & D_4^\oplus E_7^+ & 3 & 4 & D_4^\oplus E_7^+ & 3 & 4 & i_1^{11} & 19 \\
12 & 2 & D_4^\oplus o_8 & 39 & 4 & D_4^\oplus o_8 & 39 & 8 & \mbox{\cite{FGLP}} & 19 \\
13 & 2 & D_6^\oplus E_7^+ & 8 & 4 & D_6^\oplus E_7^+ & 8 & 4 & i_1^{13}  & 66\\
14 & 3 & (E_7^+)^2 & 4 & 6 & \mbox{\cite{FGLP}} & 1 & 8 & \mbox{\cite{FGLP}} & 35 \\
15 & 3 & E_7^+ o_8 & 47 & 6 & \mbox{\cite{FGLP}} & 15 & 8 & \mbox{\cite{FGLP}} & 28 \\
16 & 4 & o_8^2 & \ge 1 & 8 & C_{16} & \ge 5 & 8 & o_8^2 & \ge 5 \\
17 & 4 & C_{17} & 62 & 6 & C_{17} & \ge 17 & 8 & C_{17}&\ge 17 \\
18 & 4 & C_{18} & 66 &  8 & C_{18} & 7 &8 & C_{18}& \ge 39 \\
19 & 3 & G_{19} & \ge 1 & 6 & G_{19} & \ge 1 & 8 & G_{19}& \ge 1 \\
20 & 4 & G_{20} & \ge 1 & 8 & G_{20} & \ge 1 & 8 & G_{20}& \ge 1 \\
21 & 5 & G_{21} & 384 & 8 & G_{21} & 384 & 8 & G_{21} & \ge 384 \\
22 & 6 & G_{22} & \ge 19367 & 8 & G_{22} & \ge 19367 & 8 & G_{22} & \ge 19367 \\
23 & 7 & G_{23} & \ge 1.72 \times 10^6 & 10 & G_{23} & 30 & 12 & G_{23} & \ge 30 \\
24 & 8 & G_{24} & \ge 1.47 \times 10^8 & 12 & G_{24} & 13 & 16 & G_{24} & \ge 50 
\end{array}
$$
\label{TXZ4}
\end{table}
\subsection*{Remarks on Table~\ref{TXZ4}}
\hsp
The length 16 code $C_{16}$ is given in \cite{PLF96}, where it is called 5\_f5.
It has $|Aut (C_{16}) |= 2^{5+10} 3^2 5.7$
and generator matrix
$$\left[
\begin{array}{ccccccccccc|ccccc}
1 & 1 & 1 & 1 & 1 & 1 & 1 & 1 & 1 & 1 & 1 & 1 & 1 & 1 & 1 & 1 \\
1 & 0 & 1 & 1 & 1 & 1 & 1 & 1 & 0 & 0 & 0 & 0 & 1 & 0 & 0 & 0 \\
1 & 1 & 0 & 1 & 0 & 0 & 1 & 1 & 1 & 1 & 0 & 0 & 0 & 1 & 0 & 0 \\
1 & 1 & 1 & 0 & 1 & 0 & 1 & 0 & 1 & 0 & 1 & 0 & 0 & 0 & 1 & 0 \\
0 & 0 & 0 & 0 & 1 & 1 & 1 & 1 & 1 & 1 & 1 & 0 & 0 & 0 & 0 & 1 \\ \hline
0 & 0 & 0 & 0 & 0 & 2 & 0 & 0 & 0 & 0 & 0 & 2 & 2 & 0 & 0 & 2 \\
0 & 0 & 0 & 0 & 0 & 0 & 2 & 0 & 0 & 0 & 0 & 2 & 2 & 2 & 2 & 2 \\
0 & 0 & 0 & 0 & 0 & 0 & 0 & 2 & 0 & 0 & 0 & 0 & 2 & 2 & 0 & 2 \\
0 & 0 & 0 & 0 & 0 & 0 & 0 & 0 & 2 & 0 & 0 & 0 & 0 & 2 & 2 & 2 \\
0 & 0 & 0 & 0 & 0 & 0 & 0 & 0 & 0 & 2 & 0 & 2 & 0 & 2 & 0 & 2 \\
0 & 0 & 0 & 0 & 0 & 0 & 0 & 0 & 0 & 0 & 2 & 2 & 0 & 0 & 2 & 2
\end{array}
\right]
$$

The codes $C_{17}$ and $C_{18}$ mentioned in the table have generator
matrices
$$
\left[
\begin{array}{ccccccccccccccccc}
1 & 0 & 0 & 0 & 0 & 0 & 0 & 0 & 1 & 2 & 1 & 1 & 1 & 1 & 1 & 1 & 2 \\
0 & 1 & 0 & 0 & 0 & 0 & 0 & 0 & 1 & 1 & 0 & 0 & 0 & 1 & 2 & 2 & 0 \\
0 & 0 & 1 & 0 & 0 & 0 & 0 & 0 & 1 & 3 & 2 & 0 & 0 & 2 & 3 & 2 & 2 \\
0 & 0 & 0 & 1 & 0 & 0 & 0 & 0 & 1 & 3 & 0 & 0 & 0 & 2 & 0 & 3 & 0 \\
0 & 0 & 0 & 0 & 1 & 0 & 0 & 0 & 0 & 1 & 3 & 3 & 3 & 3 & 1 & 1 & 2 \\
0 & 0 & 0 & 0 & 0 & 1 & 0 & 0 & 0 & 0 & 0 & 3 & 1 & 0 & 2 & 0 & 3 \\
0 & 0 & 0 & 0 & 0 & 0 & 1 & 0 & 0 & 0 & 3 & 2 & 1 & 0 & 0 & 0 & 3 \\
0 & 0 & 0 & 0 & 0 & 0 & 0 & 1 & 0 & 0 & 1 & 3 & 2 & 0 & 0 & 0 & 3 \\
0 & 0 & 0 & 0 & 0 & 0 & 0 & 0 & 2 & 2 & 2 & 2 & 2 & 0 & 0 & 0 & 0 \\
\end{array}
\right]
$$
and
$$
\left[
\begin{array}{cccccccccccccccccc}
1 & 0 & 0 & 3 & 0 & 0 & 0 & 0 & 0 & 3 & 2 & 0 & 3 & 0 & 2 & 0 & 0 & 0 \\
0 & 1 & 0 & 0 & 3 & 0 & 0 & 0 & 0 & 0 & 3 & 2 & 2 & 3 & 0 & 0 & 0 & 0 \\
0 & 0 & 1 & 0 & 0 & 3 & 0 & 0 & 0 & 2 & 0 & 3 & 0 & 2 & 3 & 0 & 0 & 0 \\
0 & 0 & 0 & 1 & 0 & 0 & 3 & 0 & 0 & 0 & 0 & 0 & 3 & 2 & 0 & 3 & 0 & 2 \\
0 & 0 & 0 & 0 & 1 & 0 & 0 & 3 & 0 & 0 & 0 & 0 & 0 & 3 & 2 & 2 & 3 & 0 \\
0 & 0 & 0 & 0 & 0 & 1 & 0 & 0 & 3 & 0 & 0 & 0 & 2 & 0 & 3 & 0 & 2 & 3 \\
0 & 0 & 0 & 0 & 0 & 0 & 1 & 3 & 0 & 0 & 3 & 3 & 0 & 3 & 3 & 3 & 2 & 1 \\
0 & 0 & 0 & 0 & 0 & 0 & 0 & 1 & 3 & 3 & 0 & 3 & 3 & 0 & 3 & 1 & 3 & 2 \\
0 & 0 & 0 & 0 & 0 & 0 & 0 & 0 & 2 & 0 & 2 & 2 & 0 & 2 & 2 & 0 & 2 & 0 \\
0 & 0 & 0 & 0 & 0 & 0 & 0 & 0 & 0 & 2 & 2 & 2 & 2 & 2 & 2 & 2 & 2 & 2 \\
\end{array} ~,
\right]
$$
and automorphism groups of orders 576 and 144, respectively.

$G_{24}$ was defined in \eqn{Eq12ab}, and $G_{19}$ through $G_{23}$ are shortened versions of it.

Besides the norm-extremal codes of length 8, 12, 14--24 shown in the table,
there are also norm-extremal codes of lengths 32 and 48 obtained by lifting
binary extended quadratic residue codes to $\ZZ_4$.  The code of length 32
has minimal Lee weight 14 and minimal norm 16.  Pless and Qian \cite{PlQ97}
have shown that the code of length 48 has minimal Lee weight 18 and minimal
norm 24.

Further examples of good self-dual codes over $\ZZ_4$ may be found in
\cite{BCS95}, \cite{Betal98}, \cite{BSBM97}, \cite{Me209},
\cite{DHS97}, \cite{GHa98}, \cite{GuHa97}, \cite{GuHa97c},
\cite{Hara97c},
\cite{huff5},
\cite{PlQ97}, \cite{PSQ97}, \cite{rainsZ4}.
\section{Further topics}\label{FT}
\subsection{Decoding self-dual codes}\label{FTD}
\hsp
The problem of decoding\index{decoding!self-dual codes} self-dual codes is an extremely important one for applications, but we will not discuss it here.
Decoding the binary Golay code, in particular, has been studied in many
papers --- see
\cite{Amr94}, \cite{Me116}, \cite{Me146}, \cite[Chapter~11]{SPLAG},
\cite{Ple15},
\cite{Ple15a}, \cite{RaS95}, \cite{SnB1}, \cite{Vardy95}, \cite{VB1}.
See also \cite{RaS93}, \cite{Vardy94}, and Section~8 of Chapter ``codes-and-groups''.
\subsection{Applications to projective planes}\label{PP}\index{projective plane}
\hsp
There is a very nice application of self-dual codes to projective planes.
If $n$ is congruent to 2 $( \bmod~4)$ then the incidence matrix of a projective
plane of order $n$ generates a self-orthogonal code $C_n$, which when an 
overall
parity-check is added becomes an $[n^2 +n+2$, 
$\frac{1}{2} (n^2 + n+2) , n+2 ]$ Type II self-dual binary code (see
\cite{AsKe92}, \cite{Me13} or Chapter ``assmus'' for the proof).

It was a famous unsolved problem to decide if a projective plane of
order 10 could exist.
The weight enumerator of $C_{10}$ was initially studied in \cite{Me13} (see also \cite{Me35}).
Finally, after many years of work, Lam, Thiel and Swiercz \cite{lam89}
(see also \cite{lam91}) succeeded in completing this project and showed that
$C_{10}$ (and hence the putative plane of order 10) does not exist.

The possibility of the existence of a plane of order 18 (or 12, but then
we do not obtain a self-dual code) remains an open question.
\subsection{Automorphism groups of self-dual codes}\label{AGSD}
\hsp
Various topics concerning the automorphism groups\index{automorphism group} of self-dual codes are discussed in chapter ``codes-and-groups'',
\index{group!trivial}
e.g. the full automorphism groups of extended quadratic residue codes, the occurrence of self-dual codes with a trivial group\index{trivial group} (see
\cite{Bus3}, \cite{Me158}, \cite{Hara96}, \cite{Me76}, \cite{OrPh}, \cite{Ton6}),
and the existence of self-dual codes with any prescribed symmetry group (\cite{OrPh}).

\subsection{Open problems}\label{OP}\index{open problems}
\hsp
Do there exist $[72,36,16]$ or $[96,48,20]$ Type II self-dual binary codes?
(Cf. \cite{cp2}, \cite{Feit}, \cite{hy}, \cite{pt}, \cite{Me31}).

Fill in the other gaps in Tables~\ref{TX2a}, \ref{TX3a}, \ref{TX4a}.
No extremal Hermitian self-dual codes over $\FF_4$ of any length greater than 30 are presently known!

There is an interesting open question concerning self-dual codes of length 24.
There exists a unique $[24,12,8]$ binary code, exactly two $[24,12,9]_3$ ternary codes, and no $[24,12,10]_4$ Hermitian or Euclidean self-dual code over $\FF_4$ (\cite{lam}).
But the possibility of an {\em additive} trace-self-dual code of length 24 over $\FF_4$ with minimal distance 10 remains open (see
Table~\ref{TX5a}).
From Theorem~\ref{thBD7}, if such a code exists then it must be even.
However, all our attempts so far to construct this code have failed,
so it may not exist.

When is the first time a Type I binary code has a higher minimal distance than the best Type II code of the same length?
(No such example is presently known.)

In this regard it is worth mentioning that there is a $[32,17,8]$ binary code \cite{Me148}, which has the same minimal distance as the best self-dual codes of length 32,
yet contains twice as many codewords.
There are similar examples in the ternary case --- see Chapter ``Brouwer''.

The Nordstrom-Robinson code (see Chapter~1) is an example of a nonlinear code that has a higher minimal distance than any self-dual (or even linear) code of the same length.
However, as mentioned in Section~\ref{Wee2}, the Nordstrom-Robinson code should really be regarded as a self-dual linear code over $\ZZ_4$ (the octacode $o_8$).
When is the first time a non-self-dual $[n,n/2, d]$ binary linear code has a higher minimal
distance than any $[n,n/2,d']$ self-dual code?
This certainly happens at length 40, but may happen at length 36 or 38.

Is there any difference asymptotically, as $n \to \infty$, between
$d/n$ for the best binary codes, the best binary linear codes and the best binary self-dual codes?

Let $\Om_n$ denote the collection of binary self-dual codes that have the highest possible minimal distance at length $n$, and let $L_n$, $U_n$ be respectively the
smallest and largest orders of $Aut (C)$, $C \in \Om_n$.
When (if ever) is the first time that $L_n = U_n =1$?
Is there an infinite sequence of values of $n$ with $U_n > 1$?
Show that $L_n =1$ for all sufficiently large $n$.
\section{Self-dual codes and lattices}\label{Latt}\index{lattice}
\hsp
There are many connections and parallels between self-dual codes and lattice sphere packings.
Our original intention was to end the chapter with an account of these connections, but constraints of space and time have not permitted this.
Instead, we give a brief list of some of the parallels, to whet the
reader's appetite.
For more information about the relationship between the two fields, see
\cite{Bro74}, \cite{MS201}, \cite{Ebe94},
\cite{Forn1}, \cite{Forn2}, \cite{GHS97},
\cite{Me58}, \cite{Me64} and especially \cite{SPLAG}, \cite{Me2000}.
\index{lattice!Leech}\index{lattice!$E_8$}
\smallskip
\begin{center}
\begin{tabular}{ll} \hline
Coding concept & Lattice concept \\ \hline
Binary linear code & Lattice \\
Dual code & Dual lattice \\
Self-orthogonal code & Integral lattice \\
Self-dual code & Unimodular lattice \\
Doubly-even self-dual code & Even unimodular lattice \\
Hamming code $e_8$ & Root lattice $E_8$\index{root lattice $E_8$} (\cite{SPLAG}, p.~120) \\
Hexacode $h_6$ & Coxeter-Todd lattices $K_{12}$ (\cite{SPLAG}, p.~127) \\
Binary Golay code $g_{24}$ & Leech lattice\index{Leech lattice} $\La_{24}$ (\cite{SPLAG}, p.~131) \\
Minimal distance & Minimal norm \\
Number of minimal weight words & Kissing number \\
Weight enumerator $W(x,y)$ & Theta series\index{theta series} \\
MacWilliams identity (Eq. (33)) & Jacobi identity\index{Jacobi identity} (\cite{SPLAG}, p.~103) \\
~~~(weight enumerator of dual code in terms & ~~~(theta series of dual lattice in terms \\
~~~of weight enumerator of code) & ~~~of theta series of lattice) \\
Gleason's theorem (Theorem ~11) & Hecke's theorem\index{Hecke theorem}\index{theorem!Hecke} (\cite{SPLAG}, p.~187) \\
~~~(weight enumerator of doubly-even & ~~~(theta series of even unimodular \\
~~~code is polynomial in weight enumerators & ~~~lattice is polynomial in theta series \\
~~~of $e_8$ and $g_{24}$) & ~~~of $E_8$ and $\La_{24}$) \\ \hline
\end{tabular}
\end{center}
\smallskip

The similarity between the theorems of Gleason and Hecke is particularly
striking, and we will end the chapter by saying a little more about this.
Suppose $C$ is a binary code of length $n$.
{\em Construction A}\index{Construction A} produces an $n$-dimensional sphere packing $\La (C)$, consisting of the points $\frac{1}{\sqrt{2}}x$ for $x \in \ZZ^n$, $x~( \bmod~2) \in C$.
If $C$ is linear, $\La (C)$ is a lattice;
if $C$ is self-dual, $\La (C)$ is unimodular; and if $C$ is Type II, $\La (C)$ is an even unimodular lattice.

If $C$ is a linear code with weight enumerator $W_C (x,y)$, then $W_C (\theta_3(2z), \theta_2 (2z))$ is the theta series of $\La (C)$, where
$$\theta_3 (z) = \sum_{m=- \infty}^\infty q^{m^2} , \quad
\theta_2 (z) = \sum_{m=- \infty}^\infty q^{(m+1/2)^2} ~,
$$
where $q = e^{\pi iz}$, ${\rm Im} (z) > 0$.
This map gives an isomorphism between (a) the ring of weight enumerators
of Type I self-dual codes,
$\CC [\phi_2, \theta_8]$ (see Eq.~\ref{EqG6a}),
and the ring of theta series of even-dimensional unimodular
lattices,
$\CC[\theta_3^2, \Delta_8]$, where
$$\Delta_8 = q \prod_{m=1}^\infty \{ (1-q^{2m-1}) (1-q^{4m}) \}^8 ~;$$
and (b) the ring of weight enumerators of Type II self-dual codes,
$\CC[\phi_8 , \phi'_{24} ]$ (Theorem~\ref{ThMS3c}), and the ring of theta series
of even unimodular lattices,
$\CC [\Theta_{E_8} , \Delta_{24} ]$, where
\begin{eqnarray*}
\Theta_{E_8} (z) & = & 1+ 240 \sum_{m=1}^\infty \sigma_3 (m) q^{2m} ~, \\
\Delta_{24} & = & q^2 \prod_{m=1}^\infty (1- q^{2m} )^{24} ~,
\end{eqnarray*}
and $\sigma_3 (m)$ is the sum of the cubes of the divisors of $m$.
For further information see \cite[Chapter~7]{SPLAG}.

The bibliography also contains a number of references that are concerned with particular
constructions of lattices from self-dual codes, or of properties of lattices
that are analogous to properties of self-dual codes mentioned in this chapter:
\cite{Bach97}, \cite{BDHO}, \cite{BCS95}, \cite{Betal98},
\cite{BoSo94}, \cite{Me209}, 
\cite{ChSo96}, \cite{Me59}, \cite{Me157}, \cite{Me177}, \cite{DGH97}, \cite{GuHa97},
\cite{KiKM},
\cite{Koch86},
\cite{KV89},
\cite{KP},
\cite{Me41}, \cite{Oze1}, \cite{Oze2},
\cite{Oze4},
\cite{Oze87}, \cite{Oze89a}, 
\cite{Oze91},
\cite{RaSl97}, \cite{Rus89}, \cite{Rus93}, \cite{Me136}, \cite{Me57}.
\addcontentsline{toc}{section}{Acknowledgements}
\section*{Acknowledgements}
\hsp
Over the past 26 years NJAS has had the pleasure of collaborating with many of the people whose names are listed in the
bibliography:
he wishes to express his appreciation to all of them.
We thank Eiichi Bannai, Dave Forney, Aaron Gulliver, Masaaki Harada,
Cary Huffman, Michio Ozeki, Vera Pless and Patrick Sol\'{e} for helpful
comments on the manuscript of this chapter.
We also thank Susan Pope for a superb typing job.

\clearpage
\addcontentsline{toc}{section}{Bibliography}

The bibliography uses the
following abbreviations for journals:\\

\begin{tabular}{lll}
DCC & $=$ & {\em Designs, Codes and Cryptography} \\
~~ \\
DM & $=$ & {\em Discrete Mathematics} \\
~~ \\
JCT & $=$ & {\em Journal of Combinatorial Theory} \\
~~ \\
PGIT & $=$ & {\em IEEE Transactions on Information Theory}
\end{tabular}

\vspace*{+.1in}

\end{document}

%% file: box.pstex_t
\begin{picture}(0,0)%
\epsfig{file=box.pstex}%
\end{picture}%
\setlength{\unitlength}{0.00066700in}%
\begingroup\makeatletter\ifx\SetFigFont\undefined%
\gdef\SetFigFont#1#2#3#4#5{%
  \reset@font\fontsize{#1}{#2pt}%
  \fontfamily{#3}\fontseries{#4}\fontshape{#5}%
  \selectfont}%
\fi\endgroup%
\begin{picture}(3924,4224)(2689,-7273)
\put(4651,-6886){\makebox(0,0)[b]{\smash{\SetFigFont{10}{12.0}{\rmdefault}{\mddefault}{\updefault}$X$}}}
\put(6001,-3436){\makebox(0,0)[b]{\smash{\SetFigFont{10}{12.0}{\rmdefault}{\mddefault}{\updefault}$0$}}}
\put(3151,-5836){\makebox(0,0)[b]{\smash{\SetFigFont{10}{12.0}{\rmdefault}{\mddefault}{\updefault}$0$}}}
\put(6001,-5986){\makebox(0,0)[b]{\smash{\SetFigFont{10}{12.0}{\rmdefault}{\mddefault}{\updefault}$G_t$}}}
\put(4051,-4186){\makebox(0,0)[b]{\smash{\SetFigFont{10}{12.0}{\rmdefault}{\mddefault}{\updefault}$G_2$}}}
\put(3151,-3436){\makebox(0,0)[b]{\smash{\SetFigFont{10}{12.0}{\rmdefault}{\mddefault}{\updefault}$G_1$}}}
\end{picture}